\documentclass[11pt]{article}
\usepackage{graphicx,mathrsfs,amssymb,color}
\usepackage{amsmath,latexsym,amsbsy}
\def\ds{\displaystyle}
\def\forall{\hbox{for all}~}
\def\L{{\bf L}}

\def\ve{\varepsilon}
\def\n{\noindent}

\def\dint{\int\!\!\int}

\def\R{\mathbb{R}}

\def\vp{\varphi}

\def\vs{\vskip 2em}

\def\v{\vskip 1em}
\def\O{{\cal O}}

\def\bfH{{\bf H}}
\def\C{{\cal C}}
\def\bega{\begin{array}}
\def\enda{\end{array}}
\def\begi{\begin{itemize}}
\def\endi{\end{itemize}}
\def\ov{\overline}
\def\Tilde{\widetilde}

\def\bel{\begin{equation}\label}
\def\eeq{\end{equation}}
\def\sqr#1#2{\vbox{\hrule height .#2pt
\hbox{\vrule width .#2pt height #1pt \kern #1pt
\vrule width .#2pt}\hrule height .#2pt }}
\def\square{\sqr74}
\def\endproof{\hphantom{MM}\hfill\llap{$\square$}\goodbreak}

\newtheorem{theorem}{Theorem}[section]

\newtheorem{lemma}{Lemma}[section]

\newtheorem{remark}{Remark}[section]
\newtheorem{definition}{Definition}[section]

\usepackage{marginnote}

\begin{document}
\title{\bf Shock interactions for the Burgers-Hilbert Equation}\vs

\author{Alberto Bressan$^{(1)}$, Sondre T.~Galtung$^{(2)}$,   Katrin Grunert$^{(2)}$,  and 
Khai T.~Nguyen$^{(3)}$\\~~\\
 {\small $^{(1)}$~Department of Mathematics, Penn State University,}\\
 {\small $^{(2)}$~Department of Mathematical Sciences, NTNU Norwegian University of Science and Technology, Trondheim,}\\
 {\small $^{(3)}$~Department of Mathematics, North Carolina State University.} \\~~\\
{\small E-mails: bressan@math.psu.edu,~sondre.galtung@ntnu.no,~katrin.grunert@ntnu.no,~khai@math.ncsu.edu.}
}
\maketitle

\begin{abstract} This paper provides an asymptotic description of a solution to the Burgers-Hilbert equation
in a neighborhood of a point where two shocks interact. The solution is obtained as the sum of 
a function with  $H^2$ regularity away from the shocks plus a corrector term 
having an asymptotic behavior like $|x|  \ln |x| $ close to each shock.  
A key step in the analysis is the construction of
piecewise smooth solutions with a single shock for a general class of  initial data.\end{abstract}

\section{Introduction}
\label{sec:1}
\setcounter{equation}{0}
Consider the balance law obtained from Burgers' equation by adding
the Hilbert transform as
a source term
\bel{BH} 
u_t + \left({u^2\over 2}\right)_x~=~\bfH [u]\,.\eeq
This equation was derived in \cite{BiH} as a model for
nonlinear waves with constant frequency. Here the nonlocal source term
$$\bfH[f](x)~\doteq~\lim_{\ve\to 0+} {1\over\pi}
\int_{|y|>\ve} {f(x-y)\over y} \, dy
$$
denotes the Hilbert transform of a function $f\in \L^2(\R)$. 
It is well known \cite{S} that $\bfH$ is a linear isometry from $\L^2(\R)$
 onto itself. Given  any initial data 
\bel{IC}
u(0,\cdot)~=~\bar{u}(\cdot)
\eeq
with $\bar u\in H^2(\R)$, 
the local existence and uniqueness of  solutions to (\ref{BH}) was proved in \cite{HI}, together with a sharp estimate on the time interval where this solution remains smooth. For a general initial data $\bar u\in {\bf L}^2(\R)$, the global existence of entropy weak  solutions to 
(\ref{BH}) was proved in \cite{BN}, together with a partial uniqueness result. 
We remark that the well-posedness of the Cauchy problem for (\ref{BH}) 
remains a largely open question. 

\begin{figure}[htbp]
\centering
  \includegraphics[scale=0.42]{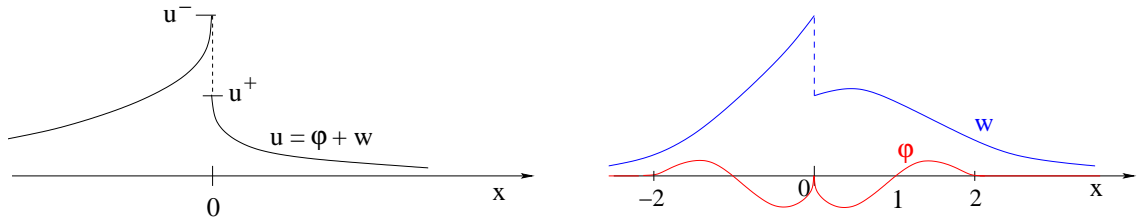}
    \caption{\small  Decomposing a solution in the form (\ref{uwp})}
\label{f:bh17}
\end{figure}

More recently, piecewise continuous solutions with a single shock  
have been
constructed in \cite{BZ}. As shown in Fig.~\ref{f:bh17}, these solutions have the form 
\bel{uwp}
u(t,x)~=~w\bigl(t,x-y(t)\bigr)+\varphi\bigl(x-y(t)\bigr),
\eeq
where $y(t)$ denotes  the location of the shock at time $t$, and $w(t,\cdot)\in 
H^2\bigl(]-\infty,0[\,\cup \,]0,+\infty[\bigr)$ for all $t\geq 0$.
Moreover,  $\varphi$ is a fixed function with compact support, describing the asymptotic behavior of the  solution near the shock. It is smooth outside the origin and satisfies
\bel{vpdef}\vp(x)~=~{2\over \pi} |x|\ln |x|\qquad\quad \hbox{for}\quad 
|x|\leq1\,.\eeq  
Remarkably, this ``corrector term" $\vp$ is universal, i.e., it does not depend on the particular
solution of (\ref{BH}).   The same analysis applies to
 solutions with finitely many, noninteracting shocks.  
In addition, the local asymptotic behavior of a solution up to the time when a new shock is formed
was investigated in \cite{Y}.    
  
The aim
of the present note is to describe the asymptotic behavior of a solution in a neighborhood of a point
where two shocks interact.  
Calling $T>0$ the time when the 
interaction takes place, our analysis splits into two parts. 
We first describe the behavior of the solution as $t\to T-$, i.e.~as the two shocks approach each other. 
In a second step, to construct the solution for $t>T$, we 
solve a Cauchy problem with initial data given at $t=T$.  

 As it turns out, the 
profile $u(T,\cdot)$ is  not ``well prepared", in the sense that it cannot be written 
in the form (\ref{uwp}).  
To explain the difficulty, we recall that the solutions constructed in \cite{BZ} had initial data
of the form
\bel{pcd1}
u(0,x)~=~\ov w(x-y_0)+\varphi(x-y_0),\eeq
for some $\ov w\in H^2\bigl(\R\setminus \{0\}\bigr)$ and $y_0\in \R$.
These data are ``well prepared", in the sense that they already contain the corrector term $\vp$.
A natural class of initial data, not considered in \cite{BZ}, is 
\bel{pcd2}
u(0,x)~=~\ov w(x-y_0)\qquad\mathrm{with}\qquad \ov w\in H^2\bigl(\R\setminus \{0\}\bigr), \phantom{a} y_0\in\R.
 \eeq
By assumption, at time $t=0$  the derivative
$u_x(0,x)~=~\ov w_x(x-y_0)$ is piecewise continuous and uniformly bounded.
However, in the solution to (\ref{BH}), (\ref{pcd2}), at each time $t>0$ we expect that $u_x(t,x)\to \pm\infty$ as $x\to y(t)\mp$.
For this reason, the local construction of this solution requires a careful analysis.
A more general class of initial data, containing both (\ref{pcd1}) and (\ref{pcd2}), as well as all 
profiles $u(T,\cdot)$ emerging from our shock interactions,
will be studied in Section~\ref{sec:2}.

We recall here the definition of entropy weak solutions used in \cite{BN}.
\begin{definition} {\it By an {\bf entropy weak solution} of (\ref{BH})-(\ref{IC}) 
we mean a function $u\in \L^1_{loc}([0,\infty[\,\times\R)$ with the following properties.
\begi
\item[(i)] The map $t\mapsto u(t,\cdot)$ is continuous with values in $\L^2(\R)$ and
satisfies the initial condition (\ref{IC}).
\item[(ii)]  For any $k\in \R$ and every nonnegative 
test function $\phi\in\C^1_c(]0,\infty[\,\times\R)$
one has
\bel{ei}\dint  \Big[ |u-k|\phi_t + \Big({u^2-k^2\over 2}\Big)\hbox{\rm sign}(u-k) \phi_x
+H[u(t)](x)\hbox{\rm sign}(u-k)\phi  \Big]\, dxdt~\geq~0.\eeq
\endi
}
\end{definition}
The present paper will be concerned with a more regular class of solutions, which are piecewise continuous and can be determined
by integrating along characteristics. These correspond to the ``broad solutions"
considered in \cite{Bbook, RY}. Throughout the sequel, the upper dot denotes a derivative w.r.t.~time.

\begin{definition}\label{d:12}
 {\it An entropy weak solution $u=u(t,x)$   of (\ref{BH})-(\ref{IC}), defined 
on the interval $t\in [0,T]$, will be called a {\bf piecewise regular solution} if there exist
 finitely many shock curves $y_1(t), \ldots, y_n(t)$ such that the following holds.
\begi
\item[(i)] For each $t\in [0,T]$, one has $u(t,\cdot)\in H^2\bigl( \R\setminus\{y_1(t),\ldots, y_n(t)\}\bigr)$.
\item[(ii)]  For each $i=1,\ldots,n$, the Rankine-Hugoniot conditions hold:
\bel{RH1} u_i^-(t)~\doteq~u(t, y_i(t)-)~>~u(t, y_i(t)+) ~\doteq~u_i^+(t),\eeq
\bel{RH2} \dot y_i(t)~=~ {u_i^-(t) +u_i^+(t)\over 2}\,.\eeq
\item[(iii)]  Along every characteristic curve $t\mapsto x(t)$ such that 
\bel{cc1} \dot x(t)~=~u\bigl(t, x(t)\bigr),\eeq
one has
\bel{cc2} {d\over dt} u\bigl(t, x(t)\bigr) ~=~\bfH[u](x(t)).\eeq
\endi
}
\end{definition}
In the above setting,  the Hilbert transform of the piecewise
regular function $u(t,\cdot)$ can be computed using an integration by parts:
\bel{HT}
\bfH[u(t)](x)~=~{1\over\pi}\int_{-\infty}^\infty u_x(t,y)\,\ln|x-y|\, dy 
+{1\over\pi}\sum_{i=1}^n\bigl[ u_i^+(t) - u_i^-(t)\bigr]\, \ln\bigl|x- y_i(t)\bigr|\,.\eeq

The remainder of the paper is organized as follows. 
In Section~\ref{sec:2} we state a local existence and uniqueness theorem for solutions to (\ref{BH}),
valid for a class of initial data containing one single shock, but more general than in \cite{BZ}. Towards the proof of Theorem~\ref{t:21}, 
Section~\ref{sec:3} develops various a priori estimates, while in Section~\ref{sec:4} 
the local solution is constructed as a limit of a convergent sequence of approximations.
As in \cite{BZ}, these are obtained by iteratively solving a sequence of linearized problems.

In the second part of the paper we study solutions of (\ref{BH}) with two shocks, up to the time
of interaction.   
In Section~\ref{sec:5} we perform some preliminary computations, motivating a particular form of the corrector term.  In Section~\ref{sec:6} we state and prove the second main result of the paper, Theorem~\ref{t:61}, providing a detailed description of solutions up to the interaction time.  
This is achieved by a change of both time and space coordinates, so that the two shocks are
located at the two points 
$$x_1(t) ~=~ t~ < ~0~ = ~x_2(t),$$
and interact at time $t=0$.
Our analysis shows that, at the interaction time, the solution profile contains a single shock and 
lies within the class of initial data covered by Theorem~\ref{t:21}.
Combining our two theorems, one thus obtains a complete description 
of the solution to (\ref{BH})
in a neighborhood of the interaction time.

\section{Solutions with one shock  and general initial data}
\label{sec:2}
\setcounter{equation}{0}
Consider a piecewise regular solution of the Burgers-Hilbert equation (\ref{BH}), with one single shock. 
By the Rankine-Hugoniot condition, the location $y(t)$ of the shock at time $t$ satisfies  
\bel{RH}
\dot{y}(t)~=~{u^{-}(t)+u^{+}(t)\over 2},\qquad u^{\pm}(t)~=~\lim_{x\to y(t)\pm}u(t,x).
\eeq
 As in \cite{BZ},  we shift the space coordinate, replacing $x$ with $x-y(t)$, so that in the new coordinate system
the shock is always located at the origin. In these new coordinates, (\ref{BH}) takes the  form
\bel{BH-1}
u_t+\left(u- {u^-(t)+u^+(t)\over 2}\right) u_x~=~{\bf H}[u].
\eeq
In \cite{BZ}, given a ``well prepared" initial data (\ref{pcd1}), a unique 
piecewise smooth entropy solution to (\ref{BH-1}) of the form 
\[
u(t,x)~=~w(t,x)+{2\eta(x)\over \pi}\cdot |x|\ln|x|\,,\qquad\qquad t>0\,
\]
was constructed.
Here  $w(t,\cdot) \in H^2\bigl(\R\backslash \{0\}\bigr)$, while   $\eta\in C^{\infty}(\R)$ 
is an even cut-off function, satisfying 
\bel{eta}
\left\{ \bega{rl}
\eta(x)~=~1&\qquad \hbox{if}~~~ |x|\leq 1 ,\\[1mm]
\eta(x)~=~0&\qquad \hbox{if}~~~ |x|\geq 2  ,\\[1mm]
\eta'(x)~\leq~ 0&\qquad \hbox{if}~~~x\geq 0\,.
\enda\right.
\eeq
For future use, it will be convenient to introduce the function
\bel{phi-b}
\phi(x,b)~\doteq~{2\eta(x)\over \pi}\cdot \Big[(|x|+b)\ln(|x|+b)- b\ln b\Big],\qquad\quad  x\in\R,~~b\geq 0.
\eeq
Observe that 
\bel{phib0}
\phi(0,b)~=~0\qquad\qquad\forall~  b\geq 0.\eeq
Our main goal in this section is to solve the Cauchy problem for (\ref{BH-1}) with initial data
\bel{id2}u(0,x)~=~\ov w(x)+\overline{\varphi}(x),
\eeq
where
\bel{o-vp}
\ov w\in H^2\bigl(\R\backslash\{0\}\bigr),\qquad \qquad
\ov\vp(x)~=~\Big(c_1\cdot \chi_{\strut ]-\infty,0[}+c_2\cdot \chi_{\strut ]0,+\infty[}\Big)\cdot\phi (x,0),
\eeq
for some constants $c_1$, $c_2\in\R$. 
Note that this reduces to (\ref{pcd1}) in the case  $c_1=c_2=1$.

To handle the more general initial data (\ref{id2})-(\ref{o-vp}), 
we  write the solution of (\ref{BH-1})     in the form
\bel{usw} u(t,x)~=~w(t,x)+\vp^{(w)}(t,x),\eeq
where the corrector term 
$\vp^{(w)}(t,x)$ now depends explicitly on time $t$ and on the strength of the jump 
\bel{JS}\sigma^{(w)}(t)~
\doteq ~ w^-(t) -w^+(t),\qquad\qquad w^{\pm}(t)\,\doteq\,w(t,0\pm).\eeq
To make an appropriate guess for the function $\vp^{(w)}$, we observe that,
by (\ref{HT}), the equation  (\ref{BH-1}) can be approximated by the simpler equation
\bel{BH-22}
u_t+\left(u- {u^-(t)+u^+(t)\over 2}\right) u_x~=~{1 \over\pi} (u^+(t)-u^-(t))\, \ln |x|. \eeq
Indeed, we expect that the solutions of (\ref{BH-1}) and (\ref{BH-22}) with the same initial data
will have the same asymptotic structure near the origin. Their difference will lie in the more regular space $H^2\bigl(\R\backslash\{0\}\bigr)$.
With this in mind, we thus make the {\it ansatz}
\bel{vp1}
\vp^{(w)}(t,x)~\doteq~\phi(x,0)+\Big((c_1-1)\cdot \chi_{\strut ]-\infty,0[}+(c_2-1)\cdot \chi_{\strut ]0,+\infty[}\Big)\cdot \phi\left(x,\ds{\sigma^{(w)}(t)t\over 2}\right).
\eeq
Inserting (\ref{usw}) into (\ref{BH-1}), we obtain an equation for the remaining component $w(t,\cdot)$. Namely
\bel{BH-2}
w_t+a(t,x,w)\cdot w_x~=~F(t,x,w),
\eeq
where $a$ and $F$ are given by
\bel{a}
a(t,x,w)~=~w(t,x)+\vp^{(w)}(t,x)-{w^-(t)+w^+(t)\over 2},
\eeq
\begin{eqnarray}
\label{F}
F(t,x,w)~&=& ~\bfH\left[\vp^{(w)}\right] (t,x) -\vp^{(w)}\vp^{(w)}_{x}(t,x)\\ \nonumber
&+&\left( \bfH\left[w\right](t,x)-\left[ \vp^{(w)}_{t}(t,x)+\left(w(t,x)-{w^-(t)+w^+(t)\over 2}\right)\cdot \vp^{(w)}_{x}(t,x)\right]\right).
\end{eqnarray}
We observe that, in the present case of a solution with a single shock, 
by (\ref{phib0}) the entropy admissibility  condition (\ref{RH1}) reduces to 
\bel{wine}w^-(t)~>~w^+(t).\eeq
Moreover, Definition~\ref{d:12}
is satisfied provided that, along every characteristic curve 
$t\mapsto x(t;t_0,x_0)\not= 0$ obtained by solving
\bel{code}
\dot{x}(t)~=~a(t,x,w),\qquad\quad x(t_0)\,=\,x_0\,,
\eeq
one has 
\bel{w-in}
w(t_0,x_0)~=~\ov w\bigl(x(0;t_0,x_0)\bigr)+\int_{0}^{t_0}F\Big(t,x(t;t_0,x_0),w\bigl(t,x(t;t_0,x_0)\bigr)
\Big)\, dt\,.
\eeq

The first main result of this paper provides the existence and uniqueness of an entropic
solution, locally in time.
\begin{theorem}\label{t:21}
For every $\overline{w}\in H^2\bigl(\R\setminus \{0\}\bigr)$ satisfying $\ov w(0-)-\ov w(0+)>0$ and every $c_1$, $c_2\in \R$,  the Cauchy problem for the Burgers-Hilbert equation (\ref{BH}), with initial condition as in (\ref{id2})-(\ref{o-vp}),  admits a unique piecewise regular solution defined for $t\in [0,T]$, for some $T>0$ sufficiently small, depending only on $M_0$, $\delta_0$, $c_1$, and $c_2$.
\end{theorem} 

The solution to the equivalent equation
(\ref{BH-2}) will be obtained as a limit of a sequence of approximations. 
Namely, consider a sequence of linear approximations constructed as follows.
As a first step, define
\bel{id1} w_1(t,x)~=~\ov w(x)
\qquad\qquad\hbox{for all}~~t\geq 0, ~~x\in\R\,.\eeq
By induction, let $w_n$ be given. We define  $w_{n+1}$ to
be the solution of the linear, non-homogeneous  Cauchy problem
\bel{wne}
w_t + a(t,x,w_n)\cdot w_x~=~F(t,x,w),\qquad w(0,\cdot)~=~\ov w(x).\eeq 
The induction argument requires three steps:
\begi
\item[(i)] Existence and uniqueness of solutions to each linear 
problem (\ref{wne}).

\item[(ii)] A priori bounds on the strong norm $\bigl\|w_n(t)\bigr\|_{H^2(\R\setminus\{0\})}$, 
uniformly valid
for $t\in [0,T]$ and all $n\geq 1$.

\item[(iii)] Convergence  in a weak norm. This will follow from the bound
\bel{conv1}\sum_{n\geq 2}\,\bigl\|w_{n}(t)-w_{n-1}(t)\bigr\|_{H^1(\R\setminus\{0\})}~<~\infty\,.\eeq
\endi
These steps will be worked out in the next two sections.

\section{Preliminary estimates}
\label{sec:3}
\setcounter{equation}{0}
To achieve the above steps (i)-(iii),  we 
establish in this section some key estimates on the right hand side of (\ref{wne}), by splitting it into three parts:
\bel{Fdec}
F(t,x,w)~=~A^{(w)}(t,x)+B^{(w)}(t,x)-C^{(w)}(t,x),
\eeq
where
\bel{AB}
A^{(w)}~\doteq~{\bf H}\left[\vp^{(w)}\right]-\vp^{(w)}\vp_x^{(w)},\qquad
\qquad  B^{(w)}~\doteq~{\bf H}[w]-\left(w-{w^-+w^+\over 2}\right)\cdot \phi_x(x,0)\,,
\eeq
\bel{C}
C^{(w)}~\doteq~\varphi^{(w)}_t+\left(w-{w^-+w^+\over 2}\right)\cdot 
\phi_x\left(x,{\sigma^{(w)}(t)\,t\over 2}\right)\cdot\Big((c_1-1)\cdot \chi_{\strut ]-\infty,0[}+(c_2-1)\cdot \chi_{\strut ]0,+\infty[}\Big).
\eeq
Consider the function  
\bel{g-b}
g_b(x)~=~\chi_{\strut [0,\infty[}(x)\cdot\phi(x,b),\qquad x\in\R, b\geq 0,
\eeq
where $\phi(x,b)$ is given by (\ref{phi-b}).  For every $b\in \big[0,{1\over 2e}\big]$  
one checks that 
the function  $g_b\in \mathcal{C}^{\infty}(\R\backslash \{0\})\cap\mathcal{C}(\R)$ is  negative and 
decreasing on the open interval  $\bigl]
0,{1\over 2e}\bigr[$. Moreover, it satisfies 
\bel{e1-gb}
\mathrm{supp}(g_b)~\subseteq~[0,2],\qquad \qquad
 |g_b(z)|~\leq~ \big|z\ln |z|\big|\quad \forall z\in \left[0,\, {1\over 2e}\right].
\eeq

The next lemma provides some bounds on the Hilbert transform of $g_b$. As usual, by the Landau symbol $\O(1)$ we shall denote a uniformly bounded quantity.

\begin{lemma}\label{b-phi} For every $0\leq b\leq {1\over 2e}$ and $|x|\leq {1\over 2e}$, one has  
\bel{Hg-b1}
\left\{\bega{rlrl}
\big|{\bf H}[g_b](x)\big|&\ds \leq~\ds\O(1),\qquad & \ds \left|{d\over dx}{\bf H}[g_b](x)\right|&
\ds\leq~\O(1)\cdot \ln^2|x|,\\[4mm]
\left|\ds{d^2\over dx^2}{\bf H}[g_b](x)\right|&\leq~\ds \O(1)\cdot \left|{\ln |x|\over x} \right|,
&\qquad  \left|\ds{d^3\over dx^3}{\bf H}[g_b](x)\right|&\leq~\ds \O(1)\cdot \left|{\ln |x|\over x^2} \right|.
\enda\right.
\eeq
Moreover, for every $\delta>0$ sufficiently small one has
\bel{Hg-b2}
\big\|{\bf H}[g_b](\cdot)\big\|_{H^{2}(\R\backslash[-\delta,\delta])}~\leq~\O(1)\cdot \delta^{-2/3},\qquad \big\|{\bf H}[g'_b](\cdot)\big\|_{H^{2}(\R\backslash[-\delta,\delta])}~\leq~\O(1)\cdot \delta^{-7/4}.
\eeq
\end{lemma}
{\bf Proof.} Fix $b\in \left[0,{1\over 2e}\right]$. By  (\ref{HT}), one has
\[
{\bf H}[g_b](x)~=~{1\over \pi}\cdot \int_{0}^{2}g_b'(y)\cdot \ln |x-y|~dy.
\]
Two cases are considered:
\medskip

\n {\bf Case 1:}  If $-\ds{1\over 2e}<x<0$ then we have the estimates 
\begin{eqnarray*}
\big|{\bf H}[g_b](x)\big|&=& {1\over \pi}\cdot \left|\int_{0}^2g_b'(y)\cdot \ln (y-x)~dy\right|~=~{1\over \pi}\cdot \left|\int_{|x|}^{2+|x|} g'_b(x+z)\ln z~dz\right|\\
&=&{1\over \pi}\cdot \left|\int_{|x|}^{2+|x|} {g_b(x+z)\over z}~dz\right|~\leq~\O(1)\cdot\left(\int_{|x|}^{{1\over 2e}}\big|\ln z\big|dz+1\right)~\leq~ \O(1),
\end{eqnarray*}
\begin{eqnarray*}
\left|\ds{d\over dx}{\bf H}[g_b](x)\right|&=& {1\over \pi}\cdot \left|\int_{0}^2{g_b'(y)\over y-x}~dy\right|~=~{1\over \pi}\cdot \left|\int_{|x|}^{2+|x|} {g'_b(x+z)\over z}~dz\right|\\
&=&{1\over \pi}\cdot \left|\int_{|x|}^{2+|x|} {g_b(x+z)\over z^2}~dz\right|~=~\O(1)\cdot\left(\int_{|x|}^{{1\over 2e}}{\big|\ln z\big|\over z}~dz+1\right)~\leq~\O(1)\cdot\ln^2|x|,
\end{eqnarray*}
and similarly
\[
\left|\ds{d^2\over dx^2}{\bf H}[g_b](x)\right|~\leq~\O(1)\cdot \left|{\ln |x|\over x} \right|,\qquad 
\qquad \left|\ds{d^3\over dx^3}{\bf H}[g_b](x)\right|~\leq~\O(1)\cdot \left|{\ln |x|\over x^2}\right|.
\]
{\bf Case 2:} If $0<x<\ds{1\over 2e}$, then we split ${\bf H}[g_b](x)$ into three parts as follows:
\begin{eqnarray*}
{\bf H}[g_b](x)&\ds =~{1\over \pi}\cdot \left(\int_{0}^{x/2}g_b'(y)\ln (x-y)~dy+\mathrm{p.v.}\int_{x/2}^{3x/2}g_b'(y) \ln |y-x|~dy\right)\\
&\qquad \ds + {1\over \pi}\cdot \int_{3x/2}^{2}g_b'(y)\ln (y-x)~dy ~\doteq
~{1\over \pi}\cdot \big(I_1+I_2+I_3\big).
\end{eqnarray*}
We first estimate 
\[
\begin{cases}
|I_1(x)|&=~\left|\ds\int_{0}^{x/2}g_b'(y)\ln (x-y)~dy\right|~\leq~\O(1)\cdot\left|\ds\int_{0}^{x/2}g_b'(y)~dy\right|\cdot \big|\ln x\big|~\leq~\O(1)\cdot x\ln^2x,\\
|I'_1(x)|&=~\left|\ds\int_{0}^{x/2}{g_b'(y)\over x-y}~dy+ {1\over 2}\cdot g_b'\left({x\over 2}\right)\ln\left({x\over 2}\right)\right|~\leq~\O(1)\cdot\ln^2x,
\end{cases}
\]
and similarly
\[
|I''_1(x)|~\leq~\O(1)\cdot\left|{\ln x\over x} \right|,\qquad |I'''_1(x)|~\leq~\O(1)\cdot\left|{\ln x\over x^2} \right|.
\]

By a similar argument, one obtains 
\[
|I_3(x)|~\leq~\O(1),\qquad \big|I'_3(x)\big|~=~\left|\int_{3x/2}^{2}{g_b'(y)\over y-x}~dy+{3 \over 2}\cdot g'_b\left({3x\over 2}\right)\cdot
\ln\left({x\over 2}\right)\right|~\leq~\O(1)\cdot\ln^2x,
\]
and
\[
\begin{cases}
\big|I''_3(x)\big|&\leq~\left|\ds\int_{3x/2}^{2}{g_b'(y)\over (y-x)^2}~dy\right|+\ds\O(1)\cdot\left|{\ln x\over x}\right|~\leq~\O(1)\cdot\left|{\ln x\over x}\right|,\\
\big|I'''_3(x)\big|&\leq~\left|\ds\int_{3x/2}^{2}{g_b'(y)\over (y-x)^3}~dy\right|+\ds\O(1)\cdot\left|{\ln x\over x^2}\right|~\leq~\O(1)\cdot\left|{\ln x\over x^2}\right|.
\end{cases}
\]

Finally, using the fact that $g_b'$ is concave, we obtain
\begin{eqnarray*}
\big|I_2(x)\big|&=& \left|\lim_{\ve\to 0+}\int_{x/2}^{x-\ve}g_b'(y)\ln (x-y)dy+\int_{x+\ve}^{3x/2}g_b'(y)\ln (y-x)dy\right|\\
&=&\left|\lim_{\ve\to 0+}\int_{\ve}^{x/2}\big[g'_b(x-z)+g'_b(x+z)\big]\cdot \ln z~dz\right|~\leq~ 2|g'_b(x)|\cdot \left|\int^{x/2}_{\ve}\ln zdz\right|~\leq~\O(1),
\end{eqnarray*}
\begin{eqnarray*}
\big|I'_2(x)\big|&=& \left|\lim_{\ve\to 0+}\int_{\ve}^{x/2}\big[g''_b(x-z)+g''_b(x+z)\big]\ln zdz+{1\over 2}[g'_b(x/2)+g'_b(3x/2)]\cdot\ln(x/2)\right|\\
&\leq&\O(1)\cdot\left( |g''_b(x/2)|\cdot \left|\int^{x/2}_{\ve}\ln zdz\right|+\ln^2 x\right)~\leq~\O(1)\cdot\ln^2 x,
\end{eqnarray*}
and
\[
\begin{cases}
\big|I''_2(x)\big|&\leq~\ds\left|\lim_{\ve\to 0+}\int_{\ve}^{x/2}\big[g'''_b(x-z)+g'''_b(x+z)\big]\ln zdz\right|+\O(1)\cdot\left|{\ln x\over x}\right|~\leq~\O(1)\cdot\left|{\ln x\over x}\right|,\cr
\big|I'''_2(x)\big|&\leq~\ds\left|\lim_{\ve\to 0+}\int_{\ve}^{x/2}\big[g''''_b(x-z)+g''''_b(x+z)\big]\ln zdz\right|+\O(1)\cdot\left|{\ln x\over x^2}\right|~\leq~\O(1)\cdot\left|{\ln x\over x^2}\right|.
\end{cases}
\]
We thus  achieve the same estimates as in Case 1, and this yields (\ref{Hg-b1}). 

Finally, the function $g_b$ is continuous with compact support and smooth 
outside  the origin. Therefore, the Hilbert transform ${\bf H}[g_b]$ is smooth outside the origin.
As $|x|\to\infty$, one has
\[
{\bf H}[g_b](x)~=~\O(1)\cdot |x|^{-1},\qquad {d^k\over dx^k} \left({\bf H}[g_b]\right)(x)~=~\O(1)\cdot x^{-(k+1)},\quad k=1,2,3.
\]
Thus,  (\ref{Hg-b1}) yields (\ref{Hg-b2}).
\endproof

\begin{remark}\label{dgb} {\rm For every $0<b\leq {1\over 2e}$, one has 
\[
{d\over db} g_b(x)~=~{2\eta(x) \over \pi}[\ln(x+b)-\ln(b)]\qquad \forall x>0.
\]
Since 
\[ \left|\ln(x+b)-\ln(b)\right|\,\leq\, {x\over b},\]
the same arguments used in the proof of Lemma \ref{b-phi} yield that, for $0<|x|\leq \ds{1\over 2e}$,
\[
\bega{c}
\ds\left|{\bf H}\left[{d\over db}g_b\right](x)\right|~=~\ds \O(1)\cdot{1\over b},\qquad\qquad \left|{d\over dx}{\bf H}\left[\left({d\over db}g_b\right)\right](x)\right|~=~ \ds \O(1)\cdot {\left |\ln(x)\right|\over b}, \\[4mm]
\ds \left|{d^2\over dx^2}{\bf H}\left[\left({d\over db}g_b\right)\right](x)\right|~=~\ds \O(1)\cdot \left |{\ln(x) \over bx}\right|.
\enda
\]
Moreover, for $\delta>0$ sufficiently small,
\[
\left\|{\bf H}\left[{d\over db}g_b\right]\right\|_{H^{1}(\R\backslash [-\delta,\delta])}\leq~\O(1)\cdot{1\over b},\qquad \left\|{\bf H}\left[{d\over db}g_b\right]\right\|_{H^{2}(\R\backslash [-\delta,\delta])}\leq~\O(1)\cdot {1\over b}\cdot\delta^{-2/3}.
\]
}
\end{remark}

The next lemma provides some a priori estimates on the function $F= F(t,x,w)$ introduced at (\ref{F}).

\begin{lemma}\label{F1-g}  Let $w: [0,T]\times\R\to\R$ such that  $w(t,\cdot)\in H^2(\R\backslash\{0\})$
for all $t\in[0,T]$, 
\[
\|w(t,\cdot)\|_{H^2(\R\backslash\{0\})}~\leq~M_0,\qquad  \sigma^{(w)}(t)~\doteq~w(t,0+)-w(t,0-)~>~0.
\]
Moreover, assume that  $0<T<\ds{1\over 4eM_0}$ and that $\sigma^{(w)}(\cdot)$ is locally Lipschitz on $]0,T]$.

 Then there exists a constant $C_1>0$, depending on $M_0$, $\delta_0$, $c_1$, and $c_2$ such that, for a.e.~$t\in [0,T]$ and $|x|<{1\over 2e}$, one has 
\bel{F1}
\begin{cases}
|F(t,x,w)|&\leq~C_1\cdot \left(\ds (1+M_0)\cdot |\ln t| +{|\dot{\sigma}^{(w)}(t)|\over \sigma^{(w)}(t)}\cdot \left| x\right|\right),\\[5mm]
|F_x(t,x,w)|&\leq~C_1\cdot \left(\ds (1+M_0)\cdot |x|^{-1/4}+{|\dot{\sigma}^{(w)}(t)|\over \sigma^{(w)}(t)}\right).
\end{cases}
\eeq
Furthermore, for every $\delta>0$ sufficiently small
\bel{F2}
\|F(t,x,w)\|_{H^2(\R\backslash [-\delta,\delta])}~\leq~C_1\cdot\left[\left(\ds 1+ M_0+{|\dot{\sigma}^{(w)}(t)|\over |\sigma^{(w)}(t)|}\right)\cdot \ds \delta^{-2/3}+(1+M_0)\cdot |\ln (t)|\right].
\eeq

\end{lemma}
{\bf Proof.}  According to (\ref{Fdec}), the function $F$ can be decomposed as the sum of three terms,
which will be estimated separately.
\v
{\bf 1.} By the analysis in  \cite[Section 3]{BZ},  for every $(t,x)\in [0,T]\times \left[-{1\over 2e},{1\over 2e}\right]$ one has
 \bel{Bb1}
\begin{cases}
\left|B^{(w)}(t,x)\right|~\leq~\O(1)\cdot M_0,\qquad \qquad\left|B^{(w)}_x(t,x)\right|~\leq~ \O(1)\cdot M_0\cdot |\ln|x||^2,\\[3mm]
\left\|B^{(w)}(t,\cdot)\right\|_{H^2(\R\backslash [-\delta,\delta])}~\leq~\O(1)\cdot M_0\cdot\delta^{-2/3}.
\end{cases}
 \eeq
{\bf 2.} Next,  we estimate  $C^{(w)}(t,x)$.  For every $0<x<\ds{1\over 2e}$, we have
\bel{Cx>0}
C^{(w)}(t,x)~=~E^{(w)}_2(t)+(c_2-1)\cdot \left[{\dot{\sigma}^{(w)}(t)t\over 2}+(w(t,x)-w(t,0+)) \right]\cdot\phi_x\left(x,{\sigma^{(w)}(t)t\over 2}\right),
\eeq
where we define
\[
E^{(w)}_2(t)~\doteq~{1-c_2\over 2}\cdot\phi_x\left({\sigma^{(w)}(t)t\over 2},0\right)\cdot\left(\dot{\sigma}^{(w)}(t)t+\sigma^{(w)}(t)\right).
\]
Since $\|w(t,\cdot)\|_{H^2(\R\backslash\{0\})}\leq M_0$, one has  $|\sigma^{(w)}(t)|\leq 2M_0$,
\bel{C(t,0)}
\big|C^{(w)}(t,0+)\big|~=~\left|{1-c_2\over 2}\cdot \phi_x\left({\sigma^{(w)}(t)t\over 2},0\right) \cdot \sigma^{(w)}(t)\right|~\leq~\O(1)\cdot (1+M_0)\cdot|\ln t|,
\eeq
and
\[
\big|C^{(w)}(t,x)\big|~\leq ~\O(1)\left(\ds (1+M_0)\cdot |\ln t| +{|\dot{\sigma}^{(w)}(t)|\over \sigma^{(w)}(t)}\cdot \left| x\right|\right).
\]
Moreover, observing that 
\bel{Ca1}
\begin{cases}
C_x^{(w)}(t,x)&=~(c_2-1)\cdot \ds{d\over dx}\left[\left({\dot{\sigma}^{(w)}(t)t\over 2}+(w(t,x)-w(t,0+))\right)\cdot \phi_x\left(x,{\sigma^{(w)}(t)t\over 2}\right)\right],\\[5mm]
C_{xx}^{(w)}(t,x)&=~(c_2-1)\cdot \ds{d^2\over dx^2}\left[\left({\dot{\sigma}^{(w)}(t)t\over 2}+(w(t,x)-w(t,0+))\right)\cdot \phi_x\left(x,{\sigma^{(w)}(t)t\over 2}\right)\right],
\end{cases}
\eeq
we estimate
\bel{C-b2}
\begin{cases}
\left|C_x^{(w)}(t,x)\right|&\leq~\O(1)\cdot \ds\left( M_0\cdot |x|^{-1/4}+{|\dot{\sigma}^{(w)}(t)|\over \sigma^{(w)}(t)}\right),\\[5mm]
\left|C_{xx}^{(w)}(t,x)\right|&\leq~\O(1)\cdot \left[\left(M_0+\ds{|\dot{\sigma}^{(w)}(t)|\over |\sigma^{(w)}(t)|}\right)\cdot |x|^{-1}+|w_{xx}(t,x)|\cdot \big|\ln|x|\big|\right].
\end{cases}
\eeq
Similarly, for every $\ds-{1\over 2e}<x<0$, we have 
\[
C^{(w)}(t,x)~=~E^{(w)}_1(t)-(c_1-1)\cdot \left[{\dot{\sigma}^{(w)}(t)t\over 2}-(w(t,x))-w(t,0-)) \right]\cdot\phi_x\left(x,{\sigma^{(w)}(t)t\over 2}\right),
\]
\[
E^{(w)}_1(t)~\doteq~{1-c_1\over 2}\cdot\phi_x\left({\sigma^{(w)}(t)t\over 2},0\right)\cdot\left(\dot{\sigma}^{(w)}(t)t+\sigma^{(w)}(t)\right).
\]
This yields the same bounds as in (\ref{Cx>0})-(\ref{C-b2}). We thus conclude
\bel{C-b3}
\left\|C^{(w)}(t,\cdot)\right\|_{H^2(\R\backslash [-\delta,\delta])}~\leq~\O(1)\cdot\left[\left(\ds M_0+{|\dot{\sigma}^{(w)}(t)|\over |\sigma^{(w)}(t)|}\right)\cdot \ds \delta^{-1/2}+M_0\cdot |\ln (t)|\right].
\eeq

\n {\bf 3.} Finally, to obtain a  bound on  $A^{(w)}$ we observe that, by \eqref{e1-gb},
$$\left|\vp^{(w)}(t,x)\right|\leq\O(1)\cdot  \big|x \ln|x|\big|,\qquad (t,x)\in [0,T]\times \left[ -{1\over 2e},\,{1\over 2e}
\right].$$ 
This leads to the estimates
\[
\left|\vp^{(w)}\vp_x^{(w)}\right|\leq\O(1)\cdot |x|^{1/2},~~ \left|\Big(\vp^{(w)}\vp_x^{(w)}\Big)_{x}\right|\leq\O(1)\cdot \ln^2|x|,~~ \left|\Big(\vp^{(w)}\vp_x^{(w)}\Big)_{xx}\right|\leq\O(1)\cdot  \left|{\ln |x|\over x} \right|\,.
\]
Thus,
\[
\left\|\vp^{(w)}(t,\cdot)\vp_x^{(w)}(t,\cdot)\right\|_{H^2(\R\backslash [-\delta,\delta])}~\leq~\O(1)\cdot\delta^{-2/3}.
\]
On the other hand, if $0<T<\ds{1\over 4eM_0}$, then  $\ds\sup_{t\in [0,T]}{\sigma^{(w)}(t)t\over 2}\leq{1\over 2e} $ and 
Lemma \ref{b-phi} implies for all $t\in [0,T]$ and $|x|\leq\ds{1\over 2e}$, that 
\bel{H-b1}
\left|{\bf H}[\vp^{(w)}(t,\cdot)](x)\right|~\leq~\O(1),\qquad \left|{d\over dx}{\bf H}[\vp^{(w)}(t,\cdot)](x)\right|~\leq~\O(1)\cdot\ln^2|x|,
\eeq
\bel{H-b2}
\left|{d^2\over dx^2}{\bf H}[\vp^{(w)}(t,\cdot)](x)\right|~\leq~\O(1)\cdot \left|{\ln |x|\over x} \right|,\qquad \left\|{\bf H}[\vp^{(w)}(t,\cdot)]\right\|_{H^2(\R\backslash [-\delta,\delta])}~\leq~\O(1)\cdot \delta^{-2/3}.
\eeq
Therefore, combining (\ref{Bb1})-(\ref{H-b2}), we  obtain (\ref{F1})-(\ref{F2}). This  completes the proof.
\endproof

Our third lemma estimates the change in the function $F= F(t,x,w)$ as $w(\cdot)$ 
takes different values.  
These estimates will play a key role in the proof of convergence of the approximations
considered at (\ref{conv1}).
 
\begin{lemma}\label{FF} Let $w_i: [0,T]\times\R\to\R$, $i=1,2$ such that for all $t\in [0,T]$,
 $w_i(t,\cdot)\in H^2(\R\backslash\{0\})$ and
\[
\|w_i(t,\cdot)\|_{H^2(\R\backslash\{0\})}~\leq~M_0,\qquad \left|\sigma^{(w_i)}(t)\right|~\geq~\delta_0. \]
Moreover, assume that $\sigma^{(w_i)}$ is  locally Lipschitz on $]0,T]$ and that there exists a function $K(t)$ such that 
\[
\qquad \left|\dot{\sigma}^{(w_i)}(t)\right|~\leq~K(t)\qquad a.e.~t\in (0,T).
\]
Set  
$$z\doteq w_2-w_1,\quad \sigma^{(z)}\doteq \sigma^{(w_2)}-\sigma^{(w_1)},\quad
 M_1(t)\doteq \|z(t,\cdot)\|_{H^1(\R\backslash\{0\})},\quad M_2(t)\doteq \|z(t,\cdot)\|_{H^2(\R\backslash\{0\})}.$$ 

Then there exists a constant $C_2>0$, depending on $M_0$, $\delta_0$, $c_1$, and $c_2$ such that, for every $ x\in \left[-{1\over 2e},{1\over 2e}\right]$ and a.e. $t\in [0,T]$, one has
\begin{equation}\label{F1-F2-0}
\bigl| F(t,x,w_2)-F(t,x,w_1)\bigr|~\leq~C_2\cdot \left[\big|\dot{\sigma}^{(z)}(t)\big|\cdot |x|+M_1(t)\cdot  \bigl(|\ln t|+K(t)\bigr)\right].
\end{equation}
Moreover, for every $\delta>0$ sufficiently small, it holds 
\begin{equation}\label{F1-F2-1}
\begin{cases}
\left\|{\bf H}\left[\vp^{(w_2)}(t,\cdot)-\vp^{(w_1)}(t,\cdot)\right]\right\|_{H^1(\R\backslash [-\delta,\delta])}+\left\|B^{(w_2)}-B^{(w_1)}\right\|_{H^1(\R\backslash[-\delta,\delta])}~\leq~\ds C_2\cdot {M_1(t)\over \delta^{1/2}},\\[3mm]
\left\|F(t,\cdot,w_2)-F(t,\cdot,w_1)\right\|_{H^2(\R\backslash [-\delta,\delta])}~\leq~\ds C_2\cdot\left(M_2(t)\cdot \Bigg(\big|\ln t\big|+{1+K(t)\over \delta^{2/3}}\right)+{|\dot{\sigma}^{(z)}(t)|\over \delta^{1/2}}\Bigg).
\end{cases}
\end{equation}
\end{lemma}
{\bf Proof.} {\bf 1.} For notational convenience, we set  
\[
{\bf A}^{(z)}~\doteq~A^{(w_2)}-A^{(w_1)},\quad {\bf B}^{(z)}~\doteq~B^{(w_2)}-B^{(w_1)},\quad {\bf C}^{(z)}~\doteq~C^{(w_2)}-C^{(w_1)}.
\]
From  \cite[Section 3]{BZ},  for every $(t,x)\in [0,T]\times \left[-{1\over 2e},{1\over 2e}\right]$, it holds
 \bel{B2-B1}
\begin{cases}
\left|{\bf B}^{(z)}(t,x)\right|=\O(1)\cdot M_1(t),\qquad  \left\|{\bf B}^{(z)}(t,\cdot)\right\|_{H^1(\R\backslash [-\delta,\delta])}~\leq~\ds\O(1)\cdot {M_1(t)\over \delta^{1/2}},\\[3mm]
\left\|{\bf B}^{(z)}(t,\cdot)\right\|_{H^2(\R\backslash [-\delta,\delta])}\leq\O(1)\cdot \ds{M_2(t)\over \delta^{2/3}}.
\end{cases}
\eeq
{\bf 2.} We now provide bounds on ${\bf C}^{(z)}(t,x)$. 
By (\ref{Cx>0})-(\ref{Ca1}), for every $0<x<\ds{1\over 2e}$ one has 
\begin{multline}\label{Cz1}
\big|{\bf C}^{(z)}(t,x)\big|~\leq~\O(1)\left[{|\dot{\sigma}^{(z)}(t)|\over\delta_0}\cdot x+M_1(t)\cdot \left({K(t)x\over\delta^2_0}+{M_0x+M_0^2+1\over \delta_0}+\big|\ln  t\big|+\big|x^{1/2}\ln x\big|\right)\right],
\end{multline}
\begin{multline}\label{Cz2}
\big|{\bf C}^{(z)}_x(t,x)\big|\leq\O(1)\cdot \Bigg[{|\dot{\sigma}^{(z)}(t)|\over\delta_0}+\left({K(t)\over\delta^2_0}+{M_0\over \delta_0}\right)\cdot \big|\sigma^{z}(t)\big|+{\big|z(t,x)-z(t,0+)\big|\over x}\\
+|z_x(t,x)|\cdot \left(\left|\ln(t)\right|+{1\over \delta_0}+M_0\right)\Bigg]\\
\leq\O(1)\cdot \left[{|\dot{\sigma}^{(z)}(t)|\over\delta_0}+\left({K(t)\over\delta^2_0}+{M_0\over \delta_0}\right)\cdot M_1(t)+{M_1(t)\over x^{1/2}}+|z_x(t,x)|\cdot \left(\left|\ln(t)\right|+{1\over \delta_0}+M_0\right)\right],
\end{multline}
\begin{multline}\label{Cz3}
\big|{\bf C}^{(z)}_{xx}(t,x)\big|\leq\O(1)\cdot \Bigg\{{|\dot{\sigma}^{(z)}(t)|\over\delta_0 x}+ |z_{xx}(t,x)|\cdot  \left(\left|\ln(t)\right|+{1\over \delta_0}+M_0\right)+{\big|z_x(t,x)\big|\over x}+{|z(t,x)-z(t,0+)|\over  x^2} \\
\qquad\qquad\qquad\qquad+\big| \sigma^{(z)}(t)\big|\cdot \left({|w_{1,xx}(t,x)|\over \delta_0}+{K(t)\over \delta^2_0 x}+{|w_{1,x}(t,x)|\over \delta_0 x}+{{\big|w_{1}(t,x)-w_{1}(t,0+)\big|}\over \delta_0 x^2}\right)\Bigg\}\\
\leq\O(1)\cdot\Bigg\{\left({|\dot{\sigma}^{(z)}(t)|\over \delta_0}+{K(t)M_1(t)\over \delta^2_0}+{M_1(t) M_0\over\delta_0}+M_2(t)\right) \cdot{1\over x}\\
~~\qquad\qquad\qquad\qquad+|z_{xx}(t,x)|\cdot \left(\left|\ln(t)\right|+{1\over \delta_0}+M_0\right)+\big|w_{1,xx}(t,x)\big|\cdot{M_1(t)\over \delta_0}\Bigg\}.
\end{multline}
For every $-{1\over 2e}<x<0$, by a similar argument, we obtain the same bounds  as in (\ref{Cz1})-(\ref{Cz3}). Therefore
\bel{Cz4}
\begin{cases}
\left\|{\bf C}^{(z)}(t,\cdot)\right\|_{H^1(\R\backslash\{0\})}&\leq~\O(1)\cdot \left(\ds M_1(t)\cdot\left(\big|\ln  t\big|+{M_0+1\over \delta_0}+{K\over \delta^2_0}\right)+{|\dot{\sigma}^{(z)}(t)|\over\delta_0}\right),\\[4mm]
\left\|{\bf C}^{(z)}(t,\cdot)\right\|_{H^2(\R\backslash [-\delta,\delta])}&\leq \O(1)\cdot \left(\ds M_2(t)\cdot\left[\big|\ln  t\big|+\left({M_0+1\over \delta_0}+{K\over \delta^2_0}\right)\cdot {1\over \sqrt{\delta}}\right]+ {|\dot{\sigma}^{(z)}(t)|\over\delta_0  \sqrt{\delta}}\right).
\end{cases}
\eeq

{\bf 3.} To achieve bound on  ${\bf A}^{(z)}$,  for $0<x<{1\over 2e}$ we compute
\[
\left|\vp^{(w_2)}\vp_x^{(w_2)}-\vp^{(w_1)}\vp_x^{(w_1)}\right|~\leq~\O(1)\cdot {\big|\sigma^{(z)}(t)\big|\over \delta_0}\cdot \big|x\ln |x|\big|~\leq~\O(1)\cdot {M_1(t)\over\delta_0} \cdot \big|x\ln |x|\big|,
\]
\[
\left|\left(\vp^{(w_2)}\vp_x^{(w_2)}-\vp^{(w_1)}\vp_x^{(w_1)}\right)_x\right|~\leq~\O(1)\cdot {\big|\sigma^{(z)}(t)\big|\over\delta_0}\cdot \big|\ln |x|\big|~\leq~\O(1)\cdot {M_1(t)\over \delta_0} \cdot  \big|\ln |x|\big|,
\]
\[
\left|\left(\vp^{(w_2)}\vp_x^{(w_2)}-\vp^{(w_1)}\vp_x^{(w_1)}\right)_{xx}\right|~\leq~\O(1)\cdot {\big|\sigma^{(z)}(t)\big|\over \delta_0}\cdot \left|{\ln |x|\over  x}\right|~\leq~\O(1)\cdot {M_1(t)\over \delta_0}  \cdot \left|{\ln |x|\over x}\right|.
\]
This yields 
\[
\begin{cases}
\left\|\vp^{(w_2)}\vp_x^{(w_2)}-\vp^{(w_1)}\vp_x^{(w_1)}\right\|_{H^1(\R\backslash\{0\})}&\leq~\ds\O(1)\cdot {M_1(t)\over \delta_0},\\[4mm]
\left\|\vp^{(w_2)}\vp_x^{(w_2)}-\vp^{(w_1)}\vp_x^{(w_1)}\right\|_{H^2(\R\backslash[-\delta,\delta])}&\ds\leq~\O(1)\cdot {M_1(t)\over \delta_0}\cdot \delta^{-2/3}.
\end{cases}
\]
On the other hand,  for $0<x<{1\over 2e}$ we observe  that 
\begin{eqnarray}\nonumber
\vp^{(w_2)}(t,x)-\vp^{(w_1)}(t,x)&=&(c_2-1)\cdot  \left[\phi\left(x,\ds{\sigma^{(w_2)}(t)t\over 2}\right)-\phi\left(x,\ds{\sigma^{(w_1)}(t)t\over 2}\right)\right]\\ \nonumber
&=&{c_2-1\over 2}\cdot \left(\int_{0}^{1}\phi_b\left(x,{\sigma^{(w_1)}(t)t\over 2}+\tau\cdot {\sigma^{(z)}(t)t\over 2} \right)d\tau\right)\cdot \sigma^{(z)}(t)t\\ \label{con:vpg}
&=&{c_2-1\over 2}\cdot \left(\int_{0}^{1}{d\over db}g_{b_{\tau}}(x)~d\tau\right)\cdot \sigma^{(z)}(t)t,
\end{eqnarray}
with $\ds b_{\tau}={\sigma^{(w_1)}(t)t\over 2}+\tau\cdot {\sigma^{(z)}(t)t\over 2}$. Thus, by Remark \ref{dgb} it follows
\[
\begin{cases}
\left|{\bf H}\left[\vp^{(w_2)}(t,\cdot)-\vp^{(w_1)}(t,\cdot)\right](x)\right|&\leq~\ds\O(1)\cdot {M_1(t)\over \delta_0},\\[4mm]
\left|\ds{d\over dx}{\bf H}\left[\vp^{(w_2)}(t,\cdot)-\vp^{(w_1)}(t,\cdot)\right](x)\right|&\leq~\ds\O(1)\cdot M_1(t)\cdot{\left|\ln |x| \right|\over \delta_0},
\end{cases}
\]
and
\[
\begin{cases}
\ds\left\|{\bf H}\left[\vp^{(w_2)}(t,\cdot)-\vp^{(w_1)}(t,\cdot)\right]\right\|_{H^1(\R\backslash [-\delta,\delta])}&\leq~\ds\ds\O(1)\cdot M_1(t)\cdot{1\over \delta_0},\\[4mm]
\ds\left\|{\bf H}\left[\vp^{(w_2)}(t,\cdot)-\vp^{(w_1)}(t,\cdot)\right]\right\|_{H^2(\R\setminus [-\delta,\delta])}&\leq~\ds\ds\O(1) \cdot {M_1(t)\over \delta_0}\cdot \delta^{-2/3}.
\end{cases}
\]
Similarly, one gets the same estimate for $\ds-{1\over 2e}<x<0$. Therefore (\ref{AB}) yields 
\bel{A2-A1}
\begin{cases}
\left|{\bf A}^{(z)}(t,x)\right|~\leq~\ds\O(1)\cdot {M_1(t)\over \delta_0},\qquad \left|{\bf A}_x^{(z)}(t,x)\right|~\leq~\O(1)\cdot M_1(t)\cdot {\big|\ln |x|\big|\over\delta_0},\\[4mm]
\left\|{\bf A}^{(z)}(t,\cdot)\right\|_{H^1(\R\backslash[-\delta,\delta])}~\leq~\ds\O(1)\cdot {M_1(t)\over \delta_0}\,,\\[4mm]
\left\|{\bf A}^{(z)}(t,\cdot)\right\|_{H^2(\R\backslash[-\delta,\delta])}~\leq~\ds\O(1)\cdot {M_1(t)\over \delta_0} \cdot\delta^{-2/3}\,.
\end{cases}
\eeq
Finally, combining the estimates (\ref{B2-B1})-(\ref{A2-A1}), we  obtain (\ref{F1-F2-0})-(\ref{F1-F2-1}).
\endproof

\section{Proof of Theorem~\ref{t:21}}
\label{sec:4}
\setcounter{equation}{0}
In this section we give a proof of Theorem~\ref{t:21} by constructing a solution to the 
Cauchy problem \eqref{BH-1}
with general initial data of the form  (\ref{id2})-(\ref{o-vp}), locally in time.   This solution will be obtained
as limit of a Cauchy sequence of approximate solutions $w_n(t,x)$, following the steps (i)--(iii) 
outlined at the end of Section~\ref{sec:2}. 
\v
{\bf Step 1.}  Consider any initial profile $\ov w\in  H^2(\R\backslash \{0\})$. Let 
$\delta_0,M_0>0$ be the constants defined by the identities
\bel{ic}
\overline{w}(0-)-\overline{w}(0+)~=~6\delta_0,\qquad \|\overline{w}\|_{H^2(\R\backslash \{0\})}~=~{M_0\over 2}
\,.
\eeq
Given two constants $c_1$, $c_2\in\R$, the corresponding  initial data of the form  (\ref{id2})-(\ref{o-vp}) is
\[
u(0,x)~=~\overline{w}(x)+\left(c_1\cdot \chi_{\strut ]-\infty,0[}+c_2\cdot \chi_{\strut ]0,+\infty[}\right)\cdot\phi (x,0).
\]
Moreover, let $w_n:[0,T]\times \R\to\R$ be a function such that 
\bel{wn-c}
|w_n(t,0\pm)-\overline{w}(0\pm)|~\leq~\delta_0,\qquad \|w_n(t,\cdot)\|_{H^2(\R\backslash\{0\})}~\leq~M_0,
\qquad t\in[0,T].
\eeq
Set  $\sigma_n(t)\doteq w_n(t,0-)-w_n(t,0+)$. As in (\ref{vp1}), the correction term associated to $w_n$ is denoted by 
\bel{vpnn}
\varphi_n(t,x)~\doteq~\phi(x,0)+\Big((c_1-1)\cdot \chi_{\strut ]-\infty,0[}+(c_2-1)\cdot \chi_{\strut ]0,+\infty[}\Big)\cdot \phi\left(x,\ds{\sigma_n(t)t\over 2}\right).
\eeq
In this  step, we will establish the existence and uniqueness of solutions to the linear  problem (\ref{wne}).

We begin by observing that
the speed of all characteristics for (\ref{wne}) is  
\[
a_n(t,x)~\doteq~a(t,x,w_n)~=~ \vp_n(t,x)+w_n(t,x) - {w_n^-(t)+w_n^+(t)\over 2}\,,
\]
where $\varphi_n(t,x)\doteq \varphi^{(w_n)}(t,x)$, the correction term associated to 
$w_n$. From (\ref{vpnn}) and (\ref{wn-c}) it follows that $\varphi_n(t,0)=0$ and 
\[
-4\delta_0~\leq~a_n(t,0+)~=~-{\sigma_n(t)\over 2}~\leq~-2\delta_0,\qquad 2\delta_0~\leq~a_n(t,0-)~=~{\sigma_n(t)\over 2}~\leq~4\delta_0.
\]
Furthermore, for any given $(t,x)\in [0,T]\times \,]0,{1\over 2e}]$, we estimate, using \eqref{e1-gb}, 
\begin{multline*}
\left|a_n(t,x)-a_n(t,0+)\right|~\leq~{2(2+|c_2|) |x\ln x| \over \pi}+ \int_{0}^{x}|w_{n,x}(t,y)|dy \\
~\leq~ {2(2+|c_2|) |x\ln x| \over \pi}+ x^{1/2}\cdot \left(\int_{0}^{x}|w_{n,x}(t,y)|^2dy\right)^{1/2}~\leq~(2+|c_2|+M_0)\cdot \sqrt{x}.
\end{multline*}
Similarly, we also have 
\[
\left|a_n(t,x)-a_n(t,0-)\right|~\leq~(2+|c_1|+M_0)\cdot \sqrt{|x|},\qquad (t,x)\in [0,T]\times \left[-{1\over 2e},0\right[\,.
\]
In particular, setting 
\bel{dt1}
\delta_1~\doteq~\ds{1\over 4}\cdot\left({\delta_0\over 4+|c_1|+|c_2|+M_0}\right)^2~\leq~{1\over 16},
\eeq
we  have 
\bel{b-c1}
\begin{cases}
-5\delta_0~\leq~a_n(t,x)~\leq~-\delta_0,& (t,x)\in [0,T]\times \,]0,2\delta_1],\\[3mm]
\delta_0~\leq~a_n(t,x)~\leq~5\delta_0,& (t,x)\in [0,T]\times [-2\delta_1,0[\,.
\end{cases}
\eeq
Next, choose
\bel{T-as1}
0~<~T~<~ \min\left\{{\delta_1\over 10\delta_0},{1\over 2e}\right\},
\eeq
and denote by  $t\mapsto x(t;t_0,x_0)$ the solution to the Cauchy problem 
\bel{char}
\dot{x}(t)~=~a_n(t,x(t)),\qquad x(t_0)~=~x_0.
\eeq
By (\ref{b-c1}) it follows 
\bel{x-b1}
\delta_0 (t_0-t)~\leq~|x(t;t_0,x_0)-x_0|~\leq~5\delta_0 (t_0-t),\qquad 
 |x_0|\leq \delta_1,~~ 0\leq t\leq t_0\leq T.
\eeq
The next lemma provides the Lipschitz continuous dependence of the characteristic curves considered at (\ref{char}).
\begin{lemma} 
Let $w_n$, $\vp_n$ be as in 
(\ref{wn-c})-(\ref{vpnn}). Then there exists a constant $K_1>0$, depending on $M_0$, $\delta_0$, $c_1$, and $c_2$,  such that, for any $x_1$, $x_2\in \,[-\delta_1,0[\,$ or $x_1$, $x_2\in \,]0,\delta_1]$, one has
\bel{Lipx}
\big|x(t;\tau,x_2)-x(t;\tau,x_1)\big|~\leq~K_1\cdot |x_2-x_1|,\qquad \forall~ 0\leq t\leq \tau \leq T.
\eeq
\end{lemma}
{\bf Proof.} We shall prove (\ref{Lipx}) for $x_1,x_2\in [-\delta_1,0[\,$, the other case 
being entirely similar. For any $-\delta_1\leq z_1<z_2<0$, it holds 
\begin{eqnarray*}
\big|a_n(t,z_2)-a_n(t,z_1)\big|&\leq& |w_n(t,z_2)-w_n(t,z_1)|+\Big|\varphi_n(t,z_2)-\varphi_n(t,z_1)\Big|\\
&\leq&\left[M_0+ \ds{4+2|c_1|\over \pi}\cdot \big(1+\big|\ln|z_2|\big|\big)\right]\cdot |z_2-z_1|.
\end{eqnarray*}
Therefore, from (\ref{x-b1}), it follows 
\[
{d\over dt}\big|x(t;\tau,x_2)-x(t;\tau,x_1)\big|~\leq~\left[M_0+ \ds{4+2|c_1|\over \pi}\cdot \big(1+\big|\ln|\delta_0(\tau-t)|\big|\big)\right]\cdot \big|x(t;\tau,x_2)-x(t;\tau,x_1)\big|,
\]
and this yields (\ref{Lipx}).
\endproof
From (\ref{b-c1}), by the same arguments as in \cite[Lemma 4.1]{BZ}, one obtains:
\begin{lemma}\label{ln-b1} Let $w_n$, $\vp_n$ be as in 
(\ref{wn-c})-(\ref{vpnn}). There exists $T>0$ sufficiently small, depending only on $M_0$, $\delta_0$, $c_1$, $c_2$, such that, for every $\tau\in [0,T]$ and any solution $v$ of the linear equation 
\[
v_t+a_n(t,x)\cdot v_x~~=0,\qquad v(0,\cdot)~=~\bar{v}~\in~ H^2\bigl(
\R\backslash [-\delta_0\tau,\delta_0\tau]\bigr),
\]
one has
\[
\|v(\tau,\cdot)\|_{H^2\bigl(\R\backslash\{0\}\bigr)}~\leq~{3\over 2}\cdot \|\bar{v}\|_{H^2
\bigl(\R\backslash [-\delta_0\tau,\delta_0\tau]\bigr)}.
\]
\end{lemma}

\medskip

{\bf Step 2.} Consider a sequence of approximate solutions $w^{(k)}$ to (\ref{wne}),
inductively defined as follows. 
\begin{itemize}
\item  $w^{(1)}(t,\cdot)~\doteq~\overline{w}(\cdot)$ for all $t\geq 0$.
\item For every $k\geq 1$, $w^{(k+1)}(t,\cdot)$ solves the linear equation 
\[
w_t+a_n(t,x)\cdot w_x~=~F^{(k)}(t,x),\qquad w(0,\cdot)~=~\overline{w}(\cdot)
\]
with $F^{(k)}(t,x)\doteq F\left(t,x,w^{(k)}\right)$.  Equivalently, $w^{(k+1)}$
satisfies the integral identiies
\bel{w-k+1}
w^{(k+1)}(t_0,x_0)~=~\overline{w}(x(0;t_0,x_0))+\int_{0}^{t_0}F^{(k)}(t,x(t;t_0,x_0))dt.
\eeq
\end{itemize}
The following lemma provides a priori estimates on  $w^{(k)}$, uniformly valid  for all $k\geq 1$. 
\begin{lemma}\label{wk-p} Let $w_n$, $\vp_n$ be as in 
(\ref{wn-c})-(\ref{vpnn}). Then there exists $T>0$ sufficiently small, depending only on $M_0$, $\delta_0$, $c_1$, $c_2$, and satisfying \eqref{T-as1} so that the following holds. 
For every $k\geq 1$ and  a.e.~$\tau\in [0,T]$, one has
\bel{wk-b2}
\left|w^{(k)}(\tau,0\pm)-\overline{w}(0\pm)\right|~\leq~\delta_0,
\eeq
\bel{dotsi} \left|\dot{\sigma}^{(k)}(\tau)\right|~\leq~4C_1 (1+M_0)\cdot|\ln \tau|,
\eeq
\bel{wk-H2}
 \left\|w^{(k)}(\tau,\cdot)\right\|_{H^2(\R\backslash \{0\})}~\leq~M_0\,,
\eeq
for some constant $C_1>0$.
\end{lemma}
{\bf Proof.} 
{\bf 1.} It is clear that  (\ref{wk-b2})-(\ref{wk-H2}) hold for $k=1$. 
By induction, 
assume that  (\ref{wk-b2}) holds for a given $k\geq 1$. 
By  the assumptions (\ref{ic}) and (\ref{wk-b2}), for all $\tau\in [0,T]$ one obtains
\begin{eqnarray*}
\sigma^{(k)}(\tau)&\geq&\overline{w}(0-)-\overline{w}(0+)- \left|w^{(k)}(\tau,0+)-\overline{w}(0+)\right|-\left|w^{(k)}(\tau,0-)-\overline{w}(0-)\right|~\geq~4\delta_0.
\end{eqnarray*}
For a fixed $\tau\in [0,T]$, let $x^{\pm}:[0,\tau]\mapsto\R$ be  the characteristics which reach the origin at time $\tau$,
from  the left and the right, respectively.
Recalling (\ref{w-k+1}), (\ref{F1}), (\ref{ic}), and (\ref{x-b1}), we estimate 
$$\bega{l}\ds
\left| w^{(k+1)}(\tau,0\pm)-\overline{w}(0\pm)\right|~\leq~\left|\overline{w}(x^\pm(0))-\overline{w}(0\pm)\right|+\int_{0}^{\tau}\left|F^{(k)}(t,x^{\pm}(t))\right|dt\\[4mm]
\qquad\ds \leq~3M_0\delta_0\tau+C_1\cdot \int_{0}^{\tau}\left(\ds (1+M_0)\cdot  |\ln t|+ {|\dot{\sigma}^{(k)}(t)|\over |\sigma^{(k)}(t)|}\cdot \big|x^{\pm}(t)|\right)dt\\[4mm]
\qquad\ds\leq~3M_0\delta_0\tau+C_1(1+M_0)\cdot \int_{0}^{\tau}\left(\ds 1+ 5C_1 (\tau-t) \right)\cdot |\ln(t)|dt \\[4mm]
\qquad\ds \leq~\O(1)\cdot (1+M_0)\cdot \big|\ln(\tau)\big|\cdot \tau~ \leq~\O(1)\cdot (1+M_0)\cdot \big|\ln(T)\big|\cdot T
\enda
$$
and this  shows that $w^{(k+1)}$ satisfies (\ref{wk-b2}),  provided that  $T>0$ is chosen sufficiently small, depending only on $M_0,\delta_0,c_1$, and $c_2$.
\medskip

{\bf 2}.  For any $\tau\in [0,T]$ and  $-\delta_1\leq \bar{x}_2<\bar{x}_1<0$,  
consider the characteristics 
$$t\,\mapsto\, x_1(t)\,=\,x( t; \tau, \bar x_1),
\qquad\qquad t\,\mapsto\, x_2(t)\,=\,x( t; \tau, \bar x_2).$$ 
Using (\ref{w-k+1}), (\ref{Lipx}), (\ref{F1}), \eqref{wk-b2}, and (\ref{x-b1}), we estimate 
$$\bega{l}\ds
\left|w^{(k+1)}(\tau,\bar{x}_2)-w^{(k+1)}(\tau,\bar{x}_1)\right|\\[3mm]
\ds\qquad\leq~\bigl|\ov w(x_2(0))-\ov w(x_1(0))\bigr|
+\int_{0}^{\tau}\Big|F^{(k)}(t,x_2(t))-F^{(k)}(t,x_1(t))\Big|dt
\\[3mm]
\ds\qquad \leq~M_0K_1\cdot |\bar{x}_2-\bar{x}_1|+C_1\cdot \int_{0}^{\tau}\left(\ds (1+M_0)\cdot |x_1(t)|^{-1/4}+{|\dot{\sigma}^{(k)}(t)|\over \sigma^{(k)}(t)}\right)\cdot |x_2(t)-x_1(t)|dt
\\[3mm]
\ds\qquad \leq~~M_0K_1\cdot\left(1+\O(1)\cdot\left(1+{1\over M_0}\right)\cdot\left[\left({\tau\over \delta_0}\right)^{1/4}+{|\tau\ln\tau|\over \delta_0}\right]\right)\cdot |\bar{x}_2-\bar{x}_1|.
\enda $$
Therefore, choosing $T>0$ sufficiently small, we obtain
\bel{wx}
\left|w^{(k+1)}_x(\tau,x)\right|~\leq~3M_0K_1\qquad\forall \tau\in [0,T], x\in [-\delta_1,0[\,.
\eeq
An entirely similar  estimate holds for $\tau\in [0,T], x\in \,]0,\delta_1]$.
\v
\n{\bf 3.} Next, given any $0\leq \tau_1<\tau_2\leq T$, denote by  $t\mapsto x_i^{\pm}(t) \doteq x(t;\tau_i,0\pm)$ the characteristic which reaches the origin at time $\tau_i$, from the positive or negative side,
respectively. Recalling (\ref{Lipx})--(\ref{wk-b2}), (\ref{F1}), and (\ref{wx}), we estimate 
$$\bega{l}\ds
\left|w^{(k+1)}(\tau_2,0\pm)-w^{(k+1)}(\tau_1,0\pm)\right|\\[3mm]
\qquad\ds \leq~\left|w^{(k+1)}\big(\tau_1,x^{\pm}_2(\tau_1)\big)-w^{(k+1)}(\tau_1,0\pm)\right|
+\int_{\tau_1}^{\tau_2}\left|F^{(k)}(t,x_2^{\pm}(t))\right|dt\\[3mm]
\qquad\ds \leq~3M_0K_1 \left|x_2^{\pm}(\tau_1)\right|+C_1\cdot \int_{\tau_1}^{\tau_2}\left(\ds (1+M_0)\cdot  |\ln t|+ {|\dot{\sigma}^{(k)}(t)|\over |\sigma^{(k)}(t)|}\cdot \big|x_2^{\pm}(t)|\right)dt\\[3mm]
\qquad\ds \leq ~15M_0K_1\delta_0(\tau_2-\tau_1)+C_1(1+M_0)  \int_{\tau_1}^{\tau_2}
 |\ln t|(1+5C_1(\tau_2-t))dt\\[3mm]
\qquad\ds \leq~\left(15M_0K_1\delta_0+C_1(1+M_0) |\ln(\tau_1)|\right)\cdot (\tau_2-\tau_1)~\leq~2C_1(1+M_0) |\ln(\tau_1)|\cdot (\tau_2-\tau_1),
\enda $$
provided that  $T>0$ is sufficiently small. In particular, we have  
\begin{eqnarray*}
\left|\sigma^{(k+1)}(\tau_2)-\sigma^{(k+1)}(\tau_1)\right|&\leq& 4C_1(1+M_0)|\ln \tau_1|\cdot (\tau_2-\tau_1).
\end{eqnarray*}
This shows that $\dot\sigma^{(k+1)}$ satisfies (\ref{dotsi}).
\v
{\bf 4.} Finally, from Lemma \ref{ln-b1}, Lemma \ref{F1-g}, (\ref{dotsi}), and  Duhamel's formula, 
 for all $\tau\in [0,T]$ we obtain 
\begin{multline*}
\left\|w^{(k+1)}(\tau,\cdot)\right\|_{H^2(\R\backslash \{0\})}~\leq~{3\over 2}\cdot \|\overline{w}\|_{H^2
\bigl(\R\backslash [-\delta_0\tau,\delta_0\tau]\bigr)}+{3\over 2}\cdot\int_{0}^{\tau}\left\|F^{(k)}(t,\cdot)\right\|_{H^2\bigl(\R\backslash [-\delta_0(\tau-t),\delta_0(\tau-t)]\bigr)}dt\\
~\leq~{3M_0\over 4}+{3\over 2}\cdot C_1 (1+M_0)\cdot \int_{0}^{\tau}\left(1+{C_1\over\delta_0}\cdot \big|\ln t\big|\right)\cdot \delta_0^{-2/3}\cdot (\tau-t)^{-2/3}+|\ln t|dt\\
~\leq~{3M_0\over 4}+{3\over 2}\cdot C_1 (1+M_0)\cdot \left[6\left(1+{C_1\over\delta_0}\cdot \big|\ln \tau\big|\right)\delta_0^{-2/3}\tau^{1/3}+\big|\tau\ln\tau\big|\right]~\leq~M_0\,,
\end{multline*}
provided that $T>0$ is sufficiently small, depending only on $M_0,\delta_0,c_1$, and $c_2$. This shows that (\ref{wk-H2}) is satisfied by $w^{(k+1)}$  as well.
\endproof

Thanks to the above estimates, we can now prove that the sequence of approximations
$w^{(k)}$ is Cauchy and converges to a solution $w$ of the linear problem (\ref{wne}).  
This is a key step toward the proof of Theorem~\ref{t:21}.

\begin{lemma}\label{l:44} Let $w_n$, $\vp_n$ be as in 
(\ref{wn-c})-(\ref{vpnn}). Then, for some $T>0$ sufficiently small,  depending only on $M_0$, $\delta_0$, $c_1$, $c_2$, such that,
the sequence of approximations $\left(w^{(k)}\right)_{k\geq 1}$ converges  to 
a limit function $w$ in ${\bf L}^{\infty}\bigl([0,T],H^2(\R\backslash \{0\})\bigr)$, i.e.,
\[
\lim_{k\to\infty}\sup_{t\in [0,T]}\left\|w^{(k)}(t,\cdot)-w(t,\cdot)\right\|_{H^2(\R\backslash \{0\})}~=~0.
\]
The function $w$ provides a solution to  the Cauchy  problem (\ref{wne}) and satisfies 
\bel{cond-wn}
\left|w(\tau,0\pm)-\overline{w}(0\pm)\right|~\leq~\delta_0,\qquad \|w(t,\cdot)\|_{H^2(\R\backslash\{0\})}~\leq~M_0,\qquad t\in [0,T]
\eeq
 Moreover, $\sigma(t)\doteq  w(t,0-)-w(t,0+)$ is locally Lipschitz in $(0,T)$  and 
 \bel{cond-wn-1}
 \left|\dot{\sigma}(t)\right|~\leq~4C_1 (1+M_0)\cdot|\ln t|,\qquad a.e.~t\in (0,T).
\eeq
\end{lemma}

{\bf Proof.} {\bf 1.} For any $k\geq 1$, we set 
\bel{set1}
\begin{cases}
\quad z^{(k)}&\doteq~w^{(k+1)}-w^{(k)},\qquad \sigma_z^{(k)}(t)~\doteq~z^{(k)}(t,0-)-z^{(k)}(t,0+),\\[4mm]
 M^{(k)}_z(t)&\doteq~\left\|z^{(k)}(t,\cdot)\right\|_{H^2(\R\backslash\{0\})},\quad \beta_k(\tau)~\doteq~\ds\sup_{t\in [0,\tau]} M^{(k)}_{z}(t),\quad \alpha_k(\tau)~\doteq~\ds\sup_{t\in [0,\tau]}\left|\dot{\sigma}_z^{(k)}(t)\right|.
\end{cases}
\eeq
Recalling Lemma \ref{FF}, Lemma \ref{ln-b1}, and Lemma \ref{wk-p}, and  using  Duhamel's formula, we estimate 
$$\bega{l}
\ds M^{(k+1)}_{z}(\tau)~\leq~{3\over 2}\cdot \int_{0}^{\tau}\left\|F^{(k+1)}(t,\cdot)-F^{(k)}(t,\cdot)\right\|_{H^2(\R\backslash [-\delta_0(\tau-t),\delta_0(\tau-t)])}dt\\[4mm]
\leq~\ds C_3\cdot\int_{0}^{\tau}\beta_k(t)\cdot|\ln(t)|\left(1+{1\over (\tau-t)^{2/3}}\right)+{\alpha_k(t)\over \sqrt{\tau-t}}~dt~\leq~C_4\cdot \left(\beta_k(\tau) \tau^{1/3}+\alpha_k(\tau)\tau^{1/2}\right),\enda $$
and this implies 
\bel{z-k1}
\beta_{k+1}(\tau)~\leq~C_4\cdot\left(\beta_k(\tau)\tau^{1/3}+\alpha_k(\tau)\tau^{1/2}\right),\qquad 
 \tau\in [0,T],~ k\geq 1,
\eeq
for some constant $C_3,C_4>0$ depending only on $M_0,\delta_0,c_1$, and $c_2$.
\medskip

{\bf 2}. We now establish a bound on $\left\|\dot{\sigma}^{(k+1)}\right\|_{{\bf L}^{\infty}([0,T])}$. Given any $0<\tau_1<\tau_2\leq T$, denote by  $t\mapsto x_i^{\pm}(t)\doteq x(t;\tau_i,0\pm)$  the characteristics, which reach the origin  at time $\tau_i$, from the positive or negative side, respectively. Using (\ref{F1-F2-0}), (\ref{dotsi}), and (\ref{x-b1}) we obtain
$$\bega{l}\ds
\left|z^{(k+1)}(\tau_2,0\pm)-z^{(k+1)}(\tau_1,0\pm)\right|\\[3mm]
\quad \ds \leq~\left|z^{(k+1)}\big(\tau_1,x^{\pm}_2(\tau_1)\big)-z^{(k+1)}(\tau_1,0\pm)\right|+\int_{\tau_1}^{\tau_2}\left|F^{(k+1)}\left(t,x_2^{\pm}(t)\right)-F^{(k)}\left(t,x_2^{\pm}(t)\right)\right|dt\\[4mm]
\quad\qquad\qquad\qquad \ds\leq~\beta_{k+1}\cdot\left|x^{\pm}_2(\tau_1)\right|+C_5\cdot \int_{\tau_1}^{\tau_2}\alpha_k(t)\cdot \big|x^{\pm}_2(t)\big|+\beta_k(t)\cdot\big|\ln t\big|~dt\\[4mm]
\qquad\qquad\quad\ds \leq~\left(5\beta_{k+1}(\tau_2)\delta_0+C_6\cdot\left[\beta_k(\tau_2)\cdot\big|\ln\tau_1\big|+\alpha_k(\tau_2)(\tau_2-\tau_1)\big|\right]\right)\cdot (\tau_2-\tau_1),
\enda
$$
for some constants $C_5>0$ and $C_6>0$ depending only on $M_0,\delta_0,c_1$, and $c_2$. Thus, for $0<T<\delta_0$ sufficiently small, we obtain
\bel{sgm1}
\alpha_{k+1}(\tau)~\leq~\ds 10\beta_{k+1}(\tau)\delta_0+2C_5\beta_k(\tau) \big|\ln \tau\big|\qquad \tau\in ]0,T], k\geq 1,
\eeq
and  (\ref{z-k1}) yields
\[
\beta_{k+1}(\tau)~\leq~C_7\cdot \left(\beta_{k}(\tau)+\beta_{k-1}(\tau)\right)\cdot \tau^{1/3}.
\]
for some $C_7>0$ depending only on $M_0,\delta_0,c_1$, and $c_2$. In particular, for $T>0$ sufficiently small, one has 
\[
\beta_{k+1}(\tau)+{1\over 2} \beta_{k}(\tau)~\leq~{3\over 4}\cdot \left(\beta_{k}(\tau)+{1\over 2} \beta_{k-1}(\tau)\right),
\]
which implies
\[
\sum_{k=1}^{\infty}\sup_{\tau\in (0,T]}\left\|z^{(k)}(\tau,\cdot)\right\|_{H^2(\R\backslash\{0\})}~<~\infty.
\]
We thus conclude that  $\left(w^{(k)}\right)_{k\geq 1}$ is a Cauchy sequence  in ${\bf L}^{\infty}([0,T],{\bf H}^{2}(\R\backslash\{0\}))$ and converges to a limit function $w\in {\bf L}^{\infty}([0,T],{\bf H}^{2}(\R\backslash\{0\}))$, which provides the solution to the linear problem (\ref{wne}), and satisfies (\ref{cond-wn}). Moreover, since $\ds\lim_{k\to\infty}w^{(k)}(\tau,0\pm)=w(\tau,0\pm)$, one has that $\ds\lim_{k\to\infty}\sigma^{(k)}(\tau)=\sigma(\tau)$ for all $\tau\in [0,T]$. Thus, from (\ref{dotsi}),  $\sigma(\cdot)$ is locally Lipscthitz in $(0,T)$  and satisfies  (\ref{cond-wn-1}).
\endproof
\medskip

We are now ready to complete the proof of  our first  main result.

{\bf Proof of Theorem \ref{t:21}.} As outlined at the end of Section~\ref{sec:2}, we construct, by induction, a sequence of approximate solutions
$(w_n)_{n\geq 1}$ where each $w_{n}$ is the solution to  the linear  problem (\ref{wne}). 
For some $T>0$ small enough, depending only on $M_0$, $\delta_0$, $c_1$, and $c_2$, we claim that 
\bel{sum1}
\sum_{n\geq 2}\,\bigl\|w_{n}(\tau)-w_{n-1}(\tau)\bigr\|_{H^1(\R\setminus\{0\})}~<~\infty\qquad\qquad\forall t\in [0,T]\,.
\eeq
For a fixed $n\geq 2$, we define
\bel{Set1}
\begin{cases}
W_n\,\doteq\,w_{n}-w_{n-1},\quad a_n(t,x)\,\doteq\, a(t,x,w_{n}),\quad A_n(t,x)\,\doteq\,a_{n}(t,x)-a_{n-1}(t,x),\\[3mm]
v_n\,\doteq \varphi^{(w_n)}-\phi(\cdot,0),\quad V_n\,\doteq v_n-v_{n-1},\quad \beta_n(t)\doteq \sup_{s\in [0,t]}\|W_n\|_{H^1(\R\backslash\{0\})}.\\[3mm]
\end{cases}
\eeq
Set $Z_n=W_n+V_n$. From the above definitions, by (\ref{wne}) it follows
\bel{Z-e}
Z_{n+1,t}+a_{n}\cdot Z_{n+1,x}~=~-(A_n w_{n,x}+A_{n+1}v_{n+1,x})+G_{n+1}
\eeq
with 
\[
G_{n+1}~\doteq~B^{(w_{n+1})}-B^{(w_n)}+{\bf H}\left[\varphi^{(w_{n+1})}-\varphi^{(w_{n})}\right]-\left(\varphi^{(w_{n+1})}-\varphi^{(w_{n})}\right)\cdot \phi_x(\cdot,0).
\]
Recalling the first inequality in (\ref{F1-F2-1}) and (\ref{cond-wn}), we estimate 
$$
\begin{cases}
\ds \bigl\|A_n(t,\cdot)\cdot w_{n,x}(t,\cdot)\bigr\|_{H^1\bigl(\R\backslash [-\delta,\delta]\bigr)}
\ds ~\leq~~C_7\cdot M_n(t)~\leq~C_7\cdot\beta_n(t),\\[3mm]
\ds \bigl\|A_{n+1}(t,\cdot)\cdot v_{n+1,x}(t,\cdot)\bigr\|_{H^1\bigl(\R\backslash [-\delta,\delta]\bigr)}
\ds ~\leq~~C_7\cdot {M_{n+1}(t)\over \delta^{1/2}}~\leq~C_7\cdot{\beta_{n+1}(t)\over \delta^{1/2}},\\[3mm]
\ds \|G_{n+1}(t,\cdot)\|_{H^1(\R\backslash [-\delta,\delta])}~\leq~\ds C_7\cdot {M_{n+1}(t)\over \delta^{1/2}}~\leq~\ds C_7\cdot {\beta_{n+1}(t)\over \delta^{1/2}}.
\end{cases}
$$
for some constant $C_7>0$ depending only on $M_0$, $\delta_0$, $c_1$, and $c_2$.  Hence, choosing $T>0$ sufficiently small,  we have, using Duhamel's formula, 
\begin{multline}\label{Z-e1}
\left\|Z_{n+1}(\tau,\cdot)\right\|_{H^1(\R\backslash\{0\})}~\leq~{3\over 2}\cdot \int_{0}^{\tau}\left\|G_{n+1}-A_n w_{n,x}-A_{n+1}v_{n+1,x}\right\|_{H^1(\R\backslash [-\delta_0(\tau-t),\delta_0(\tau-t)])}dt\\
~\leq~{3C_8\over 2}\cdot\int_{0}^{\tau}\beta_n(t)+{\beta_{n+1}(t)\over (\tau-t)^{1/2}}~dt~=~{3C_8\over 2}\cdot\left(\beta_n(\tau)\cdot \tau+2\beta_{n+1}(\tau)\tau^{1/2}\right)
\end{multline}
for some constant $C_8>0$. On the other hand, (\ref{vp1}), (\ref{Set1}), and (\ref{con:vpg}) imply
\[
\|V_{n+1}(\tau,\cdot)\|_{H^1(\R\backslash\{0\})}~\leq~C_9\cdot \|W_{n+1}(\tau,\cdot)\|_{H^1(\R\backslash\{0\})}\cdot \tau^{1/4},
\]
and (\ref{Z-e1}) yields
\[
\|W_{n+1}(\tau,\cdot)\|_{H^1(\R\backslash\{0\})}~\leq~{3C_8\over 2} \cdot\left(\beta_n(\tau)\cdot \tau+2\beta_{n+1}(\tau)\tau^{1/2}\right)+C_9\cdot \|W_{n+1}(\tau,\cdot)\|_{H^1(\R\backslash\{0\})}\cdot \tau^{1/4}.
\]
for some constant $C_9>0$ depending only on $M_0$, $\delta_0$, $c_1$, and $c_2$. In particular, for $T>0$ sufficiently small, one has that
\[
\beta_{n+1}(\tau)~\leq~{1\over 2}\cdot\beta_{n}(\tau)\qquad\forall \tau\in [0,T].
\]
Thus,  (\ref{sum1}) holds and for every $t\in [0,T]$ the sequence of approximations $w_n(t,\cdot)$ is Cauchy in the space $H^1\bigl(\R\backslash\{0\}\bigr)$, and hence it converges to a unique limit $w(t,\cdot)$. 

It remains to check that this limit function $w$ is an entropic solution, i.e., it satisfies, cf. \eqref{BH-1}, \eqref{usw}, and \eqref{id2},
\[
\left(w+\vp^{(w)}\right)(t_0,x_0)=(\ov w+ \ov \vp)(x(0;t_0,x_0))+\int_0^{t_0} {\bf H}\left[w+\vp^{(w)}\right](t, x(t;t_0,x_0))dt,
\]
where $t\mapsto x(t;t_0,x_0)$ is the characteristics curve, obtained by solving (\ref{code}). This follows from slightly rewriting \eqref{wne}, which yields 
\begin{align*}
\left(w_{n+1}+\vp^{(w_{n+1})}\right)(t_0,x_0)& =(\ov w+ \ov \vp)(x_n(0;t_0,x_0))+ \int_0^{t_0} {\bf H}\left[w_{n+1}+\vp^{(w_{n+1})}\right](t,x_n(t; t_0, x_0))dt \\
& -\int_0^{t_0} \left(Z_{n+1}-\frac{W_{n+1}^-(t)+W_{n+1}^+(t)}{2}\right) \vp_x^{(w_{n+1})}(t, x_n(t; t_0, x_0))dt,
\end{align*}
where $t\mapsto x_n(t;t_0, x_0)$ denotes the characteristic curve, obtained by solving (\ref{char}).
\medskip

 Finally, to prove uniqueness, assume that $\Tilde{w},w$ are two entropic solutions.
We then define
\[
W\,\doteq\,\Tilde{w}-w,\qquad  \beta(\tau)~\doteq~\ds\sup_{s\in [0,\tau]}\|W(s,\cdot)\|_{H^1(\R\backslash\{0\})}.
\]
The arguments used in the previous steps now yield the inequality
\begin{eqnarray*}
\beta(\tau)~\leq~{1\over 2}\cdot \beta(\tau)\qquad\forall \tau \in [0,T],
\end{eqnarray*}
which implies $\beta(\tau)=0$ for all $\tau\in [0,T]$ and completes the proof.
\endproof

\section{Two interacting shocks}
\label{sec:5}
\setcounter{equation}{0}
In this section, denote by $u(t,x)$ the solution to Burgers' equation
\bel{Burg}
u_t + uu_x~=~0\,,\qquad \qquad u(0,x)\,=\, \bar u(x),\eeq
and by $v(t,x)$ the solution to the perturbed linearized equation
\bel{BP1}
v_t + u v_x~=~\bfH[u(t,\cdot)](x), \qquad \qquad v(0,x)\,=\, \bar v(x).\eeq
By the method of characteristics, at all points where $u$ is continuous, one has 
\bel{BP2} v\bigl(\tau, y\bigr)~=~\bar v\bigl(y- \tau u(\tau,y)\bigr)+\int_0^\tau 
\bfH[ u(t,\cdot)] (y-(\tau-t)\, u(\tau,y))\, dt.\eeq
We 
expect that $v$ can provide a leading order correction term, in an 
{\it ansatz} describing the solution with two interacting shocks to the Burgers-Hilbert equation (\ref{BH}).
\medskip

To fix the ideas, consider a piecewise constant solution to Burgers' equation
containing two interacting shocks with initial data
\bel{idb}
u(0,x)~=~\bar u(x)~=~\left\{ \bega{cl} u^{\ell}\qquad &\hbox{if} ~~x<\bar x_1\,,\\[2mm]
u^m\qquad &\hbox{if} ~~\bar x_1<x<\bar x_2\,,\\[2mm]
u^r\qquad &\hbox{if} ~~\bar x_2<x\,.\enda\right.\eeq
with $u^\ell>u^m>u^r$. Setting
$$\left\{ \bega{rl}\sigma_1&= ~u^{\ell}-u^m,\\[2mm]
 \sigma_2 &=~ u^m-u^r,\enda\right. \qquad\qquad \left\{ \bega{rl} a_1&=~\ds {u^\ell+u^m\over 2}\,,
 \\[3mm]
 a_2&=~ \ds{u^m+u^r\over 2}\,,\enda\right.$$
 
$$x_1(t) \,=\, \bar x_1 + a_1 t,\qquad
\qquad x_2(t) \,=\, \bar x_2 +a_2 t,$$
we thus have
\bel{pcu}
u(t,x)~=~\left\{ \bega{cl} u^\ell\qquad &\hbox{if} ~~x<x_1(t)\,,\\[2mm]
u^m\qquad &\hbox{if} ~~ x_1(t)<x< x_2(t)\,,\\[2mm]
u^r\qquad &\hbox{if} ~~x_2(t)<x\,.\enda\right.\eeq

\begin{figure}[htbp]
\centering
  \includegraphics[scale=0.5]{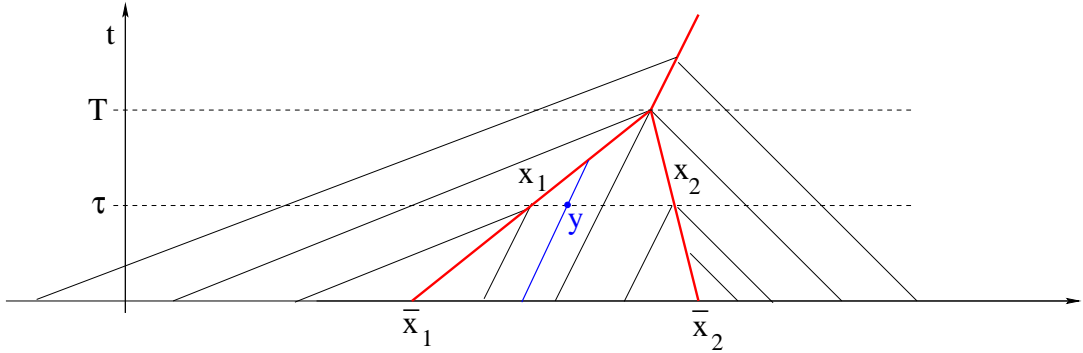}
    \caption{\small The characteristics for a solution to
    Burgers' equation with two shocks at $x_1(t)<x_2(t)$.  }
\label{f:bh30}
\end{figure}

We now compute the corresponding solution $v= v(\tau,y)$ of (\ref{BP1}).
For this purpose, consider the characteristic through the point $(\tau,y)$, namely
\bel{cty} x(t)~= ~y + (t-\tau) u(\tau, y).\eeq   
Recalling (\ref{HT}), 
we  compute  the integral
\bel{HPC}\bega{rl}\ds
I(\tau,y)&\doteq~ \ds\int_0^\tau 
\bfH[ u(t,\cdot)] (x(t))\, dt~=~
-{\sigma_1\over\pi} \int_0^\tau \ln|x_1(t)-x(t)|\, dt -{\sigma_2\over\pi} \int_0^\tau \ln|x_2(t)-x(t)|\, dt \\[4mm]
&=\ds ~{\sigma_1\over 2|\dot x_1-\dot x|} \Big( \vp(|x_1(\tau) -y|) -
\vp\bigl(|x_1(0)-x(0)|\bigr)\Big)  \\[4mm]
&\qquad\qquad \ds+{\sigma_2\over2 |\dot x_2-\dot x|} \Big( \vp(|x_2(\tau) -y|) -
\vp\bigl(|x_2(0)-x(0)|\bigr)\Big) + e(t,x).\enda
\eeq 
Here $\vp$ is given by (\ref{vpdef}), while $e=e(t,x)$ is an additional  smooth
correction term. Neglecting smooth terms, we thus consider three cases,
depending on the location of the characteristic 
$x(t)$ w.r.t.~the two shocks:
\v
CASE 1: $y<x_1(\tau)$.   We then have
\bel{C1} I(\tau,y)~\approx~{2\over\pi} \left[(x_1(\tau)-y)\ln(x_1(\tau)-y) 
+ {\sigma_2\over 2\sigma_1+ \sigma_2 } (x_2(\tau)-y)\ln(x_2(\tau)-y)\right]
.
\eeq

CASE 2: $x_1(\tau)<y<x_2(\tau)$.   We then have
\bel{C2} I(\tau,y)~\approx~{2\over\pi} \left[
(y-x_1(\tau))\ln(y-x_1(\tau)) 
+  (x_2(\tau)-y)\ln(x_2(\tau)-y)\right].
\eeq

CASE 3: $x_2(\tau)<y$.   We then have
\bel{C3} I(\tau,y)~\approx~{2\over\pi} \left[{\sigma_1\over \sigma_1+2\sigma_2} (y-x_1(\tau))\ln(y-x_1(\tau)) 
+  (y-x_2(\tau))\ln(y-x_2(\tau))\right].
\eeq

Next, consider a more general piecewise smooth 
solution $u$  of the Burgers--Hilbert equation (\ref{BH}) with two interacting shocks located at 
points $x_1(\tau)<x_2(\tau)$ with strengths
\bel{sigt}\left\{\bega{rl}\sigma_1(\tau)&=~u\bigl(\tau, x_1(\tau)-\bigr) -u\bigl(\tau, x_1(\tau)+\bigr),\\[2mm]
\sigma_2(\tau)&=~u\bigl(\tau, x_2(\tau)-\bigr) -u\bigl(\tau, x_2(\tau)+\bigr),
\enda\right.\eeq
respectively.
As the interaction time $T$ is approached, we expect that the two 
limits will coincide
$$\lim_{\tau\to T-} u(\tau, x_1(\tau)+) ~=~\lim_{\tau\to T-} u(\tau, x_2(\tau)-).$$
Furthermore, as shown in Fig.~\ref{f:bh30}, 
all characteristics located in the triangular region between the two shocks
will hit one of them within time $T$. 

To construct such solutions, we should thus try with an {\it ansatz}
of the form
\bel{an1}
u(\tau,y)~=~w(\tau,y) + \phi(\tau,y),\eeq
where 
\bel{wh2}
w(\tau,\cdot)~\in~ H^2\Big( ]-\infty, x_1(\tau)[~\cup~ ]x_1(\tau), x_2(\tau)[~
\cup~]x_2(\tau), +\infty[\Big).\eeq
Moreover, in view of (\ref{C1})--(\ref{C3}),  the correction $\phi$ should be defined as
{\small 
\bel{phidef}
\phi(\tau,y)\,=\,\left\{ \bega{l} \ds{2\over\pi} \left[
(x_1(\tau)-y) \ln (x_1(\tau)-y) + {\sigma_2(\tau)\over 2\sigma_1(\tau)+\sigma_2(\tau)} (x_2(\tau)-y) \ln (x_2(\tau)-y) 
\right]
\quad \hbox{if} ~~y<x_1(\tau),\\[4mm]
\ds{2\over\pi} \Big[
(y-x_1(\tau)) \ln (y-x_1(\tau)) + (x_2(\tau)-y) \ln (x_2(\tau)-y) 
\Big]
\qquad \quad \hbox{if} ~~x_1(\tau)<y<x_2(\tau),\\[4mm]
\ds {2\over\pi} \left[{\sigma_1(\tau)\over \sigma_1(\tau)+2\sigma_2(\tau)}
(y-x_1(\tau)) \ln (y-x_1(\tau)) +  (y-x_2(\tau)) \ln (y-x_2(\tau)) 
\right]
\quad \hbox{if} ~~x_2(\tau)<y.
\enda\right.
\eeq
}
\n Note that, for each fixed time $\tau<T$, since $x_1(\tau)<x_2(\tau)$, 
for $y<x_1(\tau) $, the term $\ln(x_2(\tau)-y)$ remains smooth.  The same 
is true for the term $\ln(y-x_1(\tau))$ in the region where $y> x_2(\tau)$.
As a consequence, the asymptotic profile of the function $\phi(\tau,\cdot)$
near both points $x_1(\tau)$ and $x_2(\tau)$ has the same ``$x\ln |x|$" singularity
that we encountered before.   However, these two additional terms cannot
be removed from the definition of $\phi$, because they are not uniformly smooth 
as $\tau\to T-$.

\v
\section{Constructing a solution with two interacting shocks}
\label{sec:6}
\begin{figure}[htbp]
\centering
  \includegraphics[scale=0.5]{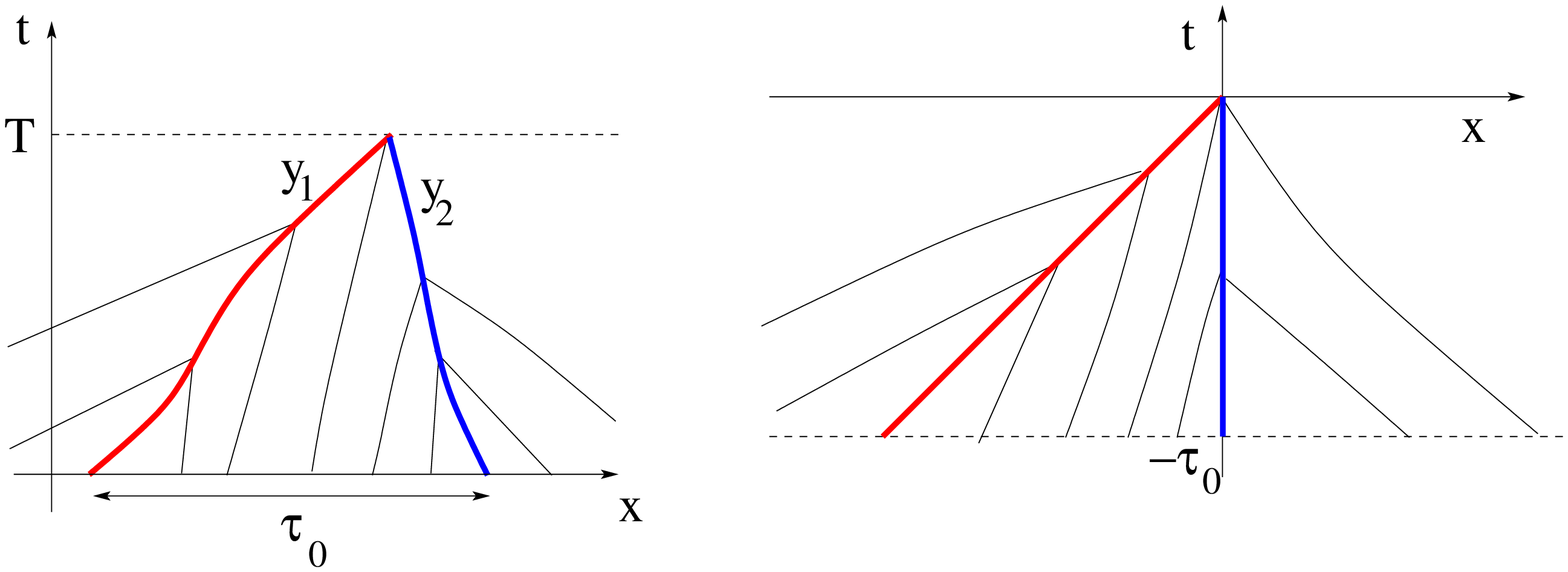}
    \caption{\small Positions of the two shocks in the original variables (left), and in the 
    adapted variables (right).  }
\label{f:bh34}
\end{figure}

\setcounter{equation}{0}
We consider here a solution of the Burgers-Hilbert equation 
(\ref{BH}), which is piecewise continuous  and which has two shocks located at  the points 
$y_1(t)<y_2(t)$. By  the Rankine-Hugoniot conditions, the time derivatives satisfy
\bel{Rankine}
\dot{y}_{i}(t)~=~{u^{-}_{i}(t)+u^{+}_i(t)\over 2},\qquad\qquad i=1,2.
\eeq
Here $u^{\pm}_i(t)\doteq u_i\bigl(t,y_i(t)\pm\bigr)$ denote 
the left and the right limits of $u(t,x)$ as $x\to y_i(t)$.  Throughout the following, we assume that 
\[
{u^{-}_{1}(t)+u^{+}_1(t)\over 2}~=~\dot{y}_1(t)~>~\dot{y}_2(t)~=~{u^{-}_{2}(t)+u^{+}_2(t)\over 2}\,.
\]
The function $\tau\doteq y_1-y_2$ is negative and monotone increasing.
It will be useful to change the space and the time variables, 
so that in the new variables $\tilde t, \tilde x$ the location of one shock is fixed, while the other moves with constant speed $1$.
For this purpose, we set 
$$\tilde x~\doteq~ x-y_2(t),\qquad \qquad \tilde t ~\doteq ~\tau(t)~<~0\,.$$
As a consequence, the two shocks, in the new coordinate system, are located at
$$y_1(\tilde t) \,=\, \tilde t,\qquad\qquad y_2(\tilde t) = 0 ,$$
and interact at the point $(\tilde t, \tilde x) = (0,0)$. Introducing the function
\bel{vde}
v\bigl(\tau(t),x\bigr)~=~u\bigl(t,x+y_2(t)\bigr),
\eeq
we define the left and right values
\bel{v12} v_1^{\pm}\bigl(\tau(t)\bigr)~\doteq~v\bigl(\tau(t),\tau(t)\pm\bigr)~=~u\bigl(t,y_1(t)\pm\bigr),
\qquad 
v_2^{\pm}\bigl(\tau(t)\bigr)~\doteq~v\bigl(\tau(t),0\pm\bigr)~=~u\bigl(t,y_2(t)\pm\bigr).\eeq
The change of variables (\ref{vde}) yields
\[
v_x(\tau,x)~=~u_x\bigl(t,x+y_2(t)\bigr),\qquad v_{\tau}(\tau,x)~=~{u_t\bigl(t,x+y_2(t)\bigr)+\dot{y}_2(t)\cdot u_x\bigl(t,x+y_2(t)\bigr)\over \dot{y}_1(t)-\dot{y}_2(t)}.
\]
Therefore, (\ref{BH}) implies 
\bel{veq}
v_{\tau}\bigl(\tau(t),x\bigr)+ {\big[v\bigl(\tau(t),x\bigr)-\dot{y}_2(t)\big]\cdot v_{x}\bigl(\tau(t),x\bigr)\over \dot{y}_1(t)-\dot{y}_2(t)}~=~{{\bf H}[v\bigl(\tau(t),\cdot\bigr)](x)\over \dot{y}_1(t)-\dot{y}_2(t)}\,.
\eeq
Thus, by (\ref{Rankine}), (\ref{v12}), and (\ref{veq}), 
we can recast the original equation (\ref{BH}) in the following equivalent form 
\bel{BH11}
u_t+{1\over a_1(t)-a_2(t)}\cdot \big[u-a_2(t)\big]\cdot u_x~=~{{\bf H}[u] \over  a_1(t)-a_2(t)}\,.
\eeq
Given $\tau_0<0$, for $t\in [\tau_0,0]$ the two  functions 
\bel{a12}
a_1(t)~\doteq~{u^-_1(t)+u^+_1(t)\over 2},\qquad\qquad a_2(t)~\doteq~{u^-_2(t)+u^+_2(t)\over 2},
\eeq
yield the speeds of the two shocks in the original coordinates, as shown in Fig.~\ref{f:bh34}.

We shall construct the solution of (\ref{BH11}) in  the form 
\bel{form-u}
u(t,x)~=~w(t,x)+\varphi(t,x)\qquad\forall (t,x)\in [\tau_0, 0]\,\times \R\,.
\eeq
Here $\varphi$ is a continuous function, which satisfies $\varphi(t,t)=\varphi(t,0)=0$, while 
\bel{H}w(t,\cdot )~\in~
H^2\Big(]-\infty,t[~\cup~]t,0[~\cup~ ]0, +\infty[ \Big)\qquad\forall t\in [\tau_0,0].
\eeq
According to (\ref{H}),  the function $w(t,\cdot)$ is continuously differentiable outside the two points $x= t$ and $x=0$.  Moreover, the distributional derivative $D_x w(t,\cdot)$ 
is an ${\bf L}^2$ function restricted to each interval $]-\infty,t[\,$, $]t,0[\,$ and $]0,+\infty[\,$.  However, both $w(t,\cdot)$ and $w_x(t,\cdot)$ can have a jump at $x=t$ and at $x=0$.  
At the points $(t,t)$ and $(t,0)$, the following traces are well defined:
\bel{traces-0}
\left\{\bega{cl} w^{-}_1(t)~\doteq~w(t,t-)~=~u^{-}_1(t),\cr\cr
w^{+}_1(t)~\doteq~w(t,t+)~=~u^{+}_1(t),
\enda\right.
 \qquad\qquad\qquad\qquad
 \left\{\bega{cl} b^-_{1}(t)~\doteq~w_x(t,t-),\cr\cr
b^+_{1}(t)~\doteq~w_x(t,t+),
 \enda\right.
 \eeq
 \bel{traces-00}
\left\{\bega{cl} w^{-}_2(t)~\doteq~w(t,0-)~=~u^{-}_2(t),\cr\cr
w^{+}_2(t)~\doteq~w(t,0+)~=~u^{+}_2(t),
\enda\right.
\qquad\qquad\qquad\qquad
 \left\{\bega{cl} b^-_{2}(t)~\doteq~w_x(t,0-),\cr\cr
b^+_{2}(t)~\doteq~w_x(t,0+).
 \enda\right.
 \eeq
For the shocks to be entropy admissible, the inequalities
 \bel{eac}
 w_1^-(t)~>~w_1^+(t),\qquad\qquad w_2^-(t)~>~w_2^+(t),
 \eeq
 will always be assumed.
  Writing  
\bel{ave-w}
a^{(w)}_1(t)~\doteq~{w^-_1(t)+w^+_1(t)\over 2}\qquad\qquad a^{(w)}_2(t)~\doteq~{w^-_2(t)+w^+_2(t)\over 2}\,,
\eeq
the equation (\ref{BH11}) reads
\bel{BH-3}
w_t+a(t,x,w)\cdot w_x~=~F(t,x,w),
\eeq
where $a$ and $F$ are given by
\bel{a33}
a(t,x,w)~=~{w(t,x)+\varphi^{(w)}(t,x)-a^{(w)}_2(t)\over a^{(w)}_1(t)-a^{(w)}_2(t)}\,,
\eeq
\bel{F3}
F(t,x,w)~=~\left[{{\bf H}[w]-(w-a^{(w)}_2(t))\cdot \varphi^{(w)}_x\over a^{(w)}_1(t)-a^{(w)}_2(t)}-\varphi^{(w)}_t\right]+{{\bf H}[\varphi^{(w)}]-\varphi^{(w)}\varphi^{(w)}_x\over a^{(w)}_1(t)-a^{(w)}_2(t)},
\eeq
respectively.
Here the function $\vp^{(w)}(t,x) $ is chosen in such a way that a cancellation between leading order terms near to the location of the two shocks at $x=t$ and at $x=0$ is achieved. More precisely, in view of  (\ref{phidef}) and recalling  (\ref{eta}) and (\ref{phi-b}), we set 
\bel{b-phi-2s}
\phi_0(x)~=~\phi(x,0)~=~{2\eta(x)\over\pi}\cdot |x|\ln|x|,
\eeq
and define 
\bel{cordef}
\vp^{(w)}(t,x)\,=\, \left\{ \bega{cl} \ds 
\phi_0(x-t) + {\sigma^{(w)}_2(t)\over 2\sigma^{(w)}_1(t)+\sigma^{(w)}_2(t)}\cdot\big(\phi_0(x)- \phi_0(t)\big) 
\quad &\hbox{if} ~~x<t,\\[6mm]
\ds \phi_0(x-t)  + \phi_0(x)- \phi_0(t)
\qquad \quad& \hbox{if} ~~t<x<0,\\[3mm]
\ds {\sigma^{(w)}_1(t)\over \sigma^{(w)}_1(t)+2\sigma^{(w)}_2(t)}\cdot
\big(\phi_0(x-t)-\phi_0(t)\big) +  \phi_0(x)
\quad &\hbox{if} ~~0<x\,.
\enda\right.
\eeq

The following theorem provides the existence of solutions to the Cauchy problem for (\ref{BH})
where the initial datum contains two shocks. In particular, the solution to (\ref{BH11}) is constructed up to the time where
the two shocks interact. Furthermore, the solution is of the form 
(\ref{form-u}), where $\vp=\vp^{(w)} $, the corrector function defined in (\ref{cordef}).

\begin{theorem}\label{t:61} For any given constants  $b$, $M_0$, $\delta_1$, $\delta_2>0$, there exists 
$\ve_0>0$ small enough and a constant $K$ such that the following holds.

Consider any $\tau_0\in [-\ve_0, 0[\,$ and any initial condition $\overline{w}\in H^2(\R\backslash\{\tau_0,0\})$ such that 
\bel{condw-w}\left\{
\bega{rll}
\|\overline{w}\|_{H^2(\R\backslash \{\tau_0,0\})}&\leq~\ds{M_0\over 4},\qquad  &\left\|\overline{w}_x\right\|_{{\bf L}^{\infty}(]\tau_0,0[)}~\leq~b,\\[4mm]
\overline{w}(\tau_0-)-\overline{w}(\tau_0+)&\geq~8\delta_{1},\qquad &\overline{w}(0-)-\overline{w}(0+)~\geq~8\delta_2.
\enda\right.
\eeq
Then the Cauchy problem (\ref{BH-3}) with initial data
\bel{id7}
w(\tau_0, \cdot)~=~\ov w~\in~H^2\bigl(\R\setminus \{\tau_0,0\}\bigr)\eeq admits a unique entropic solution, defined for $t\in [\tau_0,0]$. 
Moreover,  this solution satisfies 
\bel{soln-2s-e}
\left\{
\bega{rll}
\|w(t,\cdot)\|_{H^2(\R\backslash \{t,0\})}&\leq~M_0,\qquad  &\left\|w_x(t,\cdot)\right\|_{{\bf L}^{\infty}(]\tau_0,0[)}~\leq~Kb,\\[4mm]
w(t,t-)-w(t,t+)&\geq~\delta_{1},\qquad & w(t,0-)-w(t,0+)~\geq~\delta_2
\enda\right.
\eeq
for all $t\in [\tau_0,0[$.
\end{theorem}
\begin{remark} {\rm By (\ref{soln-2s-e}), at the interaction time $t=0$ 
the solution $u=w+\vp^{(w)}$ is the sum of 
a corrector term plus a function in $H^2\bigl(\R\setminus\{0\}\bigr)$.  
This function lies within the class of initial data
covered by our earlier Theorem~\ref{t:21}. Thus, combining Theorems~\ref{t:61} and \ref{t:21} yields the behavior of a solution to (\ref{BH}) across the interaction of two shocks.}
\end{remark}

Toward a proof of Theorem~\ref{t:61},
 solutions to (\ref{BH-3}) will be constructed by an iteration procedure. The main difference between this and the earlier case with a single shock 
is that the correction term $\vp$ now depends on time through the variable strengths $\sigma_1,\sigma_2$ of the two shocks.  Define 
\bel{init}
w^{(1)}(t,x)~=~\begin{cases}
\overline{w}(x)&\qquad\mathrm{if}~~x\in ]t,0[\cup ]0,\infty[,\\[2mm]
\overline{w}(x+\tau_0-t)&\qquad\mathrm{if}~~x\in ]-\infty,t[.
\end{cases}
\eeq
By induction, let  $w^{(n)}$ be given and satisfy (\ref{H}) for every $t\in [\tau_0, 0[\,$.
Moreover, call $\sigma_1^{(n)}(t)$ and $\sigma_2^{(n)}(t)$ the strengths of the 
two shocks at $x=t$ and $x=0$ of $w^{(n)}$, respectively. We construct the next iterate $w= w^{(n+1)}(t,x)$ by solving the linear equation
\bel{BH4}
w_t+a\left(t,x,w^{(n)}\right) \, w_x~=~F(t,x,w),
\eeq
with initial data (\ref{id7}) and $a$ as introduced in \eqref{a33}.

The induction argument requires the following steps:
\begi
\item [(i)] Given $w^{(n)}$, the equation  (\ref{BH4}) with the initial data 
$\ov w$ admits a unique solution $w$ with $w(t,\cdot)\in H^2\bigl(\R\setminus \{t,0\}\bigr)$ for all
$t\in [\tau_0,0[\,$.  
\item [(ii)] A priori bounds on the strong norm $\|w^{(n)}(t,\cdot)\|_{H^2\bigl(\R\setminus \{t,0\}\bigr)}$ for all $t\in [\tau_0,0[$, $n\geq 1$.
\item [(iii)] Convergence in a weak norm. This will follow from the bound
\[
\sum_{n\geq 1}\|w^{(n+1)}(t,\cdot)-w^{(n)}(t,\cdot)\|_{H^2\bigl(\R\setminus \{t,0\}\bigr)}~<~+\infty.
\]
\endi

\subsection{Some preliminary estimates}
\label{sec:6-s1}
To achieve the above steps (i)-(iii), we first establish some key estimates on the right hand side of (\ref{BH-3}). For any $w:[\tau_0,0]\times\R\to\R$ 
such that $w(t,\cdot)\in H^2(\R\backslash\{t,0\})$ for all $t\in [\tau_0,0]$, we write  
\bel{split-w}
w(t,\cdot)~=~v_1(t,\cdot)+v_2(t,\cdot)\eeq
with 
\bel{w2}
v_2(t,x)~\doteq~\begin{cases}  \left(w(t,0-)+x\cdot w_x(t,0-)\right)\cdot\eta(x)&\mathrm{if}~~~ x<0\,,\\[3mm]
 \left(w(t,0+)+x\cdot w_x(t,0+)\right)\cdot\eta(x)&\mathrm{if}~~~ x>0.
\end{cases}
\eeq
Recalling (\ref{traces-0})-(\ref{ave-w}) and (\ref{F3}), we split $F$ into the four parts:
\bel{F-dec}
F(t,x,w)~=~\ds{1\over a^{(w)}_1(t)-a^{(w)}_2(t)}\cdot \Big[{A^{(w)}(t,x)+B^{(w)}(t,x)+C^{(w)}(t,x)+D^{(w)}(t,x)}\Big].
\eeq
Here we take
\bel{AB-2s}
A^{(w)}~\doteq~{\bf H}\left[\vp^{(w)}\right]-\vp^{(w)}\vp_x^{(w)},\quad B^{(w)}~=~\Tilde{B}^{(v_2)}~\doteq~{\bf H}[v_2]-\left(v_2-{v_2(t,0-)+v_2(t,0+)\over 2}\right)\cdot \phi'_0(x)\,,
\eeq
\bel{C-2s}
C^{(w)}~=~\Tilde{C}^{(v_1)}~\doteq~{\bf H}[v_1]-\left(v_1-{v_1(t,t-)+v_1(t,t+)\over 2}\right)\cdot \phi'_0(x-t),
\eeq
\bel{D-2s}
D^{(w)}~\doteq~{\bf H}[w]-\left(w-a^{(w)}_2(t)\right)\cdot \varphi^{(w)}_x-\left(a^{(w)}_1(t)-a^{(w)}_2(t)\right)\cdot \varphi^{(w)}_t-B^{(w)}-C^{(w)}.
\eeq
\begin{lemma}\label{F1-g-2s} Let $w: [\tau_0,0]\times\R\to\R$ be such that  $w(t,\cdot)\in H^2(\R\backslash\{t,0\})$ for all $t\in[\tau_0,0]$,  and 
\[
\|w(t,\cdot)\|_{H^2(\R\backslash\{t,0\})}~\leq~M_0,\qquad |w_x(t,0-)|+|w_x(t,0+)|~\leq~b,\qquad \sigma_j^{(w)}(t)~\geq~\delta_0~>~0.
\]
Moreover, assume that  $\ds-{1\over 4}\cdot\min\left\{1,{\delta^2_0\over M_0^2}\right\}<\tau_0<0$ and $\sigma_{j}^{(w)}(\cdot)$ is locally Lipschitz on $[\tau_0,0[$ for $j=1,2$. Then there is a constant $C_1>0$, depending only on $M_0$, $b$, $\delta_0$ such that,  for a.e.~$t\in [\tau_0,0]$ and $|x|<{1\over 2e}$, one has 
\bel{F1-2s}
\begin{cases}
|F(t,x,w)|~\leq~C_1\cdot \Bigg[\ds {1+M_0+b\over \delta_0}+\big|\ln |t|\big |\\ \qquad\qquad\qquad\qquad\qquad\qquad\qquad\ds 
+~\chi_{\R^\backslash]t,0[}\cdot{\left|\dot{\sigma}_1^{(w)}(t)\right|+\left|\dot{\sigma}_2^{(w)}(t)\right|\over \delta_0} \cdot \left(|x|^{1/2}+|x-t|^{1/2}\right)\Bigg],\\[4mm]
|F_x(t,x,w)|~\leq~C_1\cdot \ds {1+M_0+b+\chi_{\R^\backslash]t,0[}\cdot\left(\left|\dot{\sigma}_1^{(w)}(t)\right|+\left|\dot{\sigma}_2^{(w)}(t)\right|\right)\over \delta_0}\cdot \left(|x|^{-1/2}+|x-t|^{-1/2}\right).
\end{cases}
\eeq
Furthermore, for every $\delta>0$ sufficiently small one has for all $t\in [\tau_0,0]$ 
\bel{F2-2s}
\|F(t,x,w)\|_{H^2(\R\backslash  [-\delta,\delta]\cup [t-\delta,t+\delta])}~\leq~\ds C_1\cdot\left[{1+ M_0+b+\left|\dot{\sigma}_1^{(w)}(t)\right|+\left|\dot{\sigma}_2^{(w)}(t)\right|\over\delta_0}\cdot \ds \delta^{-2/3}+|\ln (t)|\right].
\eeq
\end{lemma}

{\bf Proof.} We  observe that, for all $t\in [\tau_0,0]$, it holds
\bel{a2-a1}
a_1^{(w)}(t)-a_2^{(w)}(t)~\geq~\delta_0+w(t,t+)-w(t,0-)~\geq~\delta_0-M_0\cdot |t|^{1/2}~\geq~{\delta_0\over 2},
\eeq
\[ 
0<\sigma^{(w)}_j(t)<\sigma_1^{(w)}(t)+\sigma_2^{(w)}(t)=w(t,t-)-w(t,t+)+w(t,0-)-w(t,0+)\leq (2+\sqrt{\vert t\vert})M_0,
\]
\[
\|v_1(t,\cdot)\|_{H^2(\R\backslash\{t\})}, \|v_2(t,\cdot)\|_{H^2(\R\backslash\{0\})}~\leq~ \O(1)\cdot (M_0+b).
\]
According to (\ref{F-dec}), the function $F$ can be decomposed as the sum of four terms,
which will be estimated separately.
\medskip

{\bf 1.} Recalling \eqref{AB} and \eqref{Bb1}, \eqref{AB-2s} and \eqref{C-2s} imply that for every $(t,x)\in [\tau_0,0]\times \left[-{1\over 2e},{1\over 2e}\right]$ one has
 \bel{Bb1-2s}
\begin{cases}
\left|B^{(w)}(t,x)\right|~\leq~\O(1)\cdot (M_0+b),\qquad \left|B^{(w)}_x(t,x)\right|~\leq~ \O(1)\cdot (M_0+b)\cdot \big|\ln|x|\big|^2,\\[3mm]
 \left\|B^{(w)}(t,\cdot)\right\|_{H^2(\R\backslash [-\delta,\delta])}~\leq~\O(1)\cdot (M_0+b)\cdot\delta^{-2/3}.
\end{cases}
 \eeq
and, for every $(t,x)\in [\tau_0,0]\times \left[t-{1\over 2e},t+{1\over 2e}\right]$  
 \bel{Cc1-2s}
\begin{cases}
\left|C^{(w)}(t,x)\right|~\leq~\O(1)\cdot (M_0+b),\qquad \left|C^{(w)}_x(t,x)\right|~\leq~ \O(1)\cdot (M_0+b)\cdot \big|\ln|x-t|\big|^2,\\[3mm]
 \left\|C^{(w)}(t,\cdot)\right\|_{H^2(\R\backslash [t-\delta,t+\delta])}~\leq~\O(1)\cdot (M_0+b)\cdot\delta^{-2/3}.
\end{cases}
 \eeq
{\bf 2.} Next, we estimate $A^{(w)}$. Recalling (\ref{g-b}), i.e.,
\[
g_b(x)~=~\chi_{\strut [0,\infty[}(x)\cdot\phi(x,b),\qquad x\in\R, b\geq 0,
\]
\eqref{phi-b}, and \eqref{b-phi-2s},
we can rewrite, for $(t,x)\in [\tau_0,0]\times \left[-{1\over 2e},{1\over 2e}\right]$, 
\bel{vpd2}
\varphi^{(w)}(t,x)~=~\phi_0(x-t)+\phi_0(x)-\phi_0(t)+E_1^{(w)}(t,x)+E_2^{(w)}(t,x),
\eeq 
\bel{E} 
\begin{cases}
E_1^{(w)}(t,x)&=~- \ds{2\cdot \chi_{]-\infty,t]}\cdot \sigma^{(w)}_1(t)\over 2\sigma^{(w)}_1(t)+\sigma^{(w)}_2(t)}\cdot \big[\phi_0(x)-\phi_0(t)\big]~=~ \ds-{2 \sigma^{(w)}_1(t)\over 2\sigma^{(w)}_1(t)+\sigma^{(w)}_2(t)}\cdot g_{|t|}(t-x), \\[4mm]
E_2^{(w)}(t,x)&=~- \ds{2\cdot \chi_{[0,\infty[}\cdot\sigma^{(w)}_2(t)\over \sigma^{(w)}_1(t)+2\sigma^{(w)}_2(t)}\cdot \big[\phi_0(x-t)-\phi_0(t)\big]~=~ -{2\sigma^{(w)}_2(t)\over \sigma^{(w)}_1(t)+2\sigma^{(w)}_2(t)}\cdot g_{|t|}(x).
\end{cases}
\eeq
Thus, for $(t,x)\in [\tau_0,0]\times \left[-{1\over 2e},{1\over 2e}\right]$, Lemma \ref{b-phi} and \cite[Section 3]{BZ} imply, \bel{Hvp-s2}
\left\{\bega{rlrl}
&\big|{\bf H}\big[\varphi^{(w)}(t,\cdot)\big](x)\big|\ds~\leq~\ds\O(1),\quad  \ds \left|{d\over dx}{\bf H}\big[\varphi^{(w)}(t,\cdot)\big](x)\right|
\ds\leq~\O(1)\cdot \left(\ln^2|x|+\ln^2|x-t|\right),\\[4mm]
& \big\|{\bf H}\big[\varphi^{(w)}(t,\cdot)\big]\big\|_{H^{2}(\R\backslash([-\delta,\delta])\cup [t-\delta,t+\delta])}~\leq~\O(1)\cdot \delta^{-2/3}.
\enda\right.
\eeq
On the other hand, given $t\in [\tau_0,0]$, for every $\vert x\vert \leq \frac{1}{2e}$, \eqref{vpd2} combined with \eqref{e1-gb} yields 
\[ 
\left\vert \varphi^{(w)}(t,x)\right\vert \leq \min(\vert x \ln(\vert x\vert)\vert , \vert (x-t)\ln(\vert x-t\vert)\vert).
\]
Furthermore, we compute for every $x\in \R\backslash \{t,0\}$
\bel{vpx-s2}
\varphi_x^{(w)}(t,x)~=~~\left(1-\ds{2\cdot \chi_{[0,\infty[}\cdot\sigma^{(w)}_2(t)\over \sigma^{(w)}_1(t)+2\sigma^{(w)}_2(t)}\right)\phi'_0(x-t)+\left(1-\ds{2\cdot \chi_{]-\infty,t]}\cdot \sigma^{(w)}_1(t)\over 2\sigma^{(w)}_1(t)+\sigma^{(w)}_2(t)}\right)\cdot\phi'_0(x),
\eeq
which together with \cite[Section 3]{BZ} implies for $\vert x\vert \leq \frac1{2e}$ that 
\[
\left|\varphi_x^{(w)}(t,x)\right|~\leq~\O(1)\cdot\left(\big|\ln |x|\big|+\big|\ln|x-t|\big|\right),\qquad \left|\varphi_{xx}^{(w)}(t,x)\right|~\leq~\O(1)\cdot\left({1\over |x|}+{1\over |x-t|}\right),
\]
\[
\left|\varphi_{xxx}^{(w)}(t,x)\right|~\leq~\O(1)\cdot\left({1\over x^2}+{1\over (x-t)^2}\right).
\]
A direct computation yields, for $\vert x\vert \leq \frac{1}{2e}$, 
\[
\begin{cases}
\left|\vp^{(w)} \varphi_x^{(w)}(t,x)\right|&\leq~\O(1)\cdot \left(|x|^{1/2}+|x-t|^{1/2}\right),\\[2mm]
 \left|\ds{d\over dx}[\vp^{(w)} \varphi_x^{(w)}](t,x)\right|&\leq~\O(1)\cdot \left(\ln^2 |x|+\ln^2|x-t|\right),\\[3mm]
  \left|\ds{d^2\over dx^2}[\vp^{(w)} \varphi_x^{(w)}](t,x)\right|&\leq~\O(1)\cdot \ds\left(\left|{\ln |x|\over x} \right|+\left|{\ln |x-t|\over x-t} \right|\right).
\end{cases}
\]
and thus
\[
 \left\|\vp^{(w)} \varphi_x^{(w)}(t,\cdot)\right\|_{H^2(\R\backslash [t-\delta,t+\delta]\cup [-\delta,\delta])}~\leq~\O(1)\cdot\delta^{-2/3}.
\]
Recalling  (\ref{Hvp-s2}), we get 
 \bel{Aa1-2s}
\begin{cases}
\left|A^{(w)}(t,x)\right|~\leq~\O(1),\qquad \left|A^{(w)}_x(t,x)\right|~\leq~  \O(1)\cdot \left(\ln^2|x|+\ln^2|x-t|\right),\\[2mm]
\left\|A^{(w)}(t,\cdot)\right\|_{H^2(\R\backslash [t-\delta,t+\delta]\cup (-\delta,\delta))}~\leq~\O(1)\cdot\delta^{-2/3}.
\end{cases}
\eeq

{\bf 3.} Finally, to estimate $D^{(w)}$, we shall consider three cases:
\medskip

{\bf Case 1:} Assume that $-{1\over 2}< t-{1\over 2e}<x<t$. We have
\bel{D-sl-case1}
D^{(w)}(t,x)~=~D^{(w)}_1(t,x)+D^{(w)}_2(t,x)+ D_3^{(w)}(t,x)
\eeq
with 
\[
\begin{cases}
D_1^{(w)}(t,x)~=~ w_x(t,0-)\cdot (t-x)\cdot   \phi_0'(x-t),\\[2mm]
D_2^{(w)}(t,x)~=~\ds \ds\left({\sigma^{(w)}_2(t)\cdot[w(t,0-)-w(t,t+)+w(t,t-)-w(t,x)]\over 2\sigma^{(w)}_1(t)+\sigma^{(w)}_2(t)}+xw_{x}(t,0-) \right)\cdot \phi'_0(x),\\[4mm]
D_3^{(w)}(t,x)~=~\big[a^{(w)}_2(t)-a^{(w)}_1(t)\big]\cdot \left[\left(\ds {\sigma^{(w)}_2(t)\over 2\sigma^{(w)}_1(t)+\sigma^{(w)}_2(t)}\right)'\cdot g_{|t|}(t-x)-\ds{\sigma^{(w)}_2(t)\cdot\phi'_0(t)\over 2\sigma^{(w)}_1(t)+\sigma^{(w)}_2(t)}\right].
\end{cases}
\]
Recalling \eqref{b-phi-2s}, we estimate 
\bel{D1}
\begin{cases}
\left|D_1^{(w)}(t,x)\right|~\leq~\O(1)\cdot b\cdot |x-t|\big|\ln |x-t|\big|~\leq~O(1)\cdot b,\\[2mm]
\left|\ds{d\over dx}D_1^{(w)}(t,x)\right|~\leq~\O(1)\cdot b\cdot  \big|\ln|x-t|\big|,\qquad \left|\ds{d^2\over d^2x}D_1^{(w)}(t,x)\right|~\leq~\ds O(1)\cdot {b\over |x-t|},
\end{cases}
\eeq
\bel{D2}
\begin{cases}
\left|D_2^{(w)}(t,x)\right|~\leq~\O(1)\cdot (M_0+b)\cdot x\big|\ln |x|\big|,\quad\left|\ds{d\over dx}D_2^{(w)}(t,x)\right|~\leq~\O(1)\cdot (M_0+b)\cdot  \big|\ln|x|\big|,\\[3mm]
\left|\ds{d^2\over dx^2}D_2^{(w)}(t,x)\right|~\leq~\ds  |w_{xx}\cdot \phi'_0(x)|+\O(1)\cdot {M_0+b\over |x|},
\end{cases}
\eeq
and 
\bel{D3}
\begin{cases}
\left|D_3^{(w)}(t,x)\right|~\leq~\ds\O(1)\cdot \left|a^{(w)}_2(t)-a^{(w)}_1(t)\right|\cdot\left({\left|\dot{\sigma}^{(w)}_1(t)\right|+\left|\dot{\sigma}^{(w)}_2(t)\right|\over \delta_0} |x-t|  |\ln |x-t| |+ \big|\ln |t|\big|\right)\\[4mm]
\left|\ds{d\over dx}D_3^{(w)}(t,x)\right|~\leq~\ds\O(1)\cdot \left|a^{(w)}_2(t)-a^{(w)}_1(t)\right|\cdot{\left|\dot{\sigma}^{(w)}_1(t)\right|+\left|\dot{\sigma}^{(w)}_2(t)\right|\over \delta_0}\cdot  \big|\ln |x-t|\big|,\\[4mm]
\left|\ds{d^2\over dx^2}D_3^{(w)}(t,x)\right|~\leq~\ds\O(1)\cdot \left|a^{(w)}_2(t)-a^{(w)}_1(t)\right|\cdot{\left|\dot{\sigma}^{(w)}_1(t)\right|+\left|\dot{\sigma}^{(w)}_2(t)\right|\over \delta_0}\cdot {1\over |x-t|}\,.
\end{cases}
\eeq
Combining (\ref{D-sl-case1})-(\ref{D3}), we obtain
\bel{D-case1}
\begin{cases}
\ds{\left|D^{(w)}(t,x)\right|\over \left|a^{(w)}_2(t)-a^{(w)}_1(t)\right|}~\leq~\O(1)\cdot \left({M_0+b\over \delta_0}+|\ln (t)|+{\left|\dot{\sigma}^{(w)}_1(t)\right|+\left|\dot{\sigma}^{(w)}_2(t)\right|\over \delta_0}\cdot |x-t|^{1/2}\right),\\[4mm]
\ds{\left|\ds{d\over dx}D^{(w)}(t,x)\right|\over  \left|a^{(w)}_2(t)-a^{(w)}_1(t)\right|}~\leq~\ds\O(1)\cdot {M_0+b+\left|\dot{\sigma}^{(w)}_1(t)\right|+\left|\dot{\sigma}^{(w)}_2(t)\right|\over \delta_0}\cdot  \big|\ln |x-t|\big|,\\[4mm]
\ds{\left|\ds{d^2\over dx^2}D^{(w)}(t,x)\right|\over  \left|a^{(w)}_2(t)-a^{(w)}_1(t)\right|}~\leq~\ds\O(1)\cdot \left({M_0+b+\left|\dot{\sigma}^{(w)}_1(t)\right|+\left|\dot{\sigma}^{(w)}_2(t)\right|\over \delta_0\cdot |x-t|}+ {|w_{xx}|\over\delta_0}\cdot \big|\ln |x-t|\big|\right).
\end{cases}
\eeq

{\bf Case 2:} Assume that $\ds t<x<0$. We have
\bel{D-sl-case2}
D^{(w)}(t,x)~=~w_x(t,0-)\cdot (t-x)\cdot   \phi_0'(x-t)+(v_1(t,0)-v_1(t,x))\cdot \phi'_0(x)+\big[a^{(w)}_1(t)-a^{(w)}_2(t)\big]\cdot \phi_0'(t),
\eeq
and this yields 
\bel{D-case2}
\begin{cases}
\ds{\left|D^{(w)}(t,x)\right|\over  \left|a^{(w)}_2(t)-a^{(w)}_1(t)\right|}~\leq~\O(1)\left({M_0+b\over \delta_0}+\big|\ln|t|\big|\right),\\[4mm]
\ds{\left|\ds{d\over dx}D^{(w)}(t,x)\right|\over  \left|a^{(w)}_2(t)-a^{(w)}_1(t)\right|}~\leq~\O(1)\cdot {M_0+b\over \delta_0}\cdot \left( \big|\ln|x|\big|+\big|\ln|x-t|\big|\right),\\[4mm]
\ds{\left|\ds{d^2\over dx^2}D^{(w)}(t,x)\right|\over  \left|a^{(w)}_2(t)-a^{(w)}_1(t)\right|}~\leq~\ds\O(1)\cdot \left({M_0+b\over\delta_0}\cdot \left({1\over |x|}+{1\over |x-t|}\right)+{\big|w_{xx}(t,x)|\over \delta_0}\cdot \big|\ln |x|\big| \right).
\end{cases}
\eeq

{\bf Case 3:} Assume that $0<x<1/2$. As in Case 1, writing
\bel{D-sl-case3}
D^{(w)}(t,x)~=~D^{(w)}_1(t,x)+D^{(w)}_2(t,x)+ D_3^{(w)}(t,x)
\eeq
with 
\[
\begin{cases}
D_1^{(w)}(t,x)~=~\left[v_1(t,x)-v_1(t,t+) \right] \cdot \phi_0'(x-t)\\
\qquad\qquad+~\left[\ds{\sigma^{(w)}_1(t)\over \sigma^{(w)}_1(t)+2\sigma^{(w)}_2(t)}\cdot (w(t,t+)-w(t,0-)+w(t,0+)-w(t,x))\right]\cdot \phi_0'(x-t),\\[4mm]
D_2^{(w)}(t,x)~=~\ds \ds\left(w(t,0+)-w(t,x)+x w_x(t,0+)\right)\cdot \phi'_0(x),\\[4mm]
D_3^{(w)}(t,x)~=~\big[a^{(w)}_2(t)-a^{(w)}_1(t)\big]\cdot \left[\left(\ds {\sigma^{(w)}_1(t)\over \sigma^{(w)}_1(t)+2\sigma^{(w)}_2(t)}\right)'\cdot g_{|t|}(x)-\ds{\sigma^{(w)}_1(t)\cdot\phi'_0(t)\over \sigma^{(w)}_1(t)+2\sigma^{(w)}_2(t)}\right],
\end{cases}
\]
we estimate
\bel{D-case3}
\begin{cases}
\ds{\left|D^{(w)}(t,x)\right|\over  \left|a^{(w)}_2(t)-a^{(w)}_1(t)\right|}~\leq~\ds\O(1)\cdot \left({M_0+b\over \delta_0}+ |\ln (t)|+{\left|\dot{\sigma}^{(w)}_1(t)\right|+\left|\dot{\sigma}^{(w)}_2(t)\right|\over \delta_0}\cdot |x|^{1/2}\right),\\[4mm]
\ds{\left|\ds{d\over dx}D^{(w)}(t,x)\right|\over  \left|a^{(w)}_2(t)-a^{(w)}_1(t)\right|}~\leq~\ds\O(1)\cdot {M_0+b+\left|\dot{\sigma}^{(w)}_1(t)\right|+\left|\dot{\sigma}^{(w)}_2(t)\right|\over \delta_0}\cdot  \big|\ln |x|\big|,\\[4mm]
\ds{\left|\ds{d^2\over dx^2}D^{(w)}(t,x)\right|\over  \left|a^{(w)}_2(t)-a^{(w)}_1(t)\right|}~\leq~\ds\O(1)\cdot \left({M_0+b+\left|\dot{\sigma}^{(w)}_1(t)\right|+\left|\dot{\sigma}^{(w)}_2(t)\right|\over\delta_0}\cdot  {1\over |x|}+ {|w_{xx}|\over \delta_0}\cdot \big|\ln |x|\big| \right).
\end{cases}
\eeq
In summary, from (\ref{D-case1}), (\ref{D-case2}), and (\ref{D-case3}), given $t\in [\tau_0,0]$, for every $x\in (-1/2,1/2)\backslash\{t,0\}$, it holds that 
\bel{D-g}
\begin{cases}
\ds{\left|D^{(w)}(t,x)\right|\over  \left|a^{(w)}_2(t)-a^{(w)}_1(t)\right|}~\leq~\O(1)\cdot \Bigg({M_0+b\over \delta_0}+ |\ln (t)|\\ \qquad\qquad\qquad\qquad\qquad\qquad\qquad\quad\ds+\chi_{\R\backslash ]0,t[}\cdot{\left|\dot{\sigma}^{(w)}_1(t)\right|+\left|\dot{\sigma}^{(w)}_2(t)\right|\over \delta_0}\cdot \left(|x|^{1/2}+|x-t|^{1/2}\right)\Bigg),\\[4mm]
\ds{\left|\ds{d\over dx}D^{(w)}(t,x)\right|\over  \left|a^{(w)}_2(t)-a^{(w)}_1(t)\right|}\leq\O(1)\cdot {M_0+b+\chi_{\R\backslash]t,0[}\left(\left|\dot{\sigma}^{(w)}_1(t)\right|+\left|\dot{\sigma}^{(w)}_2(t)\right|\right)\over\delta_0}\cdot  \left(\big|\ln |x|\big|+\big|\ln |x-t|\big|\right),\\[4mm]
\ds{\left\|D^{(w)}(t,\cdot)\right\|_{H^2(\R\backslash [t-\delta,t+\delta]\cup (-\delta,\delta))}\over \left|a^{(w)}_2(t)-a^{(w)}_1(t)\right|}~\leq~\ds\O(1)\cdot \left[{M_0+b+\left|\dot{\sigma}^{(w)}_1(t)\right|+\left|\dot{\sigma}^{(w)}_2(t)\right|\over\delta_0 \delta^{2/3}}+\big|\ln |t|\big|\right].
\end{cases}
\eeq
To complete the proof, combining   (\ref{Bb1-2s}), (\ref{Cc1-2s}), (\ref{Aa1-2s}) and (\ref{D-g}), we obtain (\ref{F1-2s})-(\ref{F2-2s}).
\endproof

The next lemma estimates the change in the function $F= F(t,x,w)$ as $w(\cdot)$ 
takes different values.  
These estimates will play a key role in the proof of convergence of the approximations inductively
defined by (\ref{BH4}).
 
\begin{lemma}\label{FF-2s} Let $w_1,w_2: [\tau_0,0]\times\R\to\R$ be such that,  for $i\in\{1,2\}$  and $t\in [\tau_0,0]$,
one has $w_i(t,\cdot)\in H^2(\R\backslash\{t,0\})$ and
\[
\|w_i(t,\cdot)\|_{H^2(\R\backslash\{t,0\})}~\leq~M_0,\qquad |w_{i,x}(t,0-)|+|w_{i,x}(t,0+)|~\leq~b,\qquad \left|\sigma_j^{(w_i)}(t)\right|~\geq~\delta_0. \]
Moreover, assume that $-{1\over 4}\cdot \min\left\{ 1, {\delta_0^2\over M_0^2}\right\}< \tau_0<0$ and $\sigma_j^{(w_i)}$ is  locally Lipschitz on $[\tau_0,0[$ and there exists a function $K(t)$ such that
\[
\max\left\{\left|\dot{\sigma}_1^{(w_i)}(t)\right|, \left|\dot{\sigma}_2^{(w_i)}(t)\right|\right\}~\leq~K(t)\qquad a.e.~t\in [\tau_0,0].
\]
Set  $z\doteq w_2-w_1,~\sigma_i^{(z)}\doteq \sigma_i^{(w_2)}-\sigma_i^{(w_1)}$, and $\gamma^{(z)}(t)~\doteq~\max\left\{\left|\dot{\sigma}_1^{(z)}\right|,\left|\dot{\sigma}_2^{(z)}\right|\right\}$. 
Furthermore, let
\[
M_2(t)~\doteq~ \|z(t,\cdot)\|_{H^2(\R\backslash\{t,0\})}+|z_{x}(t,0-)|+|z_{x}(t,0+)|+\vert z(t,0-)\vert + \vert z(t,0+)\vert.
\]
Then there exists a constant $C_2>0$, depending only on $M_0$, $b$, $\delta_0$ such that, for every $\ds x\in \left[-{1\over 2e},{1\over 2e}\right]$ and a.e. $t\in [\tau_0,0]$, one has
\begin{multline}\label{F1-F2-0-2s}
\bigl| F(t,x,w_2)-F(t,x,w_1)\bigr|~\leq~{C_2\over \delta_0^2}\cdot (M_0+b)\cdot \gamma^{(z)}(t)\cdot \left(|x|^{1/2}\chi_{[0,\infty[}+|x-t|^{1/2}\chi_{[-\infty,t[}\right)\\
\qquad~~~~+{C_2\over \delta_0^2}\cdot \left[M_2(t)\cdot\left(\big|\ln|t|\big|+{M_0+b\over \delta_0}+{K(t)\over \delta_0}\cdot \left(|x|^{1/2}+|x-t|^{1/2}\right)\right)\right]\end{multline}
and for every $x\in (t,0)$
\begin{multline}\label{F1-F2-0-2ss}
\bigl| F_x(t,x,w_2)-F_x(t,x,w_1)\bigr|~\leq {C_2\over \delta_0^2}\cdot M_2(t)\cdot \left(1+\delta_0+M_0+b)(|x|^{-1/2}+|x-t|^{-1/2}\right).
\end{multline}
Moreover, for every $\delta>0$ sufficiently small, it holds 
\begin{multline}\label{F1-F2-H2-2s}
\left\|F(t,\cdot,w_2)-F(t,\cdot,w_1)\right\|_{H^2(\R\backslash [-\delta,\delta]\cup [t-\delta,t+\delta])}~\leq~{C_2\over\delta_0^2}\cdot (M_0+b)\cdot \gamma^{(z)}(t)\cdot \delta^{-2/3}\\
+{C_2\over \delta_0^2}\cdot M_2(t)\cdot \left[\left({K(t)+M_0+b+1\over \delta_0}+M_0+b+1\right)\cdot \delta^{-2/3}+\big|\ln |t|\big|\right].
\end{multline}
\end{lemma}

{\bf Proof.} {\bf 1.} For notational convenience, we set  
\bel{bF-2s}
{\bf A}^{(z)}~\doteq~A^{(w_2)}-A^{(w_1)},~~ {\bf B}^{(z)}~\doteq~B^{(w_2)}-B^{(w_1)},~~ {\bf C}^{(z)}~\doteq~C^{(w_2)}-C^{(w_1)},~~ {\bf D}^{(z)}~\doteq~D^{(w_2)}-D^{(w_1)}.
\eeq
Furthermore, let $z_j=v_{1,j}-v_{2,j}$ for $j=1$, $2$, then 
 \[ \left\|z_1(t,\cdot)\right\|_{H^2(\R\backslash \{t\})}\leq \O(1)\cdot M_2(t) \quad \left\|z_2(t,\cdot)\right\|_{H^2(\R\backslash \{0\})}\leq \O(1)\cdot M_2(t).\]
Comparing \eqref{AB} and \eqref{AB-2s} and recalling \eqref{B2-B1}, then yields, for every $(t,x)\in [\tau_0,0]\times \left[-{1\over 2e},{1\over 2e}\right]$, 
 \bel{B2-B1-2s}
\begin{cases}
\left|{\bf B}^{(z)}(t,x)\right|=\O(1)\cdot M_2(t),\quad \left|{\bf B}^{(z)}_x(t,x)\right|=\O(1)\cdot M_2(t)\cdot \vert \ln|x|\vert ^2, \\[4mm]
 \left\|{\bf B}^{(z)}(t,\cdot)\right\|_{H^2(\R\backslash [-\delta,\delta])}\leq\O(1)\cdot \ds{M_2(t)\over \delta^{2/3}}.
\end{cases}
\eeq
Similarly, for every $(t,x)\in [\tau_0,0]\times \left[t-{1\over 2e},t+{1\over 2e}\right]$, it holds 
 \bel{C2-C1-2s}
\begin{cases}
\left|{\bf C}^{(z)}(t,x)\right|=\O(1)\cdot  M_2(t),\quad \left|{\bf C}^{(z)}_x(t,x)\right|=\O(1)\cdot M_2(t)\cdot \vert \ln|x-t|\vert ^2, \\[4mm]
 \left\|{\bf C}^{(z)}(t,\cdot)\right\|_{H^2(\R\backslash [t-\delta,t+\delta])}\leq\O(1)\cdot \ds{M_2(t)\over \delta^{2/3}}.
\end{cases}
\eeq

{\bf 2.}  We now provide bounds on ${\bf A}^{(z)}(t,x)$. From (\ref{vpd2}) and (\ref{E}), it follows that
\begin{multline}\label{vp-vp}
\varphi^{(w_2)}(t,x)-\varphi^{(w_1)}(t,x)~=~\left[\ds{2 \sigma^{(w_1)}_1(t)\over 2\sigma^{(w_1)}_1(t)+\sigma^{(w_1)}_2(t)}-\ds{2 \sigma^{(w_2)}_1(t)\over 2\sigma^{(w_2)}_1(t)+\sigma^{(w_2)}_2(t)}\right]\cdot g_{|t|}(t-x)\\
+\left[{2\sigma^{(w_1)}_2(t)\over \sigma^{(w_1)}_1(t)+2\sigma^{(w_1)}_2(t)}-{2\sigma^{(w_2)}_2(t)\over \sigma^{(w_2)}_1(t)+2\sigma^{(w_2)}_2(t)}\right]\cdot g_{|t|}(x).
\end{multline}
Since
\[
\begin{cases}
\left|\ds{2 \sigma^{(w_1)}_1(t)\over 2\sigma^{(w_1)}_1(t)+\sigma^{(w_1)}_2(t)}-\ds{2 \sigma^{(w_2)}_1(t)\over 2\sigma^{(w_2)}_1(t)+\sigma^{(w_2)}_2(t)}\right|~\leq~\ds\O(1)\cdot{M_2(t)\over \delta_0}\,,\\[5mm]
\ds\left|{2\sigma^{(w_1)}_2(t)\over \sigma^{(w_1)}_1(t)+2\sigma^{(w_1)}_2(t)}-{2\sigma^{(w_2)}_2(t)\over \sigma^{(w_2)}_1(t)+2\sigma^{(w_2)}_2(t)}\right|~\leq~\O(1)\cdot{M_2(t)\over \delta_0}\,,
\end{cases}
\]
Lemma \ref{b-phi} implies for $x\not \in \{t,0\}$ and $\vert x\vert \leq {1\over 2e}$,
\[
 \begin{cases}
 \left|{\bf H}[\varphi^{(w_2)}(t,x)-\varphi^{(w_1)}(t,x)]\right|~\leq~\ds\O(1)\cdot  {M_2(t)\over \delta_0},\\[3mm]
 \left|\ds{d\over dx}{\bf H}[\varphi^{(w_2)}(t,x)-\varphi^{(w_1)}(t,x)]\right|~\leq~\ds \O(1)\cdot {M_2(t)\over \delta_0}\cdot \left(\ln^2|x|+\ln^2|x-t|\right),\\[3mm]
 \left\|{\bf H}[\varphi^{(w_2)}(t,\cdot)-\varphi^{(w_1)}(t,\cdot)]\right\|_{H^2(\R\backslash [t-\delta,t+\delta]\cup (-\delta,\delta))}~\leq~\ds \O(1)\cdot {M_2(t)\over \delta_0}\cdot\delta^{-2/3}.
 \end{cases}
\]
and from  (\ref{e1-gb}), we obtain for $x\not \in \{t,0\}$ and $\vert x\vert \leq {1\over 2e}$,
\[
 \begin{cases}
 \left|\varphi^{(w_2)}\varphi^{(w_2)}_x-\varphi^{(w_1)}\varphi^{(w_1)}_x\right|~\leq~\ds\O(1)\cdot  {M_2(t)\over \delta_0}\cdot \left(|x|^{1/2}+|x-t|^{1/2}\right),\\[3mm]
\left|\ds{d\over dx}\left(\varphi^{(w_2)}\varphi^{(w_2)}_x-\varphi^{(w_1)}\varphi^{(w_1)}_x\right)\right|~\leq~\ds\O(1)\cdot  {M_2(t)\over \delta_0}\cdot \left(\ln^2|x|+\ln^2|x-t|\right),\\[4mm]
\left|\ds{d^2\over dx^2}\left(\varphi^{(w_2)}\varphi^{(w_2)}_x-\varphi^{(w_1)}\varphi^{(w_1)}_x\right)\right|~\leq~\ds\O(1)\cdot  {M_2(t)\over \delta_0}\cdot \left({|\ln|x||\over |x|}+{|\ln|x-t||\over |x-t|}\right),\\[4mm]
 \left\|\varphi^{(w_2)}\varphi^{(w_2)}_x-\varphi^{(w_1)}\varphi^{(w_1)}_x\right\|_{H^2(\R\backslash [t-\delta,t+\delta]\cup (-\delta,\delta))}~\leq~\ds \O(1)\cdot {M_2(t)\over \delta_0}\cdot\delta^{-2/3}.
 \end{cases}
\]
Thus, (\ref{AB-2s}) yields
 \bel{bA1-2s}
\begin{cases}
\left|{\bf A}^{(z)}(t,x)\right|~\leq~\ds\O(1)\cdot {M_2(t)\over \delta_0},\qquad \left|{\bf A}^{(z)}_x(t,x)\right|~\leq~  \O(1)\cdot  {M_2(t)\over \delta_0}\cdot \left(\ln^2|x|+\ln^2|x-t|\right),\\[3mm]
\left\|{\bf A}^{(z)}(t,\cdot)\right\|_{H^2(\R\backslash [t-\delta,t+\delta]\cup (-\delta,\delta))}~\leq~\ds\O(1)\cdot  {M_2(t)\over \delta_0}\cdot\delta^{-2/3}\,.
\end{cases}
\eeq

{\bf 3.} Finally, to achieve bound on  ${\bf D}^{(z)}$, we consider three cases as in the proof of Lemma \ref{F1-g-2s}. As before, we define ${\bf D}_i^{(z)}=D_i^{(w_2)}-D_i^{(w_1)}$ for $i=1,2,3.$ 

{\bf Case 1:} Assume that $-1/2< t-{1\over 2e}<x<t$. Note, that we can write 
\[
\begin{cases}
{\bf D}^{(z)}(t,x)~=~{\bf D}_1^{(z)}(t,x)+{\bf D}_2^{(z)}(t,x)+{\bf D}_3^{(z)}(t,x)\\[2mm]
{\bf D}_1^{(z)}(t,x)~=~I^{(z)}_{1}(t,x)\cdot \phi_0'(x-t),\qquad {\bf D}_2^{(z)}(t,x)~=~I^{(z)}_2(t,x)\cdot \phi_0'(x)\,,\\[2mm]
{\bf D}_3^{(z)}(t,x)~=~I^{(z)}_{31}(t)\cdot g_{|t|}(t-x)+I^{(z)}_{32}(t)\cdot\phi_0'(t).
\end{cases}
\]
which implies 
  \[
 \begin{cases}
 \left|I^{(z)}_{1}(t,x)\right|~\leq~\O(1)\cdot M_2(t)\cdot |t-x|, \quad \left|\partial_xI^{(z)}_{1}(t,x)\right|~\leq~\O(1)\cdot M_1(t),\\[3mm]
 \partial^2_{xx}I^{(z)}_{1}(t,x)~=~0,\quad \left|I^{(z)}_{2}(t,x)\right|~\leq~\O(1)\cdot M_2(t)\cdot \ds {M_0+b\over \delta_0}\cdot |x|,\\[3mm]
   \left|\partial_x I^{(z)}_{2}(t,x)\right|~\leq~\ds\O(1)\cdot \left({M_0+b\over \delta_0}\cdot M_2(t)+|z_x(t,x)|\right),\\[3mm]  
  \left|\partial^2_{xx} I^{(z)}_{2}(t,x)\right|~\leq~\ds\O(1)\cdot \left({\vert w_{2,xx}(t,x)\vert \over \delta_0}\cdot M_2(t)+|z_{xx}(t,x)|\right),\\[3mm]
   \left|I^{(z)}_{31}(t)\right|~\leq~\ds\O(1)\cdot \left({K(t)(1+\delta_0)  M_2(t)\over \delta_0^2}+{M_0\gamma^{(z)}(t)\over \delta_0}\right),\quad \left|I^{(z)}_{32}(t)\right|~\leq~\O(1)\cdot {M_0 M_2(t) \over \delta_0}\,.
 \end{cases}
 \]
 Thus, for $t>0$ sufficiently small such that $|t|<e^{-M_0-b}$, it holds
 \begin{multline*}
 \left|{\bf D}^{(z)}(t,x)\right|~\leq~\ds\O(1)\cdot \Bigg[{M_0+b\over \delta_0}\cdot\gamma^{(z)}(t)\cdot |x-t|^{1/2}\\
 +{M_2(t)\over \delta_0}\cdot \left({K(t) \over \delta_0}\cdot |x-t|^{1/2}+(1+M_0)\cdot \big|\ln |t|\big|\right)\Bigg],
 \end{multline*}
and
\[
\begin{cases}
\left|\ds{d\over dx}{\bf D}^{(z)}(t,x)\right|~\leq~\ds\O(1)\cdot \left[{M_2(t)\over \delta_0}\cdot \left({K(t)\over \delta_0}+M_0+b+1\right)\cdot |x-t|^{-1/2}+|z_x(t,x)|\cdot\big|\ln|x|\big|\right]\\[4mm]
\qquad\qquad\qquad\qquad\qquad\qquad\qquad\qquad\qquad\qquad+\O(1)\cdot\left(\ds {M_0+b\over \delta_0}\cdot \gamma^{(z)}(t)\cdot\big| \ln|t-x|\big|\right)\\[4mm]
\left|\ds{d^2\over d^2x}{\bf D}^{(z)}(t,x)\right|~\leq~\ds\O(1)\cdot {M_2(t)\over \delta_0}\cdot\left[\left({K(t)\over \delta_0}+M_0+b+1\right)\cdot {1\over |x-t|} +|w_{2xx}(t,x)|\big|\ln|x|\big|\right]\\[4mm]
\qquad\qquad\qquad \qquad\qquad\qquad\ds+\O(1)\cdot\left(|z_{xx}(t,x)|\cdot\big|\ln|x|\big|+{M_0+b\over \delta_0}\cdot \gamma^{(z)}(t)\cdot {1\over |x-t|}\right).
\end{cases}
\]

{\bf Case 2:} Assuming $\ds t<x<0$, we have
\[
{\bf D}^{(z)}(t,x)~=~z_x(t,0-)\cdot (t-x)\cdot   \phi_0'(x-t)+(z_1(t,0)-z_1(t,x))\cdot \phi'_0(x)+\left[a^{(z)}_1(t)-a^{(z)}_2(t)\right]\cdot \phi_0'(t),
\]
with $$a^{(z)}_j(t)\doteq a^{(w_2)}_j(t)-a^{(w_1)}_j(t),$$
 which yields
\[
\begin{cases}
\ds{\left|{\bf D}^{(z)}(t,x)\right|}~\leq~\O(1)\cdot M_2(t)\cdot\big|\ln|t|\big|,\\[4mm]
\ds{\left|{d\over dx}{\bf D}^{(z)}(t,x)\right|}~\leq~\O(1)\cdot \left[M_2(t)\cdot \left(|x|^{-1/2}+|x-t|^{-1/2}\right)+ |z_{1,x}(t,x)|\cdot \big|\ln\vert x\vert\big|\right],\\[4mm]
\ds{\left|{d^2\over dx^2}{\bf D}^{(z)}(t,x)\right|}~\leq~\O(1)\cdot\left[M_2(t)\cdot \left({1\over |x|}+{1\over |x-t|}\right)+|z_{1,xx}(t,x)|\cdot \big|\ln\vert x\vert\big| \right].
\end{cases}
\]

{\bf Case 3:} Assume that $0<x<1/2$ and $\vert t\vert \leq e^{-M_0-b}$. As in Case 1, we estimate 
\begin{multline*}
\left|{\bf D}^{(z)}(t,x)\right|\\
~\leq~\O(1)\cdot \left[{M_2(t)\over \delta_0}\cdot \left({K(t)(M_0+\delta_0)\cdot x^{1/2}\over \delta_0}+(1+M_0)\big|\ln |t|\big|\right)+{M_0+b\over \delta_0}\cdot\gamma^{(z)}(t)\cdot x^{1/2}\right]
\end{multline*}
and
\[
\begin{cases}
\left|\ds{d\over dx}{\bf D}^{(z)}(t,x)\right|~\leq~\ds\O(1)\cdot \left[{M_2(t)\over \delta_0} \left({K(t)(M_0+\delta_0)\over \delta_0}+M_0+b+1\right)\cdot |x|^{-1/2}\right]\\[4mm]
\qquad\qquad\qquad\qquad\qquad\qquad\ds+\O(1)\cdot\left(|z_x(t,x)|+|z_{1,x}(t,x)|+{M_0+b\over \delta_0}\cdot \gamma^{(z)}(t)\right)\cdot\big| \ln  x \big|,\\[4mm]
\left|\ds{d^2\over dx^2}{\bf D}^{(z)}(t,x)\right|~\leq~ \O(1)\cdot\left[\ds\big(|z_{xx}(t,x)|+|z_{1,xx}(t,x)|\big)\cdot \big|\ln \vert x\vert\big|+{M_0+b\over \delta_0}\cdot \gamma^{(z)}(t)\cdot {1\over |x|} \right]\\[4mm]
\qquad\qquad+\ds\O(1)\cdot {M_2(t)\over \delta_0}\cdot\left[ \left({K(t)(M_0+\delta_0)\over \delta_0}+M_0+b+1\right)\cdot {1\over |x|}+|w_{2,xx}(t,x)|\big|\ln|x-t|\big|\right].
\end{cases}
\]
In summary, given $t\in [\tau_0,0]$ sufficiently small, for every $x\in (-1/2,1/2)\backslash\{t,0\}$, it holds that 
\bel{bD-2s}
\begin{cases}
\big|{\bf D}^{(z)}(t,x)\big|~\leq~\ds\O(1)\cdot \Big[{M_2(t)\over \delta_0}\cdot(1+M_0) \cdot\big|\ln |t|\big|\\[4mm]
\qquad \qquad \quad  +\ds\left({M_2(t)K(t)\over \delta^2_0}+M_2(t)+{M_0+b\over \delta_0}\gamma^{(z)}(t)\right) \cdot \left(|x|^{1/2}\chi_{[0,\infty[}+|x-t|^{1/2}\chi_{]-\infty,t]}\right)\Big],\\[4mm]
\left|\ds{d\over dx}{\bf D}^{(z)}(t,x)\right|~\leq~\O(1)\cdot {M_0+b\over \delta_0}\cdot \gamma^{(z)}(t)\cdot \left(|x|^{-1/4}\chi_{[0,\infty[}+|x-t|^{-1/4}\chi_{]-\infty, t]}\right)\\[4mm]\qquad\qquad~~+\ds\O(1)\cdot \left[\left({M_2(t)\over \delta_0} \left({K(t)\over \delta_0}+M_0+b+1\right)+|z_x(t,x)|\right)\cdot \left(|x|^{-1/2}+|x-t|^{-1/2}\right)\right]\\[4mm]
\big\|{\bf D}^{(z)}(t,\cdot)\big\|_{H^2(\R\backslash [t-\delta,t+\delta]\cup (-\delta,\delta))}~\leq~\ds\O(1)\cdot {M_2(t)\over \delta_0}\cdot \left[\left({K(t)\over \delta_0}+M_0+b+1\right)\cdot\delta^{-2/3}+\big|\ln|t|\big|\right]\\[4mm]
\qquad\qquad\qquad\qquad\qquad\qquad\qquad\qquad\qquad\qquad\qquad\qquad\quad\ds+\O(1)\cdot{M_0+b\over \delta_0}\cdot \gamma^{(z)}(t)\cdot \delta^{-2/3}.
\end{cases}
\eeq
Finally, combining (\ref{F-dec})-(\ref{D-2s}), Lemma \ref{F1-g-2s}, and (\ref{bF-2s})-(\ref{bD-2s}), we obtain (\ref{F1-F2-0-2s})-(\ref{F1-F2-H2-2s})
\endproof

\subsection{Proof of Theorem~\ref{t:61}}
\label{sec:6-s2}
We are now ready to give a proof of Theorem~\ref{t:61}. Given $\tau_0\in [-\ve_0,0[\,$ sufficiently small and some initial data $w(\tau_0,\cdot)=\overline{w}$ satisfying (\ref{condw-w}), we construct a solution to the  Cauchy problem (\ref{BH-3}). This solution will be obtained
as the limit of a Cauchy sequence of approximate solutions $w^{(n)}(t,x)$, following the steps (i)--(iii) 
outlined in the beginning of Section~\ref{sec:6}. 
\medskip

{\bf Step 1.} Let $b$, $M_0$, $\delta_1$, $\delta_2>0$ and $\overline{w}\in  H^2(\R\backslash \{\tau_0,0\})$ such that
\bel{ic-1-2s}
\begin{cases}
\|\overline{w}\|_{H^2(\R\backslash \{\tau_0,0\})}&\leq~\ds{M_0\over 4},\qquad \left\|\overline{w}_x\right\|_{{\bf L}^{\infty}(]t,0[)}~\leq~b,\\[4mm]
\overline{w}(\tau_0-)-\overline{w}(\tau_0+)&=~8\delta_{1},\qquad \overline{w}(0-)-\overline{w}(0+)~=~8\delta_2.
\end{cases}
\eeq
We first establish the existence and uniqueness of solutions to the linear 
problem (\ref{BH4}) with initial data $\overline w$ and a given function $w^{(n)}$ with $w^{(n)}(t,\cdot)\in H^2(\R\backslash \{t,0\})$ for all $t\in [\tau_0,0[$ and such that for all $ t\in [\tau_0,0[$,
\bel{wn-c-1}
\begin{cases}
\left\|w^{(n)}(t,\cdot)\right\|_{H^2(\R\backslash\{t,0\})}&\leq~M_0,\qquad \left\|w_x^{(n)}(t,x)\right\|_{{\bf L}^{\infty}(]t,0[)}~\leq~Kb,\\[4mm]
\left|w^{(n)}(t,t\pm)-\overline{w}(\tau_0\pm)\right|&\leq~\delta_1,\qquad \left|w^{(n)}(t,0\pm)-\overline{w}(0\pm)\right|~\leq~\delta_2\,,
\end{cases}
\eeq
for some constant $K>0$ depending only on $b$, $M_0$, $\delta_1$, $\delta_2$. Note that $w^{(1)}$, defined in \eqref{init}, satisfies all of these assumptions.

Note that if such a sequence exist, then the constant $\delta_0$ in Lemma~\ref{F1-g-2s} and Lemma~\ref{FF-2s} can be chosen as $\min(\delta_1, \delta_2)$. Accordingly, we define 
$$ \delta_0~\doteq~\min(\delta_1, \delta_2).$$ 
Assume
\bel{Tau-as1}
-{1\over 4}\cdot{\min\{\delta^2_1,\delta_2^2\}\over M^2_0}~<~\tau_0~<~ 0,
\eeq
and denote by  $t\mapsto x(t;t_0,x_0)$ the solution to the Cauchy problem 
\bel{char-s2}
\dot{x}(t)~=~a_n(t,x(t)),\qquad x(t_0)~=~x_0,
\eeq
where
\bel{a-n}
a_n(t,x)~\doteq~{w^{(n)}(t,x)+\varphi^{(n)}(t,x)-a^{(n)}_2(t)\over a^{(n)}_1(t)-a^{(n)}_2(t)}.
\eeq
Here, 
\bel{a-n-2}
\varphi^{(n)}(t,x)~\doteq~\varphi^{\left(w^{(n)}\right)}(t,x), \quad a^{(n)}_j(t)=a^{(w^{(n)})}_j(t), \quad \sigma_j^{(n)}(t)=\sigma^{(w^{(n)})}_j(t)\quad \text{ for }j=1,2.
\eeq
To begin with we study the travel direction of $x(t)$, which depends on the sign of $a_n$. Therefore observe that  (\ref{ic-1-2s}) and (\ref{wn-c-1}) imply 
\bel{sigma12}
6\delta_i~\leq ~\sigma^{(n)}_i(t)~\leq ~10\delta_i\qquad t\in [\tau_0,0[, i\in \{1,2\}.
\eeq
Furthermore,
\[
 \left|a^{(n)}_1(t)-a^{(n)}_2(t)-{1\over 2}\cdot\left(\sigma^{(n)}_1(t)+\sigma^{(n)}_2(t)\right)\right|~=~\left|w^{(n)}(t,t+)-w^{(n)}(t,0-)\right|~\leq~M_0\sqrt{|t|},
 \]
and 
\[
3(\delta_1+\delta_2)-M_0\sqrt{\vert t\vert }~\leq~\left|a^{(n)}_1(t)-a^{(n)}_2(t)\right|~\leq~5(\delta_1+\delta_2)+M_0\sqrt{\vert t\vert}.
\]
Recalling \eqref{Tau-as1} we end up with
\bel{a1-a2}
2(\delta_1+\delta_2)~\leq~\left|a^{(n)}_1(t)-a^{(n)}_2(t)\right|~\leq~6(\delta_1+\delta_2).
\eeq

For every $(t,x)\in [\tau_0,0[\times ]0,{1\over 2}[$, one has, using \eqref{cordef}, \eqref{E}, and \eqref{e1-gb},
\begin{multline}\label{a1}
\left|a_n(t,x)+{\sigma^{(n)}_2(t)\over 2( a^{(n)}_1(t)-a^{(n)}_2(t))}\right|~=~\left|{w^{(n)}(t,x)-w^{(n)}(t,0+)+\varphi^{(n)}(t,x)\over a^{(n)}_1(t)-a^{(n)}_2(t)}\right|\\
~\leq~{1\over 2(\delta_1+\delta_2)}\cdot\left(\big|w^{(n)}(t,x)-w^{(n)}(t,0+)\big|+\big|\varphi^{(n)}(t,x)\big|\right)\\
~\leq~{1\over 2(\delta_1+\delta_2)}\cdot\left(M_0x^{1/2}+2\big|x\ln x\big|\right).
\end{multline}
Similarly, for every $(t,x)\in [\tau_0,0[\times ]-1/2,t[$, 
\begin{multline}\label{a2}
\left|a_n(t,x)-1-{\sigma^{(n)}_1(t)\over 2\left(a^{(n)}_1(t)-a^{(n)}_2(t)\right)}\right|~=~\left|{w^{(n)}(t,x)-w^{(n)}(t,t-)+\varphi^{(n)}(t,x)\over a^{(n)}_1(t)-a^{(n)}_2(t)}\right|\\
~\leq~{1\over 2(\delta_1+\delta_2)}\cdot\left(M_0\cdot |x-t|^{1/2}+2\cdot\big||x-t| \ln |x-t|\big|\right),
\end{multline}
and for any $(t,x)\in [\tau_0,0[\times ]t,0[$,
\begin{multline}\label{a3}
\left|a_n(t,x)-1+{\sigma_1^{(n)}(t)\over 2\left( a^{(n)}_1(t)-a^{(n)}_2(t)\right)}\right|~=~\left|{w^{(n)}(t,x)-w^{(n)}(t,t+)+\varphi^{(n)}(t,x)\over a^{(n)}_1(t)-a^{(n)}_2(t)}\right|\\
~\leq~{1\over 2(\delta_1+\delta_2)}\cdot\left(M_0\cdot |x-t|^{1/2}+2\cdot\left(\big||x-t| \ln |x-t|\big|+\big|x\ln|x|\big|\right)\right).
\end{multline}

Since 
\[
\begin{cases}
& \ds{\sigma^{(n)}_1(t)\over \sigma^{(n)}_1(t)+10\delta_2+2M_0\sqrt{|t|}}~\leq~{\sigma^{(n)}_1(t)\over 2( a^{(n)}_1(t)-a^{(n)}_2(t))}~\leq~{\sigma_1^{(n)}(t)\over \sigma_1^{(n)}(t)+6\delta_2-2M_0\sqrt{|t|}}\,, \\[2mm]
&\ds{\sigma_2^{(n)}(t)\over \sigma^{(n)}_2(t)+10\delta_1+2M_0\sqrt{|t|}}~\leq~{\sigma^{(n)}_2(t)\over 2( a^{(n)}_1(t)-a^{(n)}_2(t))}~\leq~{\sigma^{(n)}_2(t)\over \sigma^{(n)}_2(t)+6\delta_1-2M_0\sqrt{|t|}}\,,
\end{cases}
\]
by \eqref{sigma12},
we conclude, using \eqref{sigma12} and \eqref{Tau-as1} one more, that
\[
\begin{cases}
&\ds{\delta_1\over \delta_1+2\delta_2}~\leq~{\sigma^{(n)}_1(t)\over 2\big( a^{(n)}_1(t)-a^{(n)}_2(t)\big)}~\leq~{2\delta_1\over 2\delta_1+\delta_2 }\,, \\[4mm]
&\ds{\delta_2\over 2\delta_1+\delta_2}~\leq~{\sigma^{(n)}_2(t)\over 2\big( a^{(n)}_1(t)-a^{(n)}_2(t)\big)}~\leq~{2\delta_2\over \delta_1+2\delta_2}\,,
\end{cases}
\qquad\forall t\in [\tau_0,0].
\]
Therefore,  by (\ref{a1})-(\ref{a3}) there exists $\bar{\delta}>0$ such that  
\bel{an-rb}
\begin{cases}
\ds-{5\delta_2\over 2\delta_1+4\delta_2}~\leq~a_n(t,x)~\leq~ -{\delta_2\over 4\delta_1+2\delta_2}&\forall (t,x)\in [\tau_0,0[\times [0,\bar{\delta}],\\[4mm]
\ds 1-\ds{4\delta_1+\delta_0\over 4\delta_1+2\delta_2}~\leq~a_n(t,x)~\leq~ 1-{2\delta_1-\delta_0\over 2\delta_1+4\delta_2 }&\forall t\in [\tau_0,0], x\in [t,0[\cap ([t,t+\bar \delta[\cup]-\bar\delta, 0]), \\[4mm]
\ds 1+\ds{\delta_1\over 2\delta_1+4\delta_2}~\leq~a_n(t,x)~\leq~ 1+{5\delta_1\over 4\delta_1+2\delta_2}&\forall t\in [\tau_0,0], x\in [t-\bar{\delta},t].\end{cases}
\eeq

The next lemma provides the Lipschitz continuous dependence of the characteristic curves (\ref{char-s2}).

\begin{lemma}\label{Lip-char} Let $w^{(n)}$ and $\vp^{(n)}$ be as in 
(\ref{wn-c-1}) and (\ref{a-n-2}). Given $\tau\in [\tau_0,0[$, let $x_1,x_2\in \R\backslash\{\tau,0\}$ with $\bar{x}_2<\bar{x}_1$ such that both $\bar{x}_1$ and $\bar{x}_2$ belong to $]-{1\over 2e},\tau[$, $]\tau,0[$ or $]0,{1\over 2e}[$. Then 
\bel{x1-x2}
|x(t;\tau,\bar{x}_1)-x(t;\tau,\bar{x}_2)|~\leq~K_1\cdot |\bar{x}_2-\bar{x}_1|\qquad\forall t\in [\tau_0,\tau[
\eeq
for some $K_1>0$ depending only on $M_0$, $\delta_1$, $\delta_2$, $K$, and $b$.
\end{lemma}
{\bf Proof.} We shall  prove \eqref{x1-x2} for $\bar x_1$, $\bar x_2\in ]\tau,0[$. The other cases follow the same lines as the proof of \eqref{Lipx}. 

For any $t<x_2<x_1<0$, \eqref{a2} and \eqref{a1-a2} imply 
\begin{multline*}
|a_n(t,x_2)-a_n(t,x_1)|~\leq~ {\left|w^{(n)}(t,x_2)-w^{(n)}(t,x_1)\right|+\left|\varphi^{(n)}(t,x_2)-\varphi^{(n)}(t,x_1)\right|\over \big|a^{(n)}_1(t)-a^{(n)}_2(t)\big|}\\
~\leq~{1\over 2(\delta_1+\delta_2)}\cdot \left(Kb+\big|\ln|x_2-t|\big|+\big|\ln|x_1|\big|\right)\cdot (x_1-x_2).
\end{multline*}
Setting $0\leq z(t)\doteq x(t;\tau,\bar{x}_1)-x(t;\tau,\bar{x}_2)$, we obtain 
\begin{eqnarray*}
\dot{z}(t)&\geq &-{1\over 2(\delta_1+\delta_2)}\cdot \left(Kb+\big|\ln|x(t;\tau,\bar x_2)-t|\big|+\big|\ln|x(t;\tau,\bar x_1)|\big|\right) \cdot z(t).
\end{eqnarray*}
Since \eqref{an-rb} implies for any $\bar x\in (0,t)$ that
\[
\left(1-{4\delta_1+\delta_0\over 2(2\delta_1+\delta_2)}\right)(\tau-t)\leq x(\tau; \tau, \bar x)-x(t;\tau,\bar x)\leq \left( 1-{2\delta_1-\delta_0\over 2(\delta_1+2\delta_2)}\right) (\tau-t), 
\]
one ends ups with 
\[
{\dot{z}(t)\over z(t)}\geq-{1\over 2(\delta_1+\delta_2)}\cdot \left(Kb+2\cdot\left|\ln\left(\min\left({\delta_0\over 2(2\delta_1+\delta_2)},{\delta_0\over 2(\delta_1+2\delta_2)}\right) (\tau-t)\right)\right|\right),
\]
which yields (\ref{x1-x2}).
\endproof
\medskip

Next, consider the constants 
\bel{gd}
\gamma_0~\doteq~\min\left\{\ds{\delta_0\over 2(\delta_1+2\delta_2)},\ds{\delta_0\over 2(2\delta_1+\delta_2)}\right\},\quad \gamma_1~\doteq~\max\left\{{5\max(\delta_1,\delta_2)\over2( \delta_1+2\delta_2)},{5\max(\delta_1, \delta_2)\over 2(2\delta_1+\delta_2)}\right\},
\eeq
and define
\bel{I} 
I^{\tau}_{t}~\doteq~ [t-\gamma_0 (\tau-t),t+\gamma_0 (\tau-t)]\cup [-\gamma_0 (\tau-t),\gamma_0 (\tau-t)].
\eeq
From (\ref{an-rb}), one has
\bel{char2-s2}
x(t;\tau,\bar x)~\notin~I_t^{\tau}\qquad\forall \tau_0\leq t<\tau\leq 0, \bar x\in [-1/2,1/2]\backslash\{0,\tau\}.
\eeq
Furthermore, 
for all $\tau_0\leq t<\tau\leq 0$, one has 
\bel{char2-s22} \begin{cases}
\vert (x(t; \tau,\bar x)-t)-(\bar x -\tau)\vert \leq \gamma_1 (\tau-t), & \bar x\in [-\frac12, 0[\backslash \{\tau\},\\[3mm]
\vert x(t; \tau,\bar x)-\bar x\vert \leq \gamma_1(\tau-t), & \bar x\in ]-\tau, \frac12]\backslash \{0\}.
\end{cases}
\eeq

By the same arguments used  in \cite[Lemma 4.1]{BZ}, we  now obtain
\begin{lemma}\label{ln-b2-2s} Let $w^{(n)}$ and $\vp^{(n)}$ be as in 
(\ref{wn-c-1}) and (\ref{a-n-2}). There exists $\ve_0>0$ small enough, so that for any 
$-\ve_0\leq t<\tau<0$
 and any solution $v$ of the linear equation 
\[
v_t+a_n(t,x)\cdot v_x~=~0,\qquad v(\tau_0,\cdot)~=~\bar{v}~\in~ H^2(\R\backslash I^{\tau}_{t}),
\]
one has
\[
\|v(t,\cdot)\|_{H^2(\R\backslash\{t,0\})}~\leq~{3\over 2}\cdot \|\bar{v}\|_{H^2(\R\backslash I^{\tau}_{t})}.
\]
\end{lemma}
\medskip

{\bf Step 2.} Let us now consider a sequence of approximate solutions $w^{(k)}$ to (\ref{BH4}) inductively defined as follows. 
\begin{itemize}
\item  $w^{(1)}:[\tau_0,0[\times\R\to\R$ such that for all $t\in [\tau_0,0[$,
\[
w^{(1)}(t,x)~=~\begin{cases}
\overline{w}(x)&\qquad\mathrm{if}~~x\in (t,0)\cup (0,\infty),\\[2mm]
\overline{w}(x+\tau_0-t)&\qquad\mathrm{if}~~x\in (-\infty,t),
\end{cases}
\]
where $\ov w$ satisfies \eqref{ic-1-2s}.
\item For every $k\geq 1$, $w^{(k+1)}(t,\cdot)$ solves the linear equation 
\[
w_t+a_n(t,x)\cdot w_x~=~F^{(k)}(t,x),\qquad w(\tau_0,\cdot)~=~\overline{w}(\cdot)
\]
with  $F^{(k)}(t,x)\doteq F\left(t,x,w^{(k)}\right)$ and $\bar w$ as in \eqref{ic-1-2s}. This can be rephrased as  
\bel{w-k+1-2s}
w^{(k+1)}(t_0,x_0)~=~\overline{w}(x(\tau_0;t_0,x_0))+\int_{\tau_0}^{t_0}F^{(k)}(t,x(t;t_0,x_0))dt.
\eeq
\end{itemize}
The following lemma provides a priori estimates on  $w^{(k)}$, uniformly valid  for all $k\geq 1$. 

\begin{lemma}\label{wk-p-2s} Let $w^{(n)}$ and $\vp^{(n)}$ be as in 
(\ref{wn-c-1}) and (\ref{a-n-2}). Then there exists $\ve_0>0$ sufficiently small 
so that the following holds. If $\tau_0\in [-\ve_0,0[\,$, then   for every  $k\geq 0$ and  a.e.~$\tau\in [\tau_0,0[\,$, one has
\bel{wk-b2-2s}
\left|w^{(k)}(\tau,\tau\pm)-\overline{w}(\tau_0\pm)\right|~\leq~\delta_1,\qquad \left|w^{(k)}(\tau,0\pm)-\overline{w}(0\pm)\right|~\leq~\delta_2\,,
\eeq
\bel{dotsi-2s} 
\max\left\{\left|\dot{\sigma}^{(k)}_1(\tau)\right|,~\left|\dot{\sigma}_2^{(k)}(\tau)\right|\right\}~\leq~4C_1\big|\ln |\tau|\big|,
\eeq
\bel{wk-H2-2s}
\left\|w_x^{(k)}(\tau,x)\right\|_{{\bf L}^{\infty}(]\tau,0[)}~\leq~2K_2b,\qquad \left\|w^{(k)}(\tau,\cdot)\right\|_{H^2(\R\backslash \{\tau,0\})}~\leq~M_0,
\eeq
for some positive constants $C_1$ and $K_2$.
\end{lemma}
{\bf Proof.} It is clear that  (\ref{wk-b2-2s})-(\ref{wk-H2-2s}) hold for $k=1$. By induction, assume that  (\ref{wk-b2-2s})-(\ref{wk-H2-2s}) hold for a given $k\geq 1$. 

{\bf 1.} We shall establish the first inequality in \eqref{wk-H2-2s}. 
Given $\tau\in [\tau_0,0[$ and $\tau<\bar{x}_2<\bar{x}_1<0$,  consider the characteristics $t\,\mapsto\, x_i(t)\,=\,x( t; \tau, \bar x_i)$ for $i\in\{1,2\}$,
which satisfy, cf. (\ref{char2-s2}),
\bel{x2-x1-2s}
\min\left\{|x_i(t)|,|x_i(t)-t|\right\}~\geq~ \gamma_0\cdot (\tau-t)\qquad\forall t\in [\tau_0,\tau], i\in\{1,2\}.
\eeq
Recalling (\ref{w-k+1-2s}), (\ref{x1-x2}), and (\ref{F1-2s}), we estimate
\begin{multline*}
\left|w^{(k+1)}(\tau,\bar{x}_2)-w^{(k+1)}(\tau,\bar{x}_1)\right|\\~\leq~\left|\overline{w}(x_2(\tau_0))-\overline{w}(x_1(\tau_0))\right|+\int_{\tau_0}^{\tau}\left|F^{(k)}(t,x_2(t))-F^{(k)}(t,x_1(t))\right|dt\\
~\leq~K_1\cdot\left(b+{C_1(1+2M_0+2K_2b)\over \gamma_0^{1/2}\delta_0}\cdot\int_{\tau_0}^{\tau}{1\over (\tau-t)^{1/2}}dt\right)\cdot |\bar{x}_2-\bar{x}_1|\\
~\leq~K_1\cdot \left(b+{2C_1(1+2M_0+2K_2b)(\tau-\tau_0)^{1/2}\over \gamma_0^{1/2}\delta_0}\right)\cdot |\bar{x}_2-\bar{x}_1|.
\end{multline*}
Thus, if  $0\leq -\tau_0\leq \ds \left({b\gamma^{1/2}_0\delta_0\over 2C_1(1+2M_0+2K_2b)}\right)^2$, then 
\[
\left|w^{(k+1)}(\tau,\bar{x}_2)-w^{(k+1)}(\tau,\bar{x}_1)\right|~\leq~2K_2b\cdot |\bar{x}_2-\bar{x}_1|
\]
and (\ref{wk-H2-2s}) is satisfied by $w^{(k+1)}$.
\medskip

{\bf 2.} We shall establish (\ref{dotsi-2s}) for $i=2$ and the second inequality in (\ref{wk-b2-2s}).  The other ones are quite similar.  Given any $\tau_0\leq \tau_1<\tau_2\leq 0$, let $t\mapsto x_2^{\pm}(t)\doteq x(t;\tau_2,0\pm)$ be  the characteristics,which reach the origin at time $\tau_2$ from the positive and negative side, respectively. From  (\ref{char2-s22}), it follows that
\[
\left|x^{\pm}_{2}(t)\right|~\leq~\gamma_1\cdot |\tau_2-t|\qquad\forall t\in [\tau_1,\tau_2].
\]
Furthermore, recalling  \eqref{dotsi-2s}, (\ref{wk-H2-2s}), and (\ref{F1-2s}), we have 
\begin{multline*}
\left|w^{(k+1)}(\tau_2,0\pm)-w^{(k+1)}(\tau_1,0\pm)\right|~\leq~\left|w^{(k+1)}(\tau_1, x^{\pm}_{2}(\tau_1))-w^{(k+1)}(\tau_1,0\pm)\right|\\+\int_{\tau_1}^{\tau_2} \left| F^{(k)}(t, x^{\pm}_{2}(t))\right| dt~\leq~2(K_2b+M_0)\gamma_1\cdot (\tau_2-\tau_1)\\+C_1\cdot\int_{\tau_1}^{\tau_2}{1+2M_0+2K_2b\over \delta_0}+ \left(1+{8C_1\over \delta_0}(\vert t\vert^{1/2}+ 2\gamma_1^{1/2}(\tau_2-t)^{1/2})\right)\vert \ln|t|\vert~dt\\
~\leq~\left[2(K_2b+M_0)\gamma_1+{C_1(1+2M_0+2K_2b)\over \delta_0}+C_1\left(1+{8C_1(1+2\gamma_1^{1/2})|\tau_0|^{1/2}\over \delta_0}\right) \big|\ln |\tau_2|\big|\right]\cdot (\tau_2-\tau_1).
\end{multline*}
Thus, for $|\tau_0|$ is sufficiently small, we then obtain (\ref{dotsi-2s}) for $(k+1)$ and $i=2$ by
\[
\left|\sigma_2^{(k+1)}(\tau_2)-\sigma_2^{(k+1)}(\tau_1)\right|~\leq~4C_1\big|\ln|\tau_2|\big|\cdot (\tau_2-\tau_1).
\]
Moreover, for the second inequality in \eqref{wk-b2-2s}, choose $\tau_1=\tau_0$ in the above estimate, i.e., 
\begin{multline*}
\left|w^{(k+1)}(\tau_2,0\pm)-\overline{w}(0\pm)\right|\\~\leq~\left[2(K_2b+M_0)\gamma_1+{C_1(1+2M_0+2K_2b)\over \delta_0}+C_1\left(1+{8C_1\over \delta_0}\right)\cdot \big|\ln |\tau_0|\big|\right]\cdot |\tau_0|.
\end{multline*}
and this yields \eqref{wk-b2-2s}.

\medskip

{\bf 3.} Finally, from Duhamel' formula, Lemma \ref{ln-b2-2s}, \eqref{ic-1-2s}, Lemma \ref{F1-g-2s}, (\ref{dotsi-2s}), and (\ref{I}), we obtain, for all $\tau\in [\tau_0,0]$,
\begin{multline*}
\left\|w^{(k+1)}(\tau,\cdot)\right\|_{H^2(\R\backslash \{\tau,0\})}~\leq~{3\over 2} \|\overline{w}\|_{H^2\big(\R\backslash I^{\tau}_{\tau_0}\big)}+{3\over 2}\cdot\int_{\tau_0}^{\tau}\left\|F^{(k)}(t,\cdot)\right\|_{H^2(\R\backslash I^{\tau}_{t})}dt\\
~\leq~{3M_0\over 8}+{3\over 2}C_1\cdot\int_{\tau_0}^{\tau}{1+2M_0+2K_2b+8C_1\big|\ln|t|\big|\over \delta_{0}\gamma_0^{2/3}}(\tau-t)^{-2/3}+\big|\ln|t|\big|~dt\\
~\leq~{3M_0\over 8}+\left({9C_1(25+2M_0+2K_2b+8C_1\cdot\big|\ln|\tau_0|\big|)\over  2\delta_{0}\gamma_0^{2/3}}\right)\cdot|\tau_0|^{1/3}+3C_1\cdot\big|\tau_0\ln|\tau_0|\big|.
\end{multline*}
Identifying an upper bound on $\tau_0$ such that the right hand side is less or equal than $M_0$, shows that the second bound in (\ref{wk-H2-2s}) is satisfied by $w^{(k+1)}$ as well.
\endproof
\medskip

Thanks to the above estimates, we can now prove that the sequence of approximations $w^{(k)}$
is Cauchy, and converges to a solution $w$ of the linear problem (\ref{BH4}).  This will accomplish
the inductive step, toward the proof of Theorem~\ref{t:61}.

\begin{lemma}\label{l:44-2s}  There exists $\ve_0>0$ sufficiently small so that, for all $\tau_0\in [-\ve_0, 0[\,$
 the following holds:
Let $w^{(n)}, \vp^{(n)}$ as in 
(\ref{wn-c-1}) and (\ref{a-n-2}).
 Then the sequence of approximations $\left(w^{(k)}(t,\cdot)\right)_{k\geq 1}$ converges uniformly for all $t\in [\tau_0, 0[$ to a limit function $w(t,\cdot)$ in $H^2(\R\backslash\{t,0\})$. Namely,
\[
\lim_{k\to\infty}\sup_{t\in [\tau_0,0]}\left\|w^{(k)}(t,\cdot)-w(t,\cdot)\right\|_{H^2(\R\backslash \{t,0\})}~=~0.
\]
The function $w$ provides a solution to  the Cauchy  problem (\ref{BH4}) and satisfies for all $\tau\in [\tau_0,0]$
\bel{w-b2-2s}
|w(\tau,\tau\pm)-\overline{w}(\tau_0\pm)|~\leq~\delta_1,\qquad |w(\tau,0\pm)-\overline{w}(0\pm)|~\leq~\delta_2\,,
\eeq
\bel{w-H2-2s}
\left\|w_x(\tau,x)\right\|_{{\bf L}^{\infty}(]\tau,0[)}~\leq~2K_2b,\qquad \left\|w(\tau,\cdot)\right\|_{H^2(\R\backslash \{\tau,0\})}~\leq~M_0.
\eeq
 Moreover, $\sigma_1(t)\doteq  w(t,t-)-w(t,t+)$ and $\sigma_2(t)\doteq  w(t,0-)-w(t,0+)$ are locally Lipscthitz in $(\tau_0,0)$  and 
\bel{dotsi-2s-w} 
\max\left\{\left|\dot{\sigma}_1(\tau)\right|,~\left|\dot{\sigma}_2(\tau)\right|\right\}~\leq~4C_1\big|\ln \tau|\big|\quad a.e.~\tau\in [\tau_0,0].
\eeq
\end{lemma}

{\bf Proof.}  {\bf 1.} For any $k\geq 1$, we set 
\bel{set11}
\begin{cases}
\quad z^{(k)}&\doteq~w^{(k+1)}-w^{(k)},\qquad  M^{(k)}_z(\tau)~\doteq~\left\|z^{(k)}(\tau,\cdot)\right\|_{H^2(\R\backslash\{\tau,0\})},\\[3mm]
 \sigma_{1}^{(k,z)}(\tau)&\doteq~z^{(k)}(\tau,\tau-)-z^{(k)}(\tau,\tau+),\qquad \sigma_{2}^{(k,z)}(\tau)~\doteq~z^{(k)}(\tau,0-)-z^{(k)}(\tau,0+),\\[3mm]
  \alpha_k(\tau)&\doteq~\ds\sup_{t\in [\tau_0,\tau]}\max_{i\in\{1,2\}}\left|\dot{\sigma}_i^{(k,z)}(t)\right|\,,\\[3mm]
\beta_k(\tau)&\doteq~\ds\sup_{t\in [\tau_0,\tau]} \left(M^{(k)}_{z}(t)+\left|z^{(k)}_{x}(t,0-)\right|+\left|z_x^{(k)}(t,0+)\right|+\left|z^{(k)}(t,0-)\right|+\left|z^{(k)}(t,0+)\right|\right).
\end{cases}
\eeq

Recalling Duhamel's formula, Lemma \ref{ln-b2-2s}, Lemma \ref{FF-2s}, (\ref{I}), and Lemma \ref{wk-p-2s}, for all $\tau\in [\tau_0,0]$ we  estimate 
\bel{Mk-2s}
\bega{l}
\ds M^{(k+1)}_{z}(\tau)~\leq~{3\over 2}\cdot \int_{\tau_0}^{\tau}\left\|F^{(k+1)}(t,\cdot)-F^{(k)}(t,\cdot)\right\|_{H^2(\R\backslash I^{\tau}_{t})}dt\\[4mm]
\quad\leq~\ds{3C_2\over 2\delta_0^2}\cdot\int_{\tau_0}^{\tau}\beta_{k}(t)\cdot\left[\left({4C_1\big|\ln |t|\big|\over \delta_0 \gamma_0^{2/3}}+{(2M_0+2K_2b+1)(1+\delta_0)\over \delta_0\gamma_0^{2/3}}\right)\cdot (\tau-t)^{-2/3}+\big|\ln |t|\big|\right]dt\\[5mm]
\qquad\qquad\qquad\qquad\qquad \ds +{3C_2\over 2\delta_0^2}\cdot\int_{\tau_0}^{\tau}{2M_0+2K_2b\over \gamma_0^{2/3}}\cdot \alpha_{k}(t)\cdot (\tau-t)^{-2/3}~dt\\[4mm]
\qquad\qquad~~\leq~\ds C_3\cdot \left(|\tau_0-\tau|^{1/3}\cdot \big|\ln |\tau_0-\tau|\big|\cdot \beta_{k}(\tau)+ \int_{\tau_0}^\tau \alpha_k(t) (\tau-t)^{-2/3}dt\right)
 \enda 
 \eeq
 for some constant $C_3$ depending only on $M_0$, $b$, $K_2$,  $\delta_1$, and $\delta_2$.
 \medskip
 
 {\bf 2.}   We now establish a bound on $z^{(k+1)}_x(\tau,0\pm)$. Since 
 \bel{Mk-2s2}
\left \vert z_x^{(k+1)}(\tau,0+)\right\vert ~\leq~\left\| z^{(k+1)}(\tau,\cdot)\right\|_{H^2(\R\backslash\{t,0\})}\leq M^{(k+1)}_z(\tau),
 \eeq
 it suffices to have a closer look at $z^{(k+1)}_x(\tau,0-)$. Given $\tau\in [\tau_0,0[$ and $\tau<\bar{x}_2<\bar{x}_1<0$, consider the characteristics $t\,\mapsto\, x_i(t)\,=\,x( t; \tau, \bar x_i)$ for $i\in\{1,2\}$.
Recalling (\ref{w-k+1-2s}), (\ref{F1-F2-0-2ss}), (\ref{x2-x1-2s}), and (\ref{x1-x2}), we estimate 
\begin{multline*}
\left|z^{(k+1)}(\tau,\bar{x}_2)-z^{(k+1)}(\tau,\bar{x}_1)\right|\\~=~\left|\int_{\tau_0}^{\tau}\left[F^{(k+1)}(t,x_2(t))-F^{(k+1)}(t,x_1(t))\right]-\left[F^{(k)}(t,x_2(t))-F^{(k)}(t,x_1(t))\right]dt\right|\\
~\leq~\int_{\tau_0}^{\tau}\int_0^1\left|\left(F^{(k+1)}_x-F^{(k)}_x\right)(t,(1-s)x_1(t)+sx_2(t))\right|\cdot |x_1(t)-x_2(t)|~ds~dt\\\
~\leq~{2C_2K_1\over \delta^2_0\gamma_0^{1/2}}\cdot \int_{\tau_0}^{\tau}\beta_k(t)\cdot (1+\delta_0+2M_0+2K_2b)\cdot (\tau-t)^{-1/2}~dt\cdot |\bar{x}_2-\bar{x}_1| \\
~\leq~C_4\cdot \beta_k(\tau)\cdot |\tau-\tau_0|^{1/2}\cdot |\bar{x}_2-\bar{x}_1|\
\end{multline*}
for some constant $C_4$ depending only on $M_0$, $b$, $K_1$, $\delta_1$, and $\delta_2$. This implies 
\bel{zx-2s-0}
\left|z^{(k+1)}_x(\tau,0-)\right|~\leq~C_4\cdot \beta_k(\tau)\cdot |\tau-\tau_0|^{1/2}.
\eeq

 {\bf 3.} Finally, we  establish a bound on $\alpha_{k+1}(\tau)$ for $\tau\in [\tau_0,0]$. We only present here the details for $\dot\sigma_2^{(k+1,z)}(t)$, since $\dot\sigma_1^{(k+1,z)}(t)$ can estimated in the same way. Given any $\tau_0\leq\tau_1<\tau_2< 0$, denote by  $t\mapsto x_2^{\pm}(t)\doteq x(t;\tau_2,0\pm)$  the characteristics which reach the origin  at time $\tau_2$ from the positive and negative side, respectively. Using (\ref{w-k+1-2s}), (\ref{F1-F2-0-2s}), (\ref{char2-s22}), and (\ref{x2-x1-2s}), we estimate 
 $$\bega{l}\ds
\left|z^{(k+1)}(\tau_2,0\pm)-z^{(k+1)}(\tau_1,0\pm)\right|\\[3mm]
\quad \ds \leq~\left|z^{(k+1)}\big(\tau_1,x^{\pm}_2(\tau_1)\big)-z^{(k+1)}(\tau_1,0\pm)\right|+\int_{\tau_1}^{\tau_2}\left|F^{(k+1)}\left(t,x_2^{\pm}(t)\right)-F^{(k)}\left(t,x_2^{\pm}(t)\right)\right|dt\\[4mm]
\quad\ds\leq~2\beta_{k+1}(\tau_1)\cdot\left|x^{\pm}_2(\tau_1)\right|\ds+{C_2\over \delta^2_0}\cdot \int_{\tau_1}^{\tau_2}(2M_0+2K_1b)\cdot \alpha_k(t)\cdot \gamma_1^{1/2} |\tau_2-t|^{1/2}dt\\[4mm]
\quad \ds+{C_2\over \delta^2_0}\cdot \int_{\tau_1}^{\tau_2}\beta_{k}(t)\cdot \left({4C_1\sqrt{\vert t\vert}+\delta_0\over \delta_0}\big|\ln|t|\big|+{2M_0+2K_1b\over\delta_0}+{8C_1\gamma_1^{1/2}\over \delta_0}\cdot \big|\ln |t|\big| |\tau_2-t|^{1/2}\right)dt\\[4mm]
\qquad\leq~\ds\left(2\beta_{k+1}(\tau_2)\gamma_1+C_5\cdot\left[\beta_k(\tau_2)\cdot\big|\ln|\tau_2|\big|+\alpha_k(\tau_2)\cdot  |\tau_1-\tau_2|^{1/2}\right]\right)\cdot (\tau_2-\tau_1)
\enda $$
for some constant $C_5$ depending only on $M_0$, $b$, $K_1$, $\delta_1$, and $\delta_2$. Thus, for $\tau\in [\tau_0,0]$,
\bel{sgm1-2s}
\alpha_{k+1}(\tau)~\leq~\ds 2\beta_{k+1}(\tau)\gamma_1+C_5\cdot\beta_k(\tau)\cdot\big|\ln|\tau|\big|.
\eeq
Moreover, by choosing $\tau_1=\tau_0$ and $\tau_2=\tau\in [\tau_0,0]$, we also get
\[
\left|z^{(k+1)}(\tau,0\pm)\right|~\leq~C_6\cdot\left(\beta_k(\tau)\cdot\big||\tau-\tau_0|\cdot \ln|\tau-\tau_0|\big| +\int_{\tau_0}^\tau \alpha_k(t)\cdot |t-\tau|^{1/2}dt\right),
\]
and (\ref{Mk-2s})-(\ref{sgm1-2s}) imply that
\begin{multline*}\label{beta-2s-e}
\beta_{k+1}(\tau)~\leq~C_7\cdot \left(|\tau-\tau_0|^{1/3}\cdot \big|\ln |\tau-\tau_0|\big|\cdot \beta_{k}(\tau)+\int_{\tau_0}^{\tau} \alpha_k(t)\cdot |\tau-t|^{-2/3}dt\right)\\
~\leq~C_8\cdot |\tau-\tau_0|^{1/3}\cdot \big|\ln |\tau-\tau_0|\big|\cdot \left(\beta_{k}(\tau)+\beta_{k-1}(\tau)\right).
\end{multline*}
In particular, for $\tau_0<0$ sufficiently close to $0$, we get 
\[
\beta_{k+1}(\tau)+{1\over 2}\cdot \beta_k~\leq~{3\over 4}\cdot \left(\beta_{k}(\tau)+{1\over 2}\cdot \beta_{k-1}\right),
\]
which implies
\[
\sum_{k=1}^{\infty}\sup_{\tau\in [\tau_0,0]}\left\|z^{(k)}(\tau,\cdot)\right\|_{H^2(\R\backslash\{t,0\})}~<~\sum^{\infty}_{k= 1}\beta_k(\tau)~<~\infty.
\]
We thus conclude that  $\left(w^{(k)}(\tau,\cdot)\right)_{k\geq 1}$ converges uniformly for all $\tau \in [\tau_0,0[$ to a limit function $w(\tau,\cdot)$ in $H^2(\R\backslash\{t,0\})$, which provides the solution to the linear problem (\ref{BH4}). Moreover, since $\ds\lim_{k\to\infty}w^{(k)}(\tau,0\pm)=w(\tau,0\pm)$ and $\ds\lim_{k\to\infty}w^{(k)}(\tau,\tau\pm)=w(\tau,\tau\pm)$, one has that $\ds\lim_{k\to\infty}\sigma_i^{(k)}(\tau)=\sigma_i(\tau)$ for all $\tau\in [0,\tau]$. Furthermore, $\ds\lim_{k\to\infty}w_x^{(k)}(\tau,0\pm)=w_x(\tau,0\pm)$ and $\ds\lim_{k\to\infty}w_x^{(k)}(\tau,\tau\pm)=w_x(\tau,\tau\pm)$ and hence Lemma \ref{wk-p-2s} implies that $w$ satisfies 
(\ref{w-b2-2s})-(\ref{dotsi-2s-w}).
\endproof
\medskip

We are now ready to complete the proof of  our second  main theorem, describing the asymptotic behavior of solutions up to the time when two shocks interact.
\medskip

{\bf Proof of Theorem \ref{t:61}.}   {\bf 1.} By induction, we construct a sequence of approximate solutions $\big(w^{(n)}\big)_{n\geq 1}$ where each $w^{(n+1)}$ is the solution to  the linear  problem (\ref{BH4}).  Assuming
that $\tau_0\in [-\ve_0, 0[\,$ is sufficiently close to 0, we claim that 
\bel{sum1-2s}\sum_{n\geq 1}\,\left\|w^{(n+1)}(t,\cdot)-w^{(n)}(t,\cdot)\right\|_{H^1(\R\setminus\{t,0\})}~<~\infty\qquad\forall t\in [\tau_0,0].
\eeq
For a fixed $n\geq 2$, recalling that $a_n(t,x)=\ds {w^{(n)}(t,x)+\varphi^{(n)}(t,x)-a^{(n)}_2(t)\over a^{(n)}_1(t)-a^{(n)}_2(t)}$, we define
\[
\begin{cases}
W^{(n)}\,\doteq\,w^{(n)}-w^{(n-1)},\quad  A^{(n)}(\tau,x)\,\doteq\,a_{n}(\tau,x)-a_{n-1}(\tau,x),\\[3mm]
\sigma^{(n)}_1(\tau)\,=\,W^{(n)}(\tau,\tau-)-W^{(n)}(\tau,\tau+),\quad \sigma^{(n)}_2(\tau)\,=\,W^{(n)}(\tau,0-)-W^{(n)}(\tau,0+),\\[3mm]
\ds \beta^{(n)}(\tau)\,\doteq\,\sup_{t\in [\tau_0,\tau]}\Big[\bigl\|W^{(n)}(\tau,\cdot)\bigr\|_{H^1(\R\backslash\{t,0\})}+\left|W^{(n)}(t,0-)\right|+\left|W^{(n)}(t,0+)\right|\Big].
\end{cases}
\]
Set $Z^{(n)}=W^{(n)}+V^{(n)}$ with $V^{(n)}=v^{(n)}-v^{(n-1)}$ and $v^{(n)}=\vp^{(n)}-\phi_0(x-t)-\phi_0(x)$ . From the above definitions, by (\ref{BH4}), we deduce
\bel{W11-2s}
Z^{(n+1)}_t+a_n\cdot Z^{(n+1)}_x~=~-\left(A^{(n)}w^{(n)}_x+A^{(n+1)}v^{(n+1)}_x\right) +G^{(n+1)}-G^{(n)}
\eeq
with 
\[
G^{(n)}(t,x)~=~{{\bf H}\left[w^{(n)}(t,\cdot)+\vp^{(n)}(t,\cdot)\right](x)\over a_1^{(n)}(t)-a_2^{(n)}(t) }-a_n(t,x)\cdot \big[\phi_0'(x-t)+\phi_0'(x)\big].
\]
We split  
$$w^{(n)}~=~v_{1,n}+v_{2,n},\qquad v_{2,n}(t,x)=\begin{cases} w(t,0-) \cdot\eta(x), & x<0,\\[3mm]
w(t,0+) \cdot\eta(x), & 0<x.
\end{cases}$$
Recalling the definition of $\Tilde{B}$ and $\tilde{C}$ in (\ref{AB-2s})-(\ref{C-2s}), we write 
\begin{multline*}
G^{(n)}(t,x)+\phi'_0(x-t)~=~- {\left[v_{2,n}(t,x)-v_{2,n}(t,t)\right]\cdot \phi'_0(x-t)+\left[v_{1,n}(t,x)-v_{1,n}(t,0)\right]\cdot\phi'_0(x)\over a_1^{(n)}(t)-a_2^{(n)}(t)}\\
+\frac{\Tilde{B}^{(v_{2,n})}(t,x)+\Tilde{C}^{(v_{1,n})}(t,x)+ {\bf H}[\vp^{(n)}(t,\cdot)](x)-\vp^{(n)}(t,x)\cdot [\phi_0(x-t)+\phi_0(x)]'}{a_1^{(n)}(t)-a_2^{(n)}(t)}.
\end{multline*}
%where 
%\[ v_2^{(n)}(t,x)=\begin{cases} w(t,0-) \cdot\eta(x), & x<0,\\
%w(t,0+) \cdot\eta(x), & 0<x,
%\end{cases}
%\]
%so that we can write 
%\begin{align*}
%G^{(n)}(t,x)+\phi'_0(x-t)= \frac{1}{a_1^{(n)}(t)-a_2^{(n)}(t)}& \Big( {\bf H}[v_1^{(n)}](t,x)-(v_1^{(n)}(t,x)-\frac {v_1^{(n)}(t,t-)+v_1^{(n)}(t,t+)}{2})\phi_0'(x-t) \\
%& + {\bf H}[v_2^{(n)}](t,x)-(v_2^{(n)}(t,x)-\frac{v_2^{(n)}(t,0-)+v_2^{(n)}(t,t+)}{2})\phi_0'(x)\\
%& + {\bf H}[\vp^{(n)}](t,x)- \vp^{(n)}(t,x)(\phi'_0(x-t)+\phi'_0(x))\\
%& -(v_2^{(n)}(t,x)-v_{2}^{(n)}(t,t))\phi'_0(x-t)\\
%& -(v_1^{(n)}(t,x)-v_1^{(n)}(t,0))\phi'_0(x)\Big).
%\end{align*}
Here it is important to note that $W_j^{(n)}= v_j^{(n)}-v_j^{(n-1)}$ satisfies 
\begin{equation*}
\left\| W_j^{(n)}(t, \cdot)\right\|_{H^1(\R\backslash \{t,0\})}\leq \O(1)\cdot\left(M^{(n)}(t)+\left| W^{(n)}(t,0-)\right|+\left|W^{(n)}(t,0+)\right|\right)\leq \O(1)\cdot\beta^{(n)}(t),
\end{equation*}
while, \eqref{a2-a1} implies,
\begin{equation*}
\|v_j^{(n)}(t,\cdot)\|_{H^1(\R\backslash\{t,0\})}\leq \O(1) \| w^{(n)}\|_{H^1(\R\backslash \{t,0\})}\leq \O(1)M_0.
\end{equation*}
Recalling \eqref{a2-a1}, \eqref{Bb1}, \eqref{B2-B1}, \eqref{vpd2}, \eqref{Hg-b1}, \eqref{e1-gb}, and \eqref{Hvp-s2}, we get  
\bel{G-2s}
\left\{ \bega{l}
\left|\left(G^{(n+1)}-G^{(n)}\right)(\tau,x)\right|\ds\leq~\Gamma_1\cdot \beta^{(n+1)}(\tau), \quad \tau<x<0,\\[2mm]
\ds \left\|(G^{(n+1)}-G^{(n)})(\tau,\cdot)\right\|_{H^1(\R\backslash [-\delta,\delta]\cup [\tau-\delta,\tau+\delta])}~\leq~\ds \Gamma_1\cdot\frac{\beta^{(n+1)}(\tau)}{\delta^{1/2}},
\enda\right.
\eeq
for some positive constant $\Gamma_1$. Furthermore, we have  for all $x[-\frac12 ,\frac12]\backslash \{\tau,0\}$ that 
\bel{A-2s}
\left\{ \bega{rl}\ds
\left| A^{(n)}w^{(n)}_{x}(\tau,x)\right|\ds \leq ~\Gamma_2  \beta^{(n)}(\tau),\quad \bigl| A^{(n)}v^{(n)}_{x}(\tau,x)\bigr|\ds \leq ~\Gamma_2  \beta^{(n)}(\tau)\cdot\left(|x|^{-1/2}+|x-\tau|^{-1/2}\right),\\[3mm]
\ds \bigl\|A^{(n)} w^{(n)}_{x}(\tau,\cdot)\bigr\|_{H^1(\R\backslash \{\tau,0\})}~\leq~\Gamma_2 \beta^{(n)}(\tau),~~ \bigl\|A^{(n)}v^{(n)}_x(\tau,\cdot)\bigr\|_{H^1(\R\backslash [-\delta,\delta]\cup [\tau-\delta,\tau+\delta])}~\leq~ \frac{\Gamma_2\beta^{(n)}(\tau)}{\delta^{1/2}},
\enda\right.
\eeq
for some constant $\Gamma_2>0$, dependent on $M_0$, $b$ ,$\delta_1$, and $\delta_2$. Hence, if $\tau_0<0$ is sufficiently close to $0$, we have, using Duhamel's formula and \eqref{gd}, for all $\tau\in [\tau_0,0[$ that 
\begin{multline*}
\left\|Z^{(n+1)}(\tau, \cdot)\right\|_{H^1(\R\backslash \{\tau,0\})}~\leq~\frac32 \int_{\tau_0}^{\tau}\left\|\left[G^{(n+1)}-G^{(n)}-A^{(n)} w^{(n)}_{x}-A^{(n+1)}v_x^{(n+1)}\right](t,\cdot)\right\|_{H^1(\R\backslash I^\tau_t )}dt\\
~\ds\leq~\frac32(\Gamma_1+\Gamma_2)\cdot\int_{\tau_0}^{\tau}\beta^{(n)}(t)+\beta^{(n+1)}(t)\cdot \gamma_0^{-1/2}(\tau-t)^{-1/2} dt.
\end{multline*}
Thus, there exists a constant $\Gamma_3>0$ dependent on $M_0$, $b$, $\delta_1$, and $\delta_2$ such that 
\bel{beta-n-2s}
\left\|Z^{(n+1)}(\tau, \cdot)\right\|_{H^1(\R\backslash \{\tau,0\})}~\leq~\Gamma_3\cdot \left(\beta^{(n)}(\tau)\cdot |\tau_0-\tau|+\beta^{(n+1)}(\tau)\cdot |\tau_0-\tau|^{1/2}\right).
\eeq

{\bf 2.}  We establish a bound on $\left| Z^{(n+1)}(\tau, 0\pm)\right|$. Given any $\tau_0\leq \tau\leq 0$, let $t\mapsto x_2(t)\doteq x(t;\tau,0-)$ be the characteristic, which reaches the origin at time $\tau$ from the negative side. Since 
\begin{equation}
Z^{(n+1)}(\tau_0,x)~\equiv~0,\qquad \left| Z^{(n+1)}(\tau,0+)\right|\leq \| Z^{(n+1)}(\tau, \cdot)\|_{H^1(\R\backslash \{\tau,0\})},
\end{equation}
 we have 
\begin{multline}\label{beta-n-2s2}
\left|Z^{(n+1)}(\tau,0-))\right|~\leq~\int_{\tau_0}^{\tau}\left|(-A^{(n)}w^{(n)}_{x}-A^{(n+1)}v_x^{(n+1)}+G^{(n+1)}-G^{(n)})(t,x_2(t))\right|dt\\
~\leq~(\Gamma_1+\Gamma_2)\cdot\int_{\tau_0}^{\tau}\beta^{(n)}(t)+ \beta^{(n+1)}(t)\cdot \left(1+2\gamma_0^{-1/2}(\tau-t)^{-1/2}\right)dt\\
\leq~ \Gamma_4\cdot \left(\beta^{(n)}(\tau)\cdot \vert \tau-\tau_0\vert +\beta^{(n+1)}(\tau)\cdot \vert \tau-\tau_0\vert^{1/2}\right),
\end{multline} 
where we used \eqref{char2-s22} and $\Gamma_4$ denotes a positive constant dependent on $M_0$, $b$, $\delta_1$, and $\delta_2$. 
Combining \eqref{beta-n-2s} -\eqref{beta-n-2s2}, we end up with
\begin{multline*}
\left\|Z^{(n+1)}(\tau, \cdot)\right\|_{H^1(\R\backslash\{\tau,0\})} + \left| Z^{(n+1)}(\tau,0-)\right|+ \left|Z^{(n+1)}(\tau,0+)\right|\\
~\leq~\Gamma_5\cdot\left( \beta^{(n)}(\tau)\cdot \vert \tau-\tau_0\vert + \beta^{(n+1)}(\tau)\cdot \vert \tau-\tau_0\vert^{1/2}\right),
\end{multline*}
where $\Gamma_5>0$ denotes a constant dependent on $M_0$, $b$, $\delta_1$, and $\delta_2$. 
\medskip

{\bf 3.} From \eqref{vpd2}, it holds that $W^{(n+1)}(\tau, \tau\pm)= Z^{(n+1)}(\tau,\tau\pm),  W^{(n+1)}(\tau,0\pm)= Z^{(n+1)}(\tau,0\pm)$, and 
\begin{multline*}
\left\|V^{(n+1)}(\tau,\cdot)\right\|_{H^1(\R\backslash\{\tau,0\})}\\~\leq~ \Gamma_6\cdot \left(\left| W^{(n+1)}(\tau,\tau-)\right|+\left| W^{(n+1)}(\tau,\tau+)\right|+ \left| W^{(n+1)}(\tau,0-)\right|+\left| W^{(n+1)}(\tau,0+)\right|\right)\\
~=~ \Gamma_6 \cdot \left(\left| Z^{(n+1)}(\tau,\tau-)\right|+\left| Z^{(n+1)}(\tau,\tau+)\right|+ \left| Z^{(n+1)}(\tau,0-)\right|+\left| Z^{(n+1)}(\tau,0+)\right|\right)\\
~\leq~ 3\Gamma_6\cdot \left(\|Z^{(n+1)}(\tau, \cdot)\|_{H^1(\R\backslash \{\tau,0\})}+ \left| Z^{(n+1)}(\tau,0-)\right|+\left| Z^{(n+1)}(\tau,0+)\right|\right)
\end{multline*}
for some positive constant $\Gamma_6$ on $M_0$, $b$, $\delta_1$, and $\delta_2$. Thus, we end up with 
\begin{align*}
\beta^{(n+1)}(\tau)& \leq (1+3\Gamma_6)\left(\|Z^{(n+1)}(\tau, \cdot)\|_{H^1(\R\backslash \{\tau,0\})}+ \left| Z^{(n+1)}(\tau,0-)\right|+\left| Z^{(n+1)}(\tau,0+)\right|\right)\\
& \leq \Gamma_7 \left(\beta^{(n)}(\tau)\cdot \vert \tau-\tau_0\vert + \beta^{(n+1)}(\tau) \cdot \vert \tau-\tau_0\vert^{1/2}\right).
\end{align*}
Provided that $\tau_0<0$ is sufficiently close to $0$, we obtain that 
\[
\beta^{(n+1)}(\tau)~\leq~\beta^{(n)}(\tau)/2\qquad\forall \tau\in [\tau_0,0[.
\]
 Thus, (\ref{sum1-2s}) holds for all  $\tau\in [\tau_0,0]$, and the sequence of approximations $w^{(n)}(\tau,\cdot)$ is Cauchy in the space $H^1\bigl(\R\backslash\{\tau,0\}\bigr)$, and hence it converges to a unique limit $w(\tau,\cdot)$. 
 
 It remains to check that this limit function $w$ is an entropic solution, i.e., it satisfies, cf. \eqref{BH11}, \eqref{form-u}, and \eqref{BH-3},
\[
\left(w+\vp^{(w)}\right)(t_0,x_0)=(\ov w+ \ov \vp)(x(\tau_0))+\int_{\tau_0}^{t_0} \frac{{\bf H}\left[w+\vp^{(w)}\right](t, x(t))}{a_1(t)-a_2(t)}dt,
\]
where $t\mapsto x(t;t_0,x_0)$ denotes the characteristics curve, obtained by solving $\dot x= a(t,x,w)$ with $x(t_0)=x_0$. This follows from slightly rewriting \eqref{wne}, which yields 
\begin{multline*}
\left(w^{(n+1)}+\vp^{n+1}\right)(t_0,x_0)~=~(\ov w+ \ov \vp)(x_n(\tau_0))+ \int_{\tau_0}^{t_0}\frac{ {\bf H}\left[w^{(n+1)}+\vp^{(n+1)}\right](t,x_n(t))}{a_1^{(n+1)}(t)-a_2^{(n+1)}(t)}dt\\
-\int_{\tau_0}^{t_0} \left(Z^{(n+1)}-\frac{W^{(n+1),-}(t)+W^{(n+1),+}(t)}{2}\right) \frac{\vp_x^{(n+1)}(t, x_n(t))}{a_1^{(n+1)}(t)-a_2^{(n+1)}(t)}dt\\
+ \int_{\tau_0}^{t_0} \left[{\vp^{(n+1)}_x(t,x_n(t))\over a_1^{(n+1)}(t)-a_2^{(n+1)}(t)}-{\vp^{(n+1)}_x(t,x_n(t))\over a_1^{(n)}(t)-a_2^{(n)}(t)}\right]\cdot\left(w^{(n)}-\frac{w^{(n),-}(t)+w^{(n),+}(t)}{2}+\vp^{(n)}\right)dt
\end{multline*}
where $t\mapsto x_n(t)$ denotes the characteristic curve, obtained by solving (\ref{char-s2}) with $x_n(t_0)=x_0$.

 Finally, to prove uniqueness, assume that $\Tilde{w}$ and $w$ are two entropic solutions.
We  define
\[
W\,\doteq\,\Tilde{w}-w,\qquad \beta(\tau)~\doteq~\sup_{t\in [0,\tau]} \Big[\|W(t,\cdot)\|_{H^1(\R\backslash\{t,0\})}+\left|W(t,0-)\right|+\left|W(t,0+)\right|\Big].
\]
The arguments used in the previous steps now yield the inequality
\begin{eqnarray*}
\beta(\tau)~\leq~ \beta(\tau)/2,
\end{eqnarray*}
and this implies $Z(\tau)=0$ for all $\tau\in [\tau_0,0]$, completing the proof.
\endproof
\medskip

{\bf Acknowledgments.} The research of A. Bressan was partially supported by NSF with grant DMS-2006884, ``Singularities and error bounds for hyperbolic equations". The research of K.T. Nguyen was  supported by a grant from the Simons Foundation/SFARI (521811, NTK). K. Grunert and S.T. Galtung were supported by the grant "Wave Phenomena and Stability - A Shocking Combination (WaPheS)" and K. Grunert was supported by the grant "Waves and Nonlinear Phenomena (WaNP)" both grants from the Research Council of Norway.

\end{document}